\documentclass[10.5pt]{amsart}
\usepackage[numbers]{natbib}
\usepackage[utf8]{inputenc}
\usepackage{amsthm}                             
\usepackage{ifthen}
\usepackage{hyperref}
\usepackage{cleveref}
\usepackage[shortlabels]{enumitem}                           
\usepackage[export]{adjustbox}
\DeclareMathAlphabet{\mathpzc}{OT1}{pzc}{m}{it}

\newtheorem{proposition}{Proposition}
\newtheorem{lemma}{Lemma}
\newtheorem{corollary}{Corollary}

\theoremstyle{remark}
\newtheorem{remark}{Remark}

\DeclareRobustCommand{\REALsymbol}{\mathbb{R}}
\DeclareRobustCommand{\REAL}[1][]{%
  \ifthenelse{\equal{#1}{}}{\REALsymbol}{{\REALsymbol}^{#1}}}
\DeclareRobustCommand{\NORM}[2][]{%
  \ifthenelse{\equal{#1}{}}{\left\| #2 \right\|}{\left\| #2
  \right\|_{#1}}}
\DeclareRobustCommand{\PolySymb}{\mathcal{P}}

\DeclareRobustCommand{\PONE}{\PolySymb_{1}}
     
\DeclareRobustCommand{\RTzero}{\mathcal{RT}_{0}}

\DeclareRobustCommand{\timeSymbol}{t}
\DeclareRobustCommand{\Time}{\timeSymbol}
\DeclareRobustCommand{\tempo}{\Time}
\DeclareRobustCommand{\Dt}{\Delta\tempo}        
\DeclareRobustCommand{\Der}[1][]{
  \ifthenelse{\equal{#1}{}}{\partial}{\partial_{#1}}}
\DeclareRobustCommand{\DerT}{\Der_{\tempo}\,}
\DeclareRobustCommand{\DerPar}[2]{\frac{\Der #1}{\Der #2}}

\DeclareRobustCommand{\Domain}{\Omega}                     
\DeclareRobustCommand{\SurfDomainsymb}{\Gamma}             
\DeclareRobustCommand{\SurfDomain}[1][]{%
  \ifthenelse{\equal{#1}{}}{\SurfDomainsymb}{\SurfDomainsymb_{#1}}}

\DeclareRobustCommand{\GradSymbol}{\operatorname{\mathbf{\nabla}}}
\DeclareRobustCommand{\Grad}{\GradSymbol}
\DeclareRobustCommand{\Div}{\Grad\cdot}
\DeclareRobustCommand{\Lap}{\Delta}

\DeclareRobustCommand{\GradSurf}{\GradSymbol_{\scriptscriptstyle\First}}
\DeclareRobustCommand{\DivSurf}{\GradSurf\cdot}
\DeclareRobustCommand{\LapSurf}{\Lap_{\scriptscriptstyle\First}}

\DeclareRobustCommand{\NormSymb}{\mathbf{N}}
\DeclareRobustCommand{\normalvec}{\NormSymb}               
\DeclareRobustCommand{\normalEdge}{\boldsymbol{\nu}}       
\DeclareRobustCommand{\tangentEdge}{\boldsymbol{\tau}}     
\DeclareRobustCommand{\conormalPlane}{\mathbf{n}}          

\DeclareRobustCommand{\basisGC}{e}                         
\DeclareRobustCommand{\vecBaseGC}[1][]{%
  \ifthenelse{\equal{#1}{}}{\mathbf{\basisGC}}{\mathbf{\basisGC}_{#1}}}
\DeclareRobustCommand{\basisCC}{t}                         
\DeclareRobustCommand{\vecBaseCCcv}[1][]{%
  \ifthenelse{\equal{#1}{}}{\mathbf{\basisCC}}{\mathbf{\basisCC}_{#1}}}
\DeclareRobustCommand{\hvecBaseCCcv}[1][]{
  \ifthenelse{\equal{#1}{}}{\hat{\mathbf{\basisCC}}}{\hat{\mathbf{\basisCC}}_{#1}}}

\DeclareRobustCommand{\ProjSymb}{\operatorname{\pi}}       
\DeclareRobustCommand{\AngleSymbol}{\theta}                

\DeclareRobustCommand{\From}{:}
\DeclareRobustCommand{\To}{\rightarrow}
\DeclareRobustCommand{\Mapsto}{\mapsto}
\DeclareRobustCommand{\MapUsymb}{\phi}                     
\DeclareRobustCommand{\MapU}[1][]{%
  \ifthenelse{\equal{#1}{}}{\MapUsymb}{\MapUsymb_{#1}}}
\DeclareRobustCommand{\SubsetU}[1][]{%
  \ifthenelse{\equal{#1}{}}{U}{U_{#1}}}
\DeclareRobustCommand{\SubsetV}[1][]{%
  \ifthenelse{\equal{#1}{}}{V}{V_{#1}}}

\DeclareRobustCommand{\point}[1][]{%
  \ifthenelse{\equal{#1}{}}{\mathbf{p}}{\mathbf{p}_{#1}}}
\DeclareRobustCommand{\Vector}[1]{\mathbf{#1}}
\DeclareRobustCommand{\fsymb}{f}
\DeclareRobustCommand{\scalFun}[1][]{%
  \ifthenelse{\equal{#1}{}}{\fsymb}{\fsymb_{#1}}}
\DeclareRobustCommand{\vecsymbu}{u}
\DeclareRobustCommand{\vecFunu}[1][]{%
  \ifthenelse{\equal{#1}{}}
  {\Vector{\vecsymbu}}
  {{\vecsymbu}^{#1}}
  }
\DeclareRobustCommand{\vecsymbv}{v}
\DeclareRobustCommand{\vecFunv}[1][]{%
  \ifthenelse{\equal{#1}{}}
  {\Vector{\vecsymbv}}
  {{\vecsymbv}^{#1}}
  }

\DeclareRobustCommand{\xvsymb}{x}
\DeclareRobustCommand{\xv}[1][]{%
  \ifthenelse{\equal{#1}{}}{\mathbf{\xvsymb}}{\mathbf{\xvsymb}_{\scriptscriptstyle{#1}}}}
\DeclareRobustCommand{\xvcomp}[1][]{%
  \ifthenelse{\equal{#1}{}}{\xvsymb}{\xvsymb_{\scriptscriptstyle{#1}}}}
\DeclareRobustCommand{\xcg}[1][]{%
  \ifthenelse{\equal{#1}{}}{\xvcomp[1]}{\xvcomp[1]_{\scriptscriptstyle{#1}}}}
\DeclareRobustCommand{\ycg}[1][]{%
  \ifthenelse{\equal{#1}{}}{\xvcomp[2]}{\xvcomp[2]_{\scriptscriptstyle{#1}}}}
\DeclareRobustCommand{\zcg}[1][]{%
  \ifthenelse{\equal{#1}{}}{\xvcomp[3]}{{\xvcomp[3]}_{\scriptscriptstyle{#1}}}}
\DeclareRobustCommand{\gcscomp}{\xcg,\ycg,\zcg}

\DeclareRobustCommand{\svsymb}{s}
\DeclareRobustCommand{\sv}[1][]{%
  \ifthenelse{\equal{#1}{}}{\mathbf{\svsymb}}{\mathbf{\svsymb}_{\scriptscriptstyle{#1}}}}
\DeclareRobustCommand{\svcomp}[1][]{%
  \ifthenelse{\equal{#1}{}}{\svsymb}{\svsymb_{\scriptscriptstyle{#1}}}}
\DeclareRobustCommand{\xcl}[1][]{%
  \ifthenelse{\equal{#1}{}}{\svcomp[1]}{\svcomp[1]_{\scriptscriptstyle{#1}}}}
\DeclareRobustCommand{\ycl}[1][]{%
  \ifthenelse{\equal{#1}{}}{\svcomp[2]}{\svcomp[2]_{\scriptscriptstyle{#1}}}}
\DeclareRobustCommand{\zcl}[1][]{%
  \ifthenelse{\equal{#1}{}}{\svcomp[3]}{\svcomp[3]_{\scriptscriptstyle{#1}}}}

\DeclareRobustCommand{\BSMsymbol}{\mathcal{B}}
\DeclareRobustCommand{\BSM}[1][]{%
  \ifthenelse{\equal{#1}{}}{\BSMsymbol}{\BSMsymbol_{#1}}}
\DeclareRobustCommand{\zbottom}{\zcg[\BSM]}

\DeclareRobustCommand{\Wmat}{\mathbf{W}}
\DeclareRobustCommand{\Jac}{\mathbf{J}}

\DeclareRobustCommand{\TanPlane}[2][]{%
  \ifthenelse{\equal{#1}{}}{T_{\scriptscriptstyle{\point}}#2}{T_{\scriptscriptstyle{#1}}#2}}

\DeclareRobustCommand{\FirstSymbol}{\mathcal{G}}
\DeclareRobustCommand{\First}[1][]{%
  \ifthenelse{\equal{#1}{}}{\FirstSymbol}{\FirstSymbol_{\scriptscriptstyle{#1}}}}
\DeclareRobustCommand{\metrcoefsymbol}{g}
\DeclareRobustCommand{\metrcoef}[1]{\metrcoefsymbol_{#1#1}}
\DeclareRobustCommand{\DET}[1]{\ensuremath{\operatorname{det}(#1)}}
\DeclareRobustCommand{\SCAL}[3][]{%
  \ifthenelse{\equal{#1}{}}{\left\langle{#2},{#3}\right\rangle}{\left\langle{#2},{#3}\right\rangle_{#1}}}
\DeclareRobustCommand{\scalprodSurf}[3][]{%
  \ifthenelse{\equal{#1}{}}{\left\langle {#2},{#3} \right\rangle_{\scriptscriptstyle\First}}{\left\langle {#2},{#3} \right\rangle_{{#1}}}}
\DeclareRobustCommand{\ABSsurf}[1]{\ensuremath{\left|#1\right|_{\First}}}

\DeclareRobustCommand{\meshparam}{h}           
\DeclareRobustCommand{\Tsymb}{\mathcal{T}}
\DeclareRobustCommand{\TriangH}{\Tsymb_{\meshparam}}
\DeclareRobustCommand{\Csymb}{\mathcal{K}}
\DeclareRobustCommand{\CellSet}{\Csymb_{\meshparam}}
\DeclareRobustCommand{\Esymb}{\mathcal{E}}
\DeclareRobustCommand{\EdgeSet}{\Esymb_{\meshparam}}
\DeclareRobustCommand{\Cellsymb}{K}
\DeclareRobustCommand{\Cell}[1][]{
  \ifthenelse{\equal{#1}{}}{\Cellsymb}{\Cellsymb_{#1}}
}
\DeclareRobustCommand{\Edgesymb}{E}
\DeclareRobustCommand{\Edge}[1][]{
  \ifthenelse{\equal{#1}{}}{\Edgesymb}{\Edgesymb_{#1}}
}
\DeclareRobustCommand{\lengthSymb}{\ell}
\DeclareRobustCommand{\edgeLength}[1][]{
  \ifthenelse{\equal{#1}{}}{\lengthSymb}{\lengthSymb_{#1}}
  } 
\DeclareRobustCommand{\StarSymb}{\omega}       
\DeclareRobustCommand{\Star}[1]{\StarSymb_{#1}}
\DeclareRobustCommand{\nodalBasis}[1]{\varphi_{\scriptscriptstyle{#1}}}

\DeclareRobustCommand{\Depth}{\eta}                        
\DeclareRobustCommand{\DepthP}[1][]{
  \ifthenelse{\equal{#1}{}}{\Depth_{\Domain}}{\Depth_{\Domain,#1}}}   
\DeclareRobustCommand{\kEta}{\kappa_{\Depth}}              
\DeclareRobustCommand{\kVol}{\kappa_{V}}                   
\DeclareRobustCommand{\RainVol}[1]{V_{#1}}                 
\DeclareRobustCommand{\tRec}{t_{50}^{\mathrm{rec}}}        
\DeclareRobustCommand{\Manning}{n}
\DeclareRobustCommand{\AdvVel}{\mathbf{w}}                 
\DeclareRobustCommand{\VelSymbol}{U}
\DeclareRobustCommand{\vectVel}{\boldsymbol{\VelSymbol}}   
\DeclareRobustCommand{\QSymbol}{q}
\DeclareRobustCommand{\Qvect}{\mathbf{\QSymbol}}           
\DeclareRobustCommand{\QvectP}{\Qvect_{\Domain}}           
\DeclareRobustCommand{\grav}{g}                            
\DeclareRobustCommand{\density}{\rho}                      
\DeclareRobustCommand{\FricSymbol}{f}
\DeclareRobustCommand{\vectFric}[1][]{
  \ifthenelse{\equal{#1}{}}{\mathbf{\FricSymbol}}{\mathbf{\FricSymbol}_{\scriptscriptstyle{#1}}}}
\DeclareRobustCommand{\force}{f}                           

\DeclareRobustCommand{\cellfield}{\boldsymbol{\psi}}       
\DeclareRobustCommand{\Transfer}{\mathcal{P}}              

\DeclareRobustCommand{\FVS}{FVS}                           
\DeclareRobustCommand{\AVG}{AVG}                           
\DeclareRobustCommand{\LNw}{LN\ensuremath{_{w}}}           
\DeclareRobustCommand{\LNs}{LN\ensuremath{_{s}}}           

\makeatletter
\DeclareRobustCommand{\pvset}[2]{\expandafter\def\csname pv@#1\endcsname{#2}}
\DeclareRobustCommand{\pv}[1]{\ifcsname pv@#1\endcsname\csname pv@#1\endcsname\else\textcolor{red}{\texttt{??\detokenize{#1}}}\fi}
\makeatother

\renewcommand{\citet}{\cite} 
\usepackage{booktabs,tabularx}
\pvset{t1.L}{100}
\pvset{t1.W}{20}
\pvset{t1.coincide}{$2.4\cdot10^{-14}$}
\pvset{t1.dtmax}{0.25}
\pvset{t1.exacterr.range}{0.042--0.049}
\pvset{t1.h}{2.5}
\pvset{t1.kind}{unstructured}
\pvset{t1.matched.max}{$2.1\cdot10^{-14}$}
\pvset{t1.modelgap.range}{0.6\%--36.3\%}
\pvset{t1.nc}{1000}
\pvset{t1.rain}{50}
\pvset{t1.slopes}{0.1, 0.31, 0.5, 0.85, 1}
\pvset{t1.tbarerr.range}{0.68\%--0.82\%}
\pvset{t1.tend}{30}
\pvset{t1.train}{20}
\pvset{t2.L}{100}
\pvset{t2.S}{0.6}
\pvset{t2.W}{40}
\pvset{t2.a}{$5.0\cdot10^{-3}$}
\pvset{t2.dtmax}{1}
\pvset{t2.ex.tbar.int}{140.4}
\pvset{t2.ex.tbar.proj}{114.5}
\pvset{t2.ex.tbar.ratio}{1.2259}
\pvset{t2.hlist}{5, 2.5, 1.25, 0.625}
\pvset{t2.hrun}{1.25}
\pvset{t2.proj.tria.area}{1.2290}
\pvset{t2.proj.tria.ratio}{1.2260}
\pvset{t2.proj.tria.storage}{1.2261}
\pvset{t2.rain}{50}
\pvset{t2.stru.avg.defect.coarse}{0.094}
\pvset{t2.stru.avg.err.axis}{0.010}
\pvset{t2.stru.avg.slopeerr.rate}{1.87}
\pvset{t2.stru.avgoverlnw}{0.28}
\pvset{t2.stru.cf}{160}
\pvset{t2.stru.fvs.err.axis}{0.094}
\pvset{t2.stru.fvs.slopeerr}{$0.0144$}
\pvset{t2.stru.lnw.err.axis}{0.010}
\pvset{t2.stru.lnw.err.off}{0.052}
\pvset{t2.stru.route.avgoverlnw}{0.23}
\pvset{t2.stru.run.tbar.range}{142.4--143.0}
\pvset{t2.stru.run.trise.spread}{1}
\pvset{t2.stru.transport}{$0.0262$}
\pvset{t2.stru.transport.rate}{0.81}
\pvset{t2.tbar.gap}{-18.4\%}
\pvset{t2.tbar.tria.errint.list}{5.6\%, 3.0\%, 1.6\%, 0.9\%}
\pvset{t2.tbar.tria.gapovererr}{21}
\pvset{t2.tend}{30}
\pvset{t2.train}{20}
\pvset{t2.tria.avg.defect.coarse}{0.28}
\pvset{t2.tria.avg.slopeerr}{$6.4\cdot10^{-4}$}
\pvset{t2.tria.avg.slopeerr.rate}{0.97}
\pvset{t2.tria.avgoverlnw}{5.35}
\pvset{t2.tria.fvs.slopeerr}{$2.8\cdot10^{-3}$}
\pvset{t2.tria.lns.slopeerr}{$1.4\cdot10^{-4}$}
\pvset{t2.tria.lnw.slopeerr}{$1.2\cdot10^{-4}$}
\pvset{t2.tria.lnw.slopeerr.rate}{1.03}
\pvset{t2.tria.route.avgoverlnw}{3.73}
\pvset{t2.tria.run.t10}{61}
\pvset{t2.tria.run.t50}{152}
\pvset{t2.tria.run.t90}{206}
\pvset{t2.tria.run.tbar.range}{142.5--142.7}
\pvset{t2.tria.run.trec}{94}
\pvset{t2.tria.transport}{$0.0242$}
\pvset{t2.tria.transport.rate}{0.77}
\pvset{t3.avg.vslnw}{0.052}
\pvset{t3.div.avg.defect}{0.325}
\pvset{t3.div.cf}{195}
\pvset{t3.div.jump}{0.023}
\pvset{t3.div.lnw.defect}{$1.1\cdot10^{-14}$}
\pvset{t3.dtmax}{1}
\pvset{t3.err.all.coarse.range}{0.37--0.40}
\pvset{t3.err.all.fine.range}{0.25--0.26}
\pvset{t3.exact.err.all.list}{0.44, 0.31}
\pvset{t3.exact.err.off.list}{0.088, 0.043}
\pvset{t3.exact.notrap.hlist}{0.5, 0.25}
\pvset{t3.exact.trap.hlist}{0.125}
\pvset{t3.exact.trap.rain}{41.2\%}
\pvset{t3.fam.pert.exact.trap.hlist}{0.25, 0.125}
\pvset{t3.fam.pert.exact.trap.rain}{41.1\%, 23.9\%}
\pvset{t3.fam.pert.fvs.trap.hlist}{0.25}
\pvset{t3.fam.pert.fvs.trap.rain}{41.3\%}
\pvset{t3.fvs.tmean}{52.3}
\pvset{t3.fvs.vslnw}{0.26}
\pvset{t3.hlist}{0.5, 0.25, 0.125}
\pvset{t3.hrun}{0.125}
\pvset{t3.inset.t0}{31}
\pvset{t3.inset.t1}{134}
\pvset{t3.lns.vslnw}{0.047}
\pvset{t3.lnw.err.off}{0.023}
\pvset{t3.lnw.tmean}{54.8}
\pvset{t3.mass.max}{$1.2\cdot10^{-13}$}
\pvset{t3.nclist}{2212, 8658, 34862}
\pvset{t3.nod.w.b.lnw}{$0.123$}
\pvset{t3.nod.w.c.D}{$0.145$}
\pvset{t3.nod.w.c.avg}{$0.0902$}
\pvset{t3.nod.w.c.lnw}{$0.0378$}
\pvset{t3.nod.w.d.D}{$3.6\cdot10^{-3}$}
\pvset{t3.nod.w.d.avg}{$3.6\cdot10^{-3}$}
\pvset{t3.nod.w.d.lnw}{$3.3\cdot10^{-3}$}
\pvset{t3.proj.agree}{$1.4\cdot10^{-4}$}
\pvset{t3.proj.area}{1.1035}
\pvset{t3.proj.h}{0.25}
\pvset{t3.proj.ratio}{1.0433}
\pvset{t3.proj.storage}{1.0432}
\pvset{t3.rain}{50}
\pvset{t3.refarea}{50\%}
\pvset{t3.refdisc}{0.05}
\pvset{t3.slopemax}{1.09}
\pvset{t3.slopemean}{0.35}
\pvset{t3.tend}{30}
\pvset{t3.thetamax}{47}
\pvset{t3.tmean.spread}{1.94\%}
\pvset{t3.train}{20}
\pvset{t4.area}{0.157}
\pvset{t4.avg.maxh}{0.659}
\pvset{t4.avg.tmean}{368}
\pvset{t4.avg.vslnw}{$0.0522$}
\pvset{t4.cosmin}{0.83}
\pvset{t4.div.avg.defect}{0.229}
\pvset{t4.div.avg.overstru}{2.4}
\pvset{t4.div.cf}{535}
\pvset{t4.div.cffrac}{3.0\%}
\pvset{t4.div.jump}{0.116}
\pvset{t4.dk}{$0.0418$}
\pvset{t4.dtmax}{10}
\pvset{t4.fvs.maxh}{333}
\pvset{t4.fvs.peak}{8.4\%}
\pvset{t4.fvs.trap.dz}{3.3}
\pvset{t4.fvs.trapcells}{3}
\pvset{t4.fvs.trapped}{91.6\%}
\pvset{t4.lns.tmean}{366}
\pvset{t4.lns.vslnw}{$0.0467$}
\pvset{t4.lnw.maxh}{0.521}
\pvset{t4.lnw.tmean}{377}
\pvset{t4.nc}{12560}
\pvset{t4.nv}{6468}
\pvset{t4.proj.agree}{$1.7\cdot10^{-3}$}
\pvset{t4.proj.area}{1.0701}
\pvset{t4.proj.ratio}{1.0496}
\pvset{t4.proj.storage}{1.0497}
\pvset{t4.rain}{50}
\pvset{t4.rings}{2}
\pvset{t4.slopemax}{0.67}
\pvset{t4.slopemean}{0.31}
\pvset{t4.tend}{120}
\pvset{t4.tmean.spread}{3.04\%}
\pvset{t4.train}{60}
\pvset{t4.tsteady}{0.39}
\pvset{t4c.fig.data}{depitted}
}
\graphicspath{{figures/}}

\usepackage{setspace}
\usepackage[a4paper,left=2.5cm,right=2.5cm, bottom=2.5cm]{geometry}
\title[IKW with conservative fluxes on curved terrain]{An intrinsic kinematic wave model with
  conservative fluxes for overland flow on curved terrain}

\author{Elena Bachini$^1$ and Matteo Camporese$^2$}
\address{$^1$ Department of Agronomy, Food, Natural
    resources, Animals and Environment, University of Padua, Italy
    (elena.bachini@unipd.it)\\
    $^2$ Department of Civil, Environmental, and Architectural
    Engineering, University of Padua, Italy
    (matteo.camporese@unipd.it)
}

\begin{document}
\maketitle

\begin{abstract}
  Kinematic wave (KW) models, the simplest and cheapest surface flow
  models, route overland flow by gravity alone, so on a terrain mesh
  the velocity computed from elevation data determines the flow paths,
  and there is no pressure or inertial term to redistribute flow when
  the resulting velocities misroute it. Computed from piecewise-linear
  elevation data, this velocity jumps across the edges of the mesh:
  two neighboring cells can push water into each other, and the scheme
  can hold water in ponds that the terrain does not contain, or
  concentrate it where the data do not. We formulate the KW model
  intrinsically on the bed surface, with the depth normal to the bed
  as unknown, and obtain a single-valued edge velocity, in the unique
  edge frame of the intrinsic geometry, by correcting the raw cell
  field with the star-based conservative flux of Larson and
  Niklasson. Because the field is derived from data rather than from a
  Galerkin solution, the local problems are solvable only with a
  compatible load, the weak discrete divergence of the field being
  corrected, which for the slope is a discrete mean curvature of the
  terrain. With a Neumann-compatible boundary treatment every element
  is conservative, and the corrected field reproduces, in every cell,
  the cell average of the nodal divergence of the raw field, plus the
  divergence resolved inside the curved cells. Since any single-valued
  flux conserves mass, the correction is about consistency and
  robustness: the steady state of the upwind scheme solves an M-matrix
  system exactly when every cell drains to an outlet, and spurious
  ponds are closed classes of a transfer directed graph, computable
  before any simulation. Four rainfall tests, from a plane to a real
  catchment, are designed to isolate one aspect each: the difference
  in storage and timing between the intrinsic and the classical KW,
  the effect of the mesh on the consistency of the edge velocity, the
  robustness at a degenerate stationary point, and the application to
  terrain data. In the convergent-valley test, the projected model
  underestimates the intrinsic mean response time by 18.4\%, while
  conservative reconstruction reduces the slope-induced error of plain
  edge averaging by a factor of 5.35 on the finest unstructured
  mesh. On the real catchment, raw-field flux-vector splitting
  produces a numerical trap intercepting 91.6\% of the rainfall,
  whereas the tested single-valued edge-flux methods avoid this
  trap. These results quantify the sensitivity of kinematic-wave
  predictions to both geometric formulation and numerical flux
  construction, without establishing the physical superiority of
  either model.
\end{abstract}


\subjclass{
}

\section{Introduction}\label{sec:intro}
Surface runoff is modeled in practice using mainly three approaches
belonging to the same family: the Shallow Water (SW) equations, the
Diffusive Wave (DW) approximation, and the Kinematic Wave (KW)
approximation. The SW equations are a consistent two-dimensional
depth-integrated approximation of the Navier-Stokes equations, with
the DW approach being derived from SW by dropping the inertial
terms. A further approximation, at the basis of the KW model, neglects
the pressure terms so that the flow is driven by gravity and with a
direction that follows the steepest descent.  These three models are
used by different communities and their domains of validity are
neither always clear nor analyzed with the appropriate mathematical
rigor.  Classical applicability criteria exist in one dimension on
uniform slopes \citep{WoolhiserLiggett1967,Ponce1978}, while for
two-dimensional flows over curved terrain the question is largely
open.  Thin-layer models in topography-fitted coordinates, with the
depth measured along the normal to the bed, are classical for granular
and debris flows \citep{SavageHutter1989,DenlingerIverson2004}, and
were derived for shallow water over variable topography by
\citep{BouchutWestdickenberg2004, FentEtAl2018, BachiniPutti2020}.  A
geometrically intrinsic formulation of the shallow water model was
presented by \cite{BachiniPutti2020}. Here the bed geometry enters the
model through a local curvilinear coordinate system anchored on the
bed.  Further approximations of the Intrinsic Shallow Water
Equations (ISWE) can produce corresponding intrinsic versions of the
diffusive wave and kinematic wave equations.  Our long term goal is to
analyze these models, the importance of the geometry, and their
validity in different scenarios, and the current work is a first step
in this direction.

The KW approach is the natural entry point and the focus for this
work. It is the simplest member of the hierarchy and the one most
sensitive to geometry. Without inertia and pressure, the terrain slope
is its only driver, so the discrete velocity field arising from the
terrain gradient determines the flow paths, and no term of the model
can redirect a flow that these velocities misroute.  This is the
central idea here: on a triangulated terrain, the consistency of the
discrete model with the elevation data and its robustness are largely
decided, before any water is moved, by the way a single-valued edge
velocity is constructed from the data.  We start by deriving the model
on the bed surface itself as the zero-inertia and zero-pressure limit
of the ISWE, and we discretize it with a first-order finite volume
(FV) scheme built on the intrinsic geometry \citep{BachiniPutti2020}.
Despite its simplicity, the KW model and its discretization have
requirements that are rarely stated explicitly.  The first concerns
the terrain: the topography must be \emph{KW-admissible} (i).  Without
a pressure term, both a closed depression (pit) and a flat region trap
water forever, and a boundary segment towards which the flow is
directed must be an outflow.  Real Digital Elevation Models (DEMs)
violate these requirements, and the usual remedies, obtained by
filling \citep{JensonDomingue1988,PlanchonDarboux2002} or breaching
\citep{Lindsay2016} of the depressions, modify the elevations and can
introduce new artifacts, such as flats and parallel flow paths.  The
other two requirements concern the numerical scheme.  An FV-based
transport solver needs a \emph{single-valued} and \emph{locally
  conservative} normal velocity on every edge (ii), a condition
violated by a velocity field computed directly from the surface mesh.
Indeed, the tangential gradient of a piecewise-linear approximation of
the height function describing the topographic surface, and hence the
kinematic field built from it, is piecewise constant and jumps across
edges and is generally not locally conservative.  On a curved surface,
single-valuedness also requires a common reference frame on each edge
(iii), since the two cells sharing an edge have different tangent
planes.
For this study, we assume requirement (i) to be satisfied by the
input.  Requirement (iii) is met by the unique edge frame of the
intrinsic FV scheme \citep{BachiniPutti2020}, while requirement (ii)
is one of the main subjects of this work.
We obtain the single-valued edge velocity feeding the finite volume
scheme by a locally conservative correction of the piecewise-constant
kinematic field, and we analyze what this construction imposes on the
velocity field and what it implies for the robustness of the scheme.

This correction is the discrete counterpart of a situation well known
in porous-media flow, where a continuous Galerkin (CG) velocity must
be made compatible with a transport scheme.  Local post-processing
techniques recover, from the CG solution, edge fluxes that are
single-valued and conservative on every element
\citep{LarsonNiklasson2004,SunWheeler2006,Cockburn2007,Odsaeter2017,DengGinting2018}. A
global alternative, a constrained least-squares projection of the
$\PONE$ gradients onto the $\RTzero$ mixed finite element space with
element-wise conservation, was proposed for particle tracking in
\citep{SelzerCirpka2020}.  In hydrological practice, overland flow on
terrain meshes often bypasses the velocity field altogether.
Flow-routing algorithms on grids and triangulations assign transfer
fractions between neighboring cells directly from the elevations
\citep{Quinn1991,Freeman1991,Tarboton1997,Tucker2001,OrlandiniMoretti2009}.
More recent approaches have focused on the identification of surface
flow paths and the preservation of topographic features controlling
drainage \citep{OrlandiniMoretti2009,MorettiOrlandini2023}.  In
particular, in \citep{MorettiOrlandini2023} the authors proposed the
extraction of thalweg and ridge networks from unaltered elevation data
as a basis for terrain partitioning. Such approaches address the
representation of DEM-consistent terrain features and drainage
pathways, but do not provide a consistent velocity field for a
PDE-based surface routing model.
On surfaces, locally conservative discretizations exist, but
conservation is obtained by changing the discretization rather than by
post-processing a CG field: hybridizable discontinuous Galerkin (DG)
and mixed methods, including surface Raviart-Thomas spaces
\citep{CockburnDemlow2016}, interior-penalty DG
\citep{DednerMadhavanStinner2013}, cut DG for the Laplace-Beltrami
operator \citep{BurmanCutDG2017}, conservative discontinuous cut
finite element method (FEM) for coupled bulk-interface problems, with
conservation on macro elements \citep{LarsonZahedi2023}, and finite
volumes \citep{JuDu2009,GiesselmannMuller2014}.
The absence of flux recovery on surfaces is also visible in a
posteriori error estimation, which for surface FEM is, as far as we
could find, residual-based
\citep{DemlowDziuk2007,BonitoDemlowNochetto2020}. An equilibrated
(Prager-Synge type) estimator, which would require precisely a
conservative flux recovery, has not been reported, whereas on flat
domains such estimators are built from vertex-patch problems with a
compatible load \citep{BraessPillweinSchoeberl2009,ErnVohralik2015}.

Among these techniques we adopt the local, star-based conservative
flux of Larson and Niklasson (LN) \citep{LarsonNiklasson2004}.  It
corrects the average of the two one-sided fluxes on each edge with the
solution of a small problem on every vertex star, and it only requires
the field to be piecewise constant on the cells, not to come from a
particular solver.  On a surface it has a further property that makes
it the natural choice: it works with one scalar per edge and never
compares vectors across edges, so the geometry enters only through
per-cell and per-edge scalars, and a common frame on each edge,
requirement (iii), is all it needs from the geometry.
In subsurface hydrology, the LN post-processing has been applied to
variably saturated flow \citep{Kees2008}, and in further hydrological
applications \citep{Cheng2010,Kees2011,Povich2013}. For
surface-subsurface applications it reconstructs mass-conservative
velocities from $\PONE$ Galerkin solutions of Richards equation to
drive a finite volume advective transport model \citep{Scudeler2016},
where locally conservative velocities were found necessary for
accurate solute transport.  In \citep{Camporese2010}, the surface
routing is path-based, and LN is applied to the subsurface solute
transport only.  In all these applications, as in
\citet{LarsonNiklasson2004}, the corrected field is the gradient of
the Galerkin solution of a PDE with a given source, and Galerkin
orthogonality makes every star problem solvable. Here the field is
computed from the terrain height, which is \emph{data}: there is no
given source term, and the compatibility of the star problems does not
hold automatically.

The combination of LN and FV is therefore not new. However, to our
knowledge, the LN construction has not been used for terrain-driven
overland flow, and star-based conservative post-processing of surface
CG fields with a compatible divergence load has not been analyzed yet.

The three main contributions of this work are the following:
\begin{enumerate}[label=\arabic*.]
\item \emph{Intrinsic kinematic wave equation.}  We derive the
  kinematic wave approximation of the intrinsic shallow water
  equations~\citep{BachiniPutti2020} and investigate its relationship
  with the classical projected formulation, with particular attention
  to storage and transient response.
\item \emph{Conservative edge velocities from data on a curved
    surface.}  We extend the Larson-Niklasson flux reconstruction to
  terrain-derived gradient fields on curved surfaces, addressing the
  compatibility of the local problems, the role of surface geometry,
  and the reconstruction of locally conservative single-valued edge
  velocities.
\item \emph{Discrete traps as a graph property.} We investigate the
  relationship between the discrete velocity field and the existence
  of artificial flow traps, using the transfer directed graph
  associated with the finite-volume discretization to characterize the
  conditions for a well-defined steady state.
\end{enumerate}

Four rainfall tests isolate one effect each: a plane with the exact
event hydrograph, a steep convergent valley with a closed-form steady
state on structured and unstructured meshes, a fully three-dimensional
surface with a degenerate stationary point, and a real catchment.  The
same numerical scheme is also applied to the \emph{projected} model,
i.e., the classical KW equation where the depth is measured along the
vertical, to investigate the influence of the geometric formulation on
storage and transient response.  The comparison also builds on the
previous studies of intrinsic and projected shallow water formulations
presented in \citep{FentEtAl2018,BachiniPutti2020}.

The paper is organized as follows.  Section~\ref{sec:model} states the
intrinsic and the projected KW approximations, their exact relations
and the requirements of the kinematic wave model.
Section~\ref{sec:geometry} describes the discrete surface geometry and
what the method requires of a terrain input.  Section~\ref{sec:fv}
introduces the intrinsic FV scheme, its exact steady state and the
characterization of discrete traps.  Section~\ref{sec:ln} presents the
conservative correction of the slope and of the kinematic field, its
compatible load and its boundary treatment.
Section~\ref{sec:experiments} reports our numerical experiments, and
finally, Section~\ref{sec:conclusions} draws the conclusions and
outlines future work.

\section{The kinematic wave on the bed surface}\label{sec:model}

\subsection{Intrinsic kinematic wave}\label{sec:intrinsic}

Consider the bed surface to be a regular surface
$\SurfDomain\subset\REAL[3]$.  We assume it to be covered by local
charts where local Monge parametrizations can be employed.  On each
chart the local Monge parametrization reads:
\begin{equation*}
  \MapU\From\SubsetU\subset\REAL[2]
  \To\SubsetV\cap\SurfDomain\subset\REAL[3],
  \qquad
  (\xcg,\ycg)\Mapsto(\xcg,\ycg,\BSM(\xcg,\ycg))=\xv ,
\end{equation*}
where $\xv=(\gcscomp)$ are the global Cartesian coordinates and $\BSM$
is the bed height over the chart.  A terrain surface is a graph over
its vertical projection, so in all the cases considered here a single
chart suffices. We will call $\Domain$ the projection of $\SurfDomain$
on $\SubsetU$.
Following \citep{BachiniPutti2020}, we define a local curvilinear
reference system (LCS) on $\SurfDomain$ and carry out all computations
with respect to it.  Since we consider only surface processes, all
vectors are tangent vectors, and the setting reduces to the tangent
plane
$\TanPlane{\SurfDomain}=\langle\vecBaseCCcv[1],\vecBaseCCcv[2]\rangle$.
The natural tangent vectors of the parametrization are:
\begin{equation*}
  \hvecBaseCCcv[1]=\left(1,0,\Der_x\BSM\right),\qquad
  \hvecBaseCCcv[2]=\left(0,1,\Der_y\BSM\right),
\end{equation*}
and, to ensure numerical stability, we decide to work with orthogonal
frames.  Thus, we set $\vecBaseCCcv[1]=\hvecBaseCCcv[1]$ and obtain
$\vecBaseCCcv[2]$ by Gram-Schmidt orthogonalization of
$\hvecBaseCCcv[2]$ with respect to $\vecBaseCCcv[1]$.  The local
coordinates associated with this orthogonal frame are indicated with
$\sv=(\xcl,\ycl)$. We let $\Wmat=\Jac^{+}\hat{\Jac}$, with
$\Jac=[\vecBaseCCcv[1],\vecBaseCCcv[2]]$ and
$\hat{\Jac}=[\hvecBaseCCcv[1],\hvecBaseCCcv[2]]$ the Jacobians of the
changes between the embedding coordinates $\xv$ and the local
coordinates $\sv$.
The $\REAL[3]$ frame in the LCS is completed by the unit normal
$\vecBaseCCcv[3]=\normalvec$, oriented upward so that the right-hand
rule holds, with associated coordinate $\zcl$.
A tangent vector is always written in contravariant components,
$\vecFunu=\vecFunu[1]\vecBaseCCcv[1]+\vecFunu[2]\vecBaseCCcv[2]$.  By
orthogonality, the metric tensor $\First$ is diagonal and given by:
\begin{equation}\label{eq:first-FF}
  \First=
  \begin{pmatrix}
    \NORM{\vecBaseCCcv[1]}^2 & 0\\ 0 & \NORM{\vecBaseCCcv[2]}^2
  \end{pmatrix}
  =
  \begin{pmatrix}
    \metrcoef{1} & 0\\ 0 & \metrcoef{2}
  \end{pmatrix},
\end{equation}
with scalar product
$\scalprodSurf{\vecFunu}{\vecFunv}=\SCAL{\vecFunu}{\First \vecFunv}$
and norm $\ABSsurf{\cdot}$.  The intrinsic gradient and divergence are
\begin{equation}\label{eq:operators}
  \GradSurf\scalFun=\First^{-1}\Grad_{\sv}\scalFun
  =\First^{-1}\Wmat\Grad_{\xv}\scalFun,\qquad
  \DivSurf\vecFunu=\frac{1}{\sqrt{\DET{\First}}}\,
  \Grad_{\sv}\cdot\big(\sqrt{\DET{\First}}\,\vecFunu\big),
\end{equation}
and the Laplace-Beltrami operator is $\LapSurf=\DivSurf\GradSurf$,
where the derivatives of the metric coefficients in $\DivSurf$ carry
the Christoffel symbols. The operators are intrinsic, i.e.,
independent of the chosen frame: the orthogonal frame is the one used
in the computations and in the discretization, while closed forms on a
graph surface are most easily written in the natural coordinates, with
(full) metric $\hat{\First}$, where
$\hat g_{ij}=\delta_{ij}+\partial_i\BSM\,\partial_j\BSM$ and
$\DET{\hat \First}=1+|\Grad\BSM|^2$.

Let $\zbottom$ denote the height of the bed surface,
$\zbottom(\xv)=\BSM(\xcg,\ycg)$, and let $\AngleSymbol$ be the local
inclination of the bed, i.e.\ the angle between the vertical and the
normal direction:
\begin{equation}\label{eq:cos}
  \cos\AngleSymbol=\DerPar{\zcg}{\zcl}=\normalvec\cdot\vecBaseGC[3]
  =\frac{1}{\sqrt{1+|\Grad\BSM|^2}} .
\end{equation}
Two identities are used throughout. Seen as a vector of $\REAL[3]$,
the intrinsic gradient of the height is the tangential part of the
vertical unit vector,
\begin{equation}\label{eq:gradG}
  \GradSurf\zbottom=\vecBaseGC[3]-\cos\AngleSymbol\,\normalvec,\qquad
  \ABSsurf{\GradSurf\zbottom}=\sin\AngleSymbol
  =\frac{|\Grad\BSM|}{\sqrt{1+|\Grad\BSM|^2}},
\end{equation}
whose horizontal projection is $\Grad\BSM/(1+|\Grad\BSM|^2)$, parallel
to $\Grad\BSM$. In the natural coordinates, a direct computation gives:
\begin{equation}\label{eq:LB}
  \LapSurf\zbottom=\frac{1}{\sqrt{\DET{\hat \First}}}\,
  \partial_i\big(\sqrt{\DET{\hat \First}}\,\hat g^{ij}\partial_j\BSM\big)
  =\cos\AngleSymbol\,H,\qquad
  H=\Div\Big(\frac{\Grad\BSM}{\sqrt{1+|\Grad\BSM|^2}}\Big),
\end{equation}
where $H$ is the mean curvature (sum of the principal curvatures,
positive in valleys with this sign convention).  For small slopes
$\cos\AngleSymbol\,H\approx\Lap\BSM$.  Equivalently,
$\LapSurf\xv=H\normalvec$ for the position vector, of which
\eqref{eq:LB} is the vertical component.

In this geometric framework, the ISWE model derived in \cite{BachiniPutti2020} reads
\begin{subequations}\label{eq:iswe}
  \begin{align}
    &\DerPar{\Depth}{\tempo}+\DivSurf\Qvect=\force, \label{eq:iswe-cont}\\
    &\DerPar{\Qvect}{\tempo}
      +\DivSurf\left(\frac{1}{\Depth}\,\Qvect\otimes\Qvect
      +\frac{\grav\Depth^2}{2}\cos\AngleSymbol\,\First^{-1}\right)
      +\frac{\grav\Depth^2}{2}\GradSurf\cos\AngleSymbol
      +\grav\Depth\GradSurf\zbottom-\frac{\vectFric_{_{\BSM}}}{\density}=0,
      \label{eq:iswe-mom}
  \end{align}
\end{subequations}
where the unknown is the depth $\Depth$ of the water layer measured
along the normal to the bed, $\vectVel$ is the depth-averaged tangent
velocity, and $\Qvect=\Depth\vectVel$ [L$^2$/T] the discharge per unit
width.  The symbol $\force$ denotes a source term in the volume
equation, and the bottom friction is given by Manning law,
\begin{equation}\label{eq:manning}
  \frac{\vectFric_{_{\BSM}}}{\density}
  =-\frac{\grav\Manning^2\vectVel\ABSsurf{\vectVel}}{\Depth^{1/3}},
\end{equation}
with $\Manning$ the Manning coefficient [T/L$^{1/3}$]. 

The intrinsic kinematic wave is obtained by retaining
in~\eqref{eq:iswe-mom} only the balance between gravity and friction:
the time derivative and the advective flux (inertia) are neglected,
together with the two terms proportional to $\grav\Depth^2/2$, the
hydrostatic pressure and its curvature correction. This leaves
\begin{equation*}
  \grav\Depth\GradSurf\zbottom-\frac{\vectFric_{_{\BSM}}}{\density}=0
  \quad\Longrightarrow\quad
  \vectVel\ABSsurf{\vectVel}
  =-\frac{\Depth^{4/3}}{\Manning^2}\,\GradSurf\zbottom ,
\end{equation*}
for $\Depth>0$. Taking the norm of both sides gives
$\ABSsurf{\vectVel}^2=\Depth^{4/3}\ABSsurf{\GradSurf\zbottom}/\Manning^2$,
hence
\begin{equation}\label{eq:kwvel}
  \vectVel
  =-\frac{\Depth^{2/3}}{\Manning\,\ABSsurf{\GradSurf\zbottom}^{1/2}}
  \,
  \GradSurf\zbottom ,
\end{equation}
and substituting into \eqref{eq:iswe-cont} gives the \emph{Intrinsic
  Kinematic Wave (IKW)} equation
\begin{equation}\label{eq:kwint}
  \DerPar{\Depth}{\tempo}+\DivSurf\Qvect=\force,\qquad \mbox{with }
  \Qvect=-\frac{\Depth^{5/3}}{\Manning}\,
  \frac{\GradSurf\zbottom}{\ABSsurf{\GradSurf\zbottom}^{1/2}} .
\end{equation}
By~\eqref{eq:gradG} the friction slope is $\sin\AngleSymbol$ and the
flux is directed along the steepest descent on $\SurfDomain$.  We
rewrite
\begin{equation}\label{eq:wfield}
  \Qvect=\Depth^{5/3}\AdvVel,\qquad
  \AdvVel=-\frac{\GradSurf\zbottom}{\Manning\,\ABSsurf{\GradSurf\zbottom}^{1/2}},
  \qquad \vectVel=\Depth^{2/3}\AdvVel ,
\end{equation}
where we define $\AdvVel$ as the \emph{Manning-weighted kinematic
  field}, which is the velocity per unit $\Depth^{2/3}$ [L$^{1/3}$/T]
depending on the geometry only. In other words, it is the vector field
with which the scheme transports water, and one of the quantities on
which the conservative correction of Section~\ref{sec:ln} is applied,
the other being the slope. From now on, for simplicity, we will refer
to $\AdvVel$ as kinematic field.

We consider the volumetric source to be rainfall. Its intensity $r$ is
measured as a volume per unit time and unit horizontal area (vertical
rain), and its conversion into a depth rate depends only on the
measure of integration: per unit bed area it is
$\force=r\cos\AngleSymbol$, and the volume received by a region
$K\subset\SurfDomain$ is
\begin{equation}\label{eq:rainvol}
  \int_K r\cos\AngleSymbol\,dA_{\SurfDomain}=\int_{\ProjSymb(K)}r\,dA ,
\end{equation}
with $\ProjSymb$ the vertical projection (area formula). 

\subsection{Intrinsic versus projected model}\label{sec:equiv}

The classical KW model used in hydrological practice is written on the
horizontal projection $\Domain$, with the vertical depth $\DepthP$ as
unknown:
\begin{equation}\label{eq:kw}
  \DerPar{\DepthP}{\tempo}+\Div\QvectP=\force=r,\qquad  \mbox{with }
  \QvectP=-\frac{\DepthP^{5/3}}{\Manning}\,\frac{\Grad\BSM}{|\Grad\BSM|^{1/2}} ,
\end{equation}
so that the friction slope is
$\tan\AngleSymbol=|\Grad\BSM|$. Formally, \eqref{eq:kw} is
\eqref{eq:kwint} written on the flat domain $\Domain$, with $\First$
the identity, $\cos\AngleSymbol$ replaced by $1$ and the bed height
$\BSM$ acting as an external potential.  On a sloping bed, the two
formulations differ in their definitions of water depth, friction
slope, and rainfall-receiving area, even at their small-slope limit.
The discretization of Sections~\ref{sec:geometry}-\ref{sec:ln} is
written once for both cases, and the projected model serves as the
model to compare with, not as a limit to be verified.

\subsubsection*{Steady-state and hydrograph response}
Let us assume a dry bed at time $t=0$ and let $r$ be constant in time.
We consider an admissible terrain where the downslope flow paths are
well defined and the boundary conditions allow all water to drain out
of the domain.  At steady state, the unit discharge (volume per unit
time and unit contour length) of the two models coincides:
$\ABSsurf{\Qvect}(\xv)=|\QvectP|(\ProjSymb(\xv))$ at every regular
point.  Indeed, the two models have the same projected flow directions
and receive the same rainfall volume per unit horizontal area.
Moreover, elevation contour lines are horizontal, so their lengths are
identical in the intrinsic and projected geometries.
Consequently, it is immediate to show that the steady depths satisfy, pointwise,
\begin{equation}\label{eq:depthratio}
 \frac{\Depth}{\DepthP}=\cos^{-3/10}\AngleSymbol=: \kEta,
\end{equation}
and the stored volume per unit horizontal area satisfies
\begin{equation}\label{eq:volration}
  \frac{\Depth/\cos\AngleSymbol}{\DepthP}=\cos^{-13/10}\AngleSymbol=:\kVol.
\end{equation}
We call $\kEta$ and $\kVol$ the \emph{depth factor} and the \emph{storage factor}, respectively.

If now we consider to let every point drain to the outlet, we obtain
that the centroid of the rising limb,
\begin{equation*}
  \bar t=\int_0^\infty\Big(1-\frac{Q(t)}{Q_s}\Big)\,dt=\frac{V_s}{Q_s},
\end{equation*}
is the steady stored volume over the steady outflow (mass balance
$\tfrac{d}{dt}V=Q_s-Q$). Since $Q_s=r|\ProjSymb(\SurfDomain)|$ in both
models,
\begin{equation*}
  \frac{\bar t_{\SurfDomain}}{\bar t_\Domain}
  =\frac{V_{s,\SurfDomain}}{V_{s,\Domain}}
  =\frac{\int_\Domain \DepthP[s]\,\kVol\,d\xv}
  {\int_\Domain\DepthP[s]\,d\xv},
\end{equation*}
the mean of the local storage factor $\kVol$ of \eqref{eq:volration}
weighted by the projected steady storage $\DepthP[s]$, not by area.
Storage concentrates where the slope is gentle (valley floors), so the
ratio lies below the area mean of $\kVol$.
The identity $\bar t=V_s/Q_s$ holds for any conservative scheme
started from a dry bed, so $\bar t$ is a time metric that does not
depend on the sampling of the hydrograph.  Two consequences frame the
comparison between the models: steady state tests only the depth
closure, i.e.\ flow directions and contributing areas, while the
geometric effects on routing appear in transients, where the local
time-scale factor varies in space.

\begin{table}\centering \small
  \begin{tabular}{@{}lcccccc@{}}
    \toprule
    slope $\tan\AngleSymbol$ & 0.05 & 0.1 & 0.31 & 0.5 & 0.85 & 1.0 \\ \midrule
    depth, $\kEta=\cos^{-3/10}\AngleSymbol$ & 1.0004 & 1.0015 & 1.014 & 1.034 & 1.085 & 1.110 \\
    storage and $t_e$, $\kVol=\cos^{-13/10}\AngleSymbol$ & 1.0016 & 1.0065 & 1.062 & 1.156 & 1.424 & 1.569 \\
    \bottomrule
  \end{tabular}
  \caption{Values of $\kEta$ and $\kVol$ on a rectangular sloping
    plane at different slopes ($\tan\AngleSymbol$)}\label{tab:factors}
\end{table}

\begin{remark}[Time scales]\label{rem:timescale}
  On a plane of inclination $\AngleSymbol$ and horizontal length $L$,
  per unit contour width, the two models can be recasted as
  \begin{equation*}
    \DerT\Depth+\alpha_{\SurfDomain}\,\Der[\xi](\Depth^{5/3})
    =r\cos\AngleSymbol,\ \ \xi\in[0,L/\cos\AngleSymbol],
    \qquad
    \DerT\DepthP+\alpha_\Domain\,\Der[x](\DepthP^{5/3})=r,\ \ x\in[0,L],
  \end{equation*}
  with $\alpha_{\SurfDomain}=\sin^{1/2}\AngleSymbol/\Manning$ and
  $\alpha_\Domain=\tan^{1/2}\AngleSymbol/\Manning$.  The substitution
  $\Depth(\xi,t)=\kEta\,\DepthP(\xi\cos\AngleSymbol,\,t/\kVol)$ turns
  the intrinsic equation into $\cos\AngleSymbol$ times the projected
  one, since $\kEta/\kVol=\cos\AngleSymbol$ and
  $\kEta^{5/3}\sin^{1/2}\AngleSymbol=\tan^{1/2}\AngleSymbol$, and the
  outlet discharges are equal.
  Hence, the intrinsic hydrograph is the projected one stretched in
  time, $Q_{\SurfDomain}(t)=Q_\Domain(t/\kVol)$, and
  $t_{e,\SurfDomain}/t_{e,\Domain}=\kVol$, provided the rain duration
  is stretched too.
  The factors are summarized in Table~\ref{tab:factors}.
  Below slopes of about $0.1$ all intrinsic effects are below $1\%$,
  while at slopes typical of hillslopes ($0.3$-$0.85$) the storage
  differs by $6$-$42\%$.
  Qualitatively similar differences between intrinsic and projected
  formulations have been reported for the shallow water equations
  \citep{BachiniPutti2020,FentEtAl2018}.
\end{remark}

\subsection{Model requirements}
\label{sec:admissibility}\label{sec:requirements}

The intrinsic KW rests on two approximations, and each restricts model
applicability.
The first approximation comes from the shallow water hypothesis.  The
ISWE are derived in \citep{BachiniPutti2020} for a layer of fluid
whose depth is small compared with the horizontal length scales and
for a bed with \emph{small curvatures}, resulting in an approximation
parameter $\epsilon_{\First}\ll1$ that is geometry dependent.  No
assumption is made on the slope, since all the quantities are
tangential to the bed surface.  This is why the depth normal to the
bed is the natural unknown as it is the quantity in which the
thin-layer assumption is stated and remains small on steep beds. On
the other hand, the vertical depth of the same layer,
$\Depth/\cos\AngleSymbol$ on a plane, grows with the slope and the
pressure and the hydrostatic assumption is no longer accurate along
the vertical line.  In fact, classical SWE models require, in
addition, a mildly sloping bed and the projected KW \eqref{eq:kw}
inherits this restriction.

The second level of approximation in the IKW model comes from
neglecting the inertia and pressure gradient terms with respect to
gravity and friction. In one dimension and uniform slopes, classical
applicability criteria are based on the kinematic number
$k=S_0L/(H_0F_0^2)$ \citep{WoolhiserLiggett1967}, where $S_0$ is the
bed slope, $L$ is the length of the plane, and $H_0$ and $F_0$ are the
normal depth and Froude number at the outlet, respectively.  The
kinematic approximation is generally considered applicable for
sufficiently large $k$ (e.g., $k\geq20$ \citep{WoolhiserLiggett1967}),
with an additional criterion $kF_0^2\geq5$ for low-Froude-number flows
\citep{MorrisWoolhiser1980}.  Alternative criteria involve a
dimensionless wave period \citep{Ponce1978}.

Thus, the applicability of the kinematic approximation depends on the
relative magnitude of the neglected inertia and pressure terms, rather
than on the bed slope alone.  At the same time, the geometric
differences between the intrinsic and projected formulations increase
with slope, as quantified by the factors in
Remark~\ref{rem:timescale}.  However, these differences do not
establish which formulation provides a more accurate description of
the actual physical process.  For two-dimensional flows over curved
terrain, no comparable applicability criterion has been established,
and we do not address this question here.
The absence of a pressure term adds requirements on the terrain, since
without this term topographic gradient completely determines the flow
directions, which are the same both for IKW in \eqref{eq:kwint} and
the (projected) KW in \eqref{eq:kw}, since the stationary points of
$\zbottom$ on $\SurfDomain$ are those of $\BSM$.

We conclude this section with a summary of the admissibility
conditions on the terrain, while their discrete counterparts, i.e.\
the admissibility of a terrain mesh, are discussed in
Section~\ref{sec:nodalgeom}.
\begin{itemize}
\item \emph{No pits.} Every point other than the outlet must have a
  strictly descending direction, otherwise water accumulates without
  bound.  This is a limitation of the KW model due to the absence of
  the pressure terms, which are essential in filling depressions. As a
  result, the DW and ISWE models do not have this restriction.
\item \emph{No flats.} Where $\Grad\BSM=0$ on a set of positive
  measure, the flow direction is undefined.
  \emph{Isolated stationary points}, i.e., $\Grad\BSM=0$ at a point,
  are admissible if non-degenerate (peaks, saddles).  In fact, in this
  case, only a set of flow lines of measure zero reaches them, they
  carry no discharge, and the steady depth satisfies
  $\Depth\propto q^{3/5}|\Grad\BSM|^{-3/10}$.  The degenerate ones,
  where the gradient vanishes to higher order, are not harmless.
  Since flow lines are not unique there, the steady kinematic wave
  does not determine how this discharge is redistributed downstream
  and the problem is ill posed downstream of the point.  Thus, the
  discrete split is mesh dependent.
\item \emph{Boundary rule.} Every boundary segment must be either a
  watershed divide ($-\Grad\BSM\cdot\conormalPlane\le0$, flow tangent
  or inward) or an outflow.  Free outflow on the whole boundary is the
  consistent default, since on true divides the flux is inward and
  remains zero.  A no-flux wall towards which the flow points traps
  water, and thus is not admissible.
\end{itemize}
%

\section{Discrete surface geometry}\label{sec:geometry}

This section describes how the bed surface is discretized, which
geometric quantities the method uses, and how they are obtained from
elevation data.  The construction follows \citet{BachiniPutti2020},
the main addition is the treatment of real terrain data.

\subsection{Surface triangulation}\label{sec:triangulation}

A triangulation $\TriangH=\{\CellSet,\EdgeSet\}$ of the bed surface
$\SurfDomain$ with mesh size $\meshparam$ is given by a set $\CellSet$
of cells $K$, with vertices $\xv_i\in\SurfDomain$, and $\EdgeSet$ of
edges $E$ covering $\SurfDomain$ without overlaps.  Each cell is a
curved triangle lying on $\SurfDomain$ and each edge is the geodesic
curve on $\SurfDomain$ connecting the vertices.  At every vertex a
frame is built from the nodal gradient of the height, as in
Section~\ref{sec:intrinsic}. Cell and edge frames at the cell centers
$M_{\Cell}$ and edge midpoints $M_{\Edge}$ are obtained by linear
interpolation of the nodal frames followed by re-orthogonalization.
The tangent $\tangentEdge_{\Edge}$ of the edge curve at $M_{\Edge}$
gives the edge conormal
$\normalEdge_{\Edge}=\tangentEdge_{\Edge}\times\normalvec_{\Edge}$.
Cell areas $|\Cell|$ and edge lengths $\edgeLength[\Edge]$ are those
of the flat triangle with the same vertices and of the chord.  All
these quantities are second-order accurate \citep{BachiniPutti2020},
except the edge length, whose relative error against the curvilinear
one is of first order.

The projected model \eqref{eq:kw}, used only for comparison, is discretized on the projected triangulation $\ProjSymb(\TriangH)$, with cell areas $|\ProjSymb(\Cell)|$ and the bed height as external potential.
The two models share the topology and the rain volume of every cell, so that their differences isolate the effects of the geometry.

\begin{remark}[Unique tangent plane on edges]\label{rem:edgeframe}
  The two cells $\Cell$ and $\Cell'$ sharing an edge $\Edge$ have different tangent planes. Thus, the normal component of a tangent vector across $\Edge$ is not defined until a common plane is chosen.
  Every edge carries a \emph{unique} frame at $M_{\Edge}$, shared by both cells, with unit conormal $\normalEdge_{\Edge}$ and tangent $\tangentEdge_{\Edge}$.
  When seen from $\Cell$, the conormal is oriented outward, $\normalEdge_{\Cell,\Edge}$, so that $\normalEdge_{\Cell',\Edge}=-\normalEdge_{\Cell,\Edge}$.
  Edge lengths and normal fluxes are then single-valued up to sign, which is requirement (iii) of Section~\ref{sec:intro}.
\end{remark}

\paragraph{Transfer between frames.}\label{sec:transfer}
We first note that, while in flat domains a piecewise constant vector
field is well defined, this is not the case on a surface.  In our
case, all quadratures are made with the midpoint rule for both cells
and edges.
Thus, following \cite{zhang2006vector}, we define a piecewise constant
vector field on a surface triangulation $\TriangH$ as a function that
assigns a tangent vector to the midpoints of each triangle or edge of
$\TriangH$.

Piecewise-constant tangent fields, such as the gradient of the
FEM-$\PONE$-interpolated height and the kinematic field $\AdvVel$ of
\eqref{eq:wfield}, are stored in the cell frames, due to the chosen FV
framework.  Whenever such a field has to be evaluated on an edge, it
is first transferred from the tangent plane of $\Cell$ to that of
$\Edge$.  We use the \emph{minimal rotation} taking the cell normal
$\normalvec_{\Cell}$ to the edge normal $\normalvec_{\Edge}$, i.e.\
the rotation about $\normalvec_{\Cell}\times\normalvec_{\Edge}$ by the
angle between the two normals, which for a tangent vector $\vecFunu$
($\vecFunu\cdot\normalvec_{\Cell}=0$) reads
\begin{equation}\label{eq:transfer}
  \Transfer_{\Cell\to\Edge}\,\vecFunu
  =\vecFunu-
  \frac{\vecFunu\cdot\normalvec_{\Edge}}
  {1+\normalvec_{\Cell}\cdot\normalvec_{\Edge}}\,
  \big(\normalvec_{\Cell}+\normalvec_{\Edge}\big).
\end{equation}
This transformation preserves lengths, maps
$\TanPlane[M_{\Cell}]{\SurfDomain}$ onto
$\TanPlane[M_{\Edge}]{\SurfDomain}$, and is a discrete parallel
transport. On flat triangles it reduces to the rotation about the edge
by half the dihedral angle, i.e., the Levi-Civita transport of
polyhedral surfaces \citep{PolthierSchmies1998}, and in the projected
geometry it is the identity.

With this transfer, every quantity used by the FV scheme and by the
conservative correction reduces to scalars per cell and edge.  For a
piecewise-constant tangent field $\cellfield$ we can define the normal
and tangent edge fluxes $g_{\Cell,\Edge}$ and $g^\tau_{\Cell,\Edge}$,
and the sum of normal edge fluxes $d_K$ as:
\begin{equation}\label{eq:gKE}
  g_{\Cell,\Edge}(\cellfield)=\edgeLength[\Edge]\,
  (\Transfer_{\Cell\to\Edge}\cellfield_{\Cell})
  \cdot\normalEdge_{\Cell,\Edge},\qquad
  g^\tau_{\Cell,\Edge}(\cellfield)=
  (\Transfer_{\Cell\to\Edge}\cellfield_{\Cell})
  \cdot\tangentEdge_{\Edge},
  \qquad
  d_{\Cell}(\cellfield)=
  \sum_{\Edge\subset\partial \Cell}g_{\Cell,\Edge}(\cellfield).
\end{equation}
\begin{remark}[Discrete (cellwise) divergence theorem]
  \label{rem:discrete-div-thm}
  The number $d_{\Cell}$ is defined on cell $\Cell$ and is the flux of
  the piecewise constant field $\cellfield$ through $\partial\Cell$,
  i.e.\ the part of its divergence resolved inside the cell.  Thus, it
  can be interpreted as a discrete divergence theorem.  On flat cells,
  as in the projected geometry, $d_{\Cell}=0$, i.e., all curvature
  sits on the edges.
  On curved cells no tangent field is constant, and $d_K\neq0$.
  
  For an exact surface and the Levi-Civita transport from $M_K$, the
  divergence theorem gives a flux of order $\meshparam^4$, while here
  the computed $d_K$ is dominated by the mismatch between the edge
  frames and the chord lengths and by the midpoint quadrature on the
  edges, so that $d_K/|K|=O(\meshparam)$.  The transfer contributes at
  higher order, since the minimal rotation coincides with the
  Levi-Civita transport along the geodesic from $M_K$ to $M_E$ up to a
  rotation of angle $O(\meshparam^3)$.  Note that the transfer affects
  the consistency of \eqref{eq:gKE}, never conservation.
\end{remark}

\begin{remark}[Linearized parallel transport]\label{rem:transfers}
  A linearized parallel transport based on the connection form of the
  Monge frame, in the spirit of \citep[eq.~(3.12)]{BachiniPutti2020},
  requires nodal second derivatives and differs from the minimal
  rotation by $O(\meshparam^2)$. Since the total error is dominated by
  the transport error, i.e.\ $O(\meshparam)$, our choice is based on
  consistency.
\end{remark}

\subsection{From terrain data to nodal geometry}
\label{sec:nodalgeom}\label{sec:mesh}

On real terrain the input is a surface triangulation given just by
surface points and connectivity (approximated by piecewise linear
cells), and only the heights at the mesh vertices are known as
geometric information.  From this input, we obtain the nodal gradient
by a weighted quadratic least-squares fit in the Monge form, pinned at
each vertex, over the $k$-ring of the vertex in the mesh graph.  By
definition, this fit is exact on quadratic surfaces.  Denser elevation
points can be added to the fit when available, resulting in a more
accurate reconstruction procedure.  The discrete surface is defined
once as vertices plus nodal frames and then every quantity is computed
from it. We call a k-ring related to a node $i$ of the triangulation
the set of vertices directly connected to $i$ by k edges along (the
shortest path).  On real data we use the 2-ring, to be able to recover
second order derivatives. In fact, the pinned quadratic has five
coefficients and the 2-ring can determine it also at the boundary
vertices, where the 1-ring has fewer than five points.
The fit is well determined when the mesh size is at least two to three
times the spacing of the elevation source, which is also the regime in
which sampling the elevation at the vertices stops aliasing the
raster.
A consequence of this kind of fitting technique is that the fit window
plays the role of a regularization parameter, and its choice can
amplify the artifacts of the input data, such as noise or meshing
stripes.

Beyond these reconstruction rules, the terrain must satisfy the
requirements discussed in Section~\ref{sec:admissibility} at the
discrete level.  Vertex-level conditions (no vertex other than the
outlets is a local minimum, no two neighbors have equal height) are
necessary but not sufficient.  Whether water can leave every cell
depends on the discrete velocity field and on the geometry in which it
is computed.  In Section~\ref{sec:traps}, we show that the discrete
traps of a given mesh can be computed in advance.  The input must also
resolve the terrain at the mesh scale.  Refining a piecewise-linear
interpolant without changing the surface turns its edges into V-shaped
creases, along which a single-valued normal flux cannot carry the
water.
Real elevation data and quantities derived from it often violate this
assumption. Indeed, depitted rasters, triangulations with meshing
stripes and piecewise-linear creases all can produce traps that no
choice of velocity field removes.  Bringing such data to match the KW
model assumptions without distorting the flow timing is a problem of
its own, which we leave to future work.
In this work we therefore assume a resolved and admissible input.

\section{Finite volume discretization}\label{sec:fv}

This section discusses the numerical approach, the choices dictated by
the kinematic wave, and what the scheme requires of the edge velocity.

\subsection{The scheme}\label{sec:fvscheme}

We discretize \eqref{eq:kwint} with a cell-centered, first-order
finite volume scheme on the surface triangulation $\TriangH$, in the
framework of the intrinsic FV scheme of \citep{BachiniPutti2020}.  The
unknowns are the cell averages $\Depth_{\Cell}$ of the depth, and the
geometry enters only through the quantities of
Section~\ref{sec:triangulation}.
The cell averages are advanced by the explicit Euler method:
\begin{equation}\label{eq:fvupdate}
  \Depth_{\Cell}^{n+1}=\Depth_{\Cell}^n
  +\Dt^n\Big(\force^n_{\Cell}
  -\frac{1}{|\Cell|}\sum_{\Edge\subset\partial\Cell}Q^n_{\Cell,\Edge}\Big),
  \qquad \force^n_{\Cell}=\frac{1}{|\Cell|}\int_{\Cell}\force^n\,dA ,
\end{equation}
where the cell average $\force_{\Cell}$ of the source is an input of
the scheme.  Note that, for vertical rain, $\force=r\cos\AngleSymbol$,
and by \eqref{eq:rainvol} the rain volume of a cell,
$\RainVol{\Cell}=|\Cell|\force_{\Cell}=\int_{\ProjSymb(\Cell)}r\,dA$,
is the one received by its vertical projection. As a result, the
intrinsic and the projected schemes receive exactly the same volume in
every cell and the same total inflow.  The global time step must
satisfy a CFL condition:
\begin{equation}\label{eq:dt}
  \Dt^n=\min\Big(
    \Dt_{\max},\;c\,\min_{\Cell}
    \frac{|\Cell|\,\Depth_{\Cell}^n}{\mathrm{out}_{\Cell}^n}
  \Big),
  \qquad
  \mathrm{out}_{\Cell}^n=\sum_{\Edge\subset\partial\Cell}
  \big(Q^n_{\Cell,\Edge}\big)^+ ,
\end{equation}
and we choose $c=1/2$ (CFL constant).  As standard, the simulation
starts from a prescribed initial step, $\Dt^0$ on the dry bed, and
\eqref{eq:dt} adapts the step at every time level.  The chosen CFL
guarantees positivity, $\Depth_{\Cell}^{n+1}\ge(1-c)\Depth_{\Cell}^n$
for $c\le1$, and monotonicity of the update for $c\le3/5$ (since the
outflow of $\Cell$ is proportional to $\Depth_{\Cell}^{5/3}$, the
derivative of $\Depth_{\Cell}^{n+1}$ with respect to
$\Depth_{\Cell}^n$ is at least $1-\tfrac53c$).
Boundary conditions consider free outflow on the whole boundary, as is
the consistent default with our assumptions.

The definition of edge fluxes $ Q_{\Cell,\Edge}$ is now the key point
in the construction of the scheme.  Let
$w_{\Cell,\Edge}=-w_{\Cell',\Edge}$ be a single-valued outward normal
component of the kinematic field $\AdvVel$ of \eqref{eq:wfield}, i.e.,
an approximation of the edge average of
$\AdvVel\cdot\normalEdge_{\Cell,\Edge}$, on every edge $\Edge$ shared
by $\Cell$ and $\Cell'$.  Then, the volumetric flux across $\Edge$,
the discrete counterpart of $\Qvect=\Depth^{5/3}\AdvVel$, is upwinded
along the sign of $w_{\Cell,\Edge}$:
\begin{equation}\label{eq:fvflux}
  Q_{\Cell,\Edge}=
  \edgeLength[\Edge]\,\Depth_{\mathrm{up}}^{5/3}\,w_{\Cell,\Edge},
  \qquad
  \Depth_{\mathrm{up}}=
  \begin{cases}
    \Depth_{\Cell} & w_{\Cell,\Edge}>0\\
    \Depth_{\Cell'} & w_{\Cell,\Edge}\le0
  \end{cases}
\end{equation}
and $U_{\Cell,\Edge}=\Depth_{\mathrm{up}}^{2/3}w_{\Cell,\Edge}$ is the
physical edge velocity, the counterpart of
$\vectVel=\Depth^{2/3}\AdvVel$.
By construction, the scheme is exactly mass conservative for every
single-valued edge velocity but not all of these are locally
conservative elementwise velocity fields. How to obtain a
single-valued locally conservative edge velocity is the subject of the
next section.

\subsection{Edge velocities}\label{sec:sources}

As previously noted, the scheme needs one normal velocity per edge,
whereas the kinematic field $\AdvVel_{\Cell}$ computed from the
$\PONE$ interpolant $\zbottom^h$ of the vertex heights gives two
one-sided values $g_{\Cell,\Edge}(\AdvVel)/\edgeLength[\Edge]$ and
$-g_{\Cell',\Edge}(\AdvVel)/\edgeLength[\Edge]$, which differ in
general on the surface triangulation.  For a piecewise-constant
tangent field $\cellfield$, the simplest single-valued edge quantity
is the average normal flux in the unique edge frame,
\begin{equation}\label{eq:avgflux}
  \bar g_{\Cell,\Edge}(\cellfield)
  =\frac12\big(g_{\Cell,\Edge}(\cellfield)-g_{\Cell',\Edge}(\cellfield)\big)
  =\frac{\edgeLength[\Edge]}{2}
  \big(\Transfer_{\Cell\to\Edge}\cellfield_{\Cell}
  +\Transfer_{\Cell'\to\Edge}\cellfield_{\Cell'}\big)
  \cdot\normalEdge_{\Cell,\Edge},
\end{equation}
single-valued up to sign,
$\bar g_{\Cell',\Edge}=-\bar g_{\Cell,\Edge}$.  It differs from the
one-sided value by half the \emph{normal jump} of the field across
$\Edge$,
\begin{equation}\label{eq:jump}
  \bar g_{\Cell,\Edge}(\cellfield)=
  g_{\Cell,\Edge}(\cellfield)-\frac12 J_{\Edge}(\cellfield),
  \qquad
  J_{\Edge}(\cellfield)=g_{\Cell,\Edge}(\cellfield)+g_{\Cell',\Edge}(\cellfield),
\end{equation}
which is symmetric in $\Cell$ and $\Cell'$, vanishes when the two
cells agree on the flux across $\Edge$, and is positive when the field
converges on $\Edge$.
We compare three strategies for constructing edge velocities, based
either on the raw cell fields or on their conservative reconstruction:
\begin{description}
\item[\FVS.]
  The raw (uncorrected) cell field, with flux-vector splitting,
  \begin{equation}\label{eq:fvs}
    Q_{\Cell,\Edge}=
    \Depth_{\Cell}^{5/3}\,g_{\Cell,\Edge}(\AdvVel)^+
    -\Depth_{\Cell'}^{5/3}\,g_{\Cell',\Edge}(\AdvVel)^+ .
  \end{equation}
  The flux is single-valued and globally conservative, but each cell
  sends across $\Edge$ what its own velocity pushes out.  It is
  possible to have \emph{counterflow} edges, i.e.\ edges along which
  both one-sided fluxes are positive, and thus the two cells exchange
  water and local conservation is lost.
\item[\AVG.] The plain average of the slope: with the normal slope
  $s^\nu_{\Edge}=\bar
  g_{\Cell,\Edge}(\GradSurf\zbottom^h)/\edgeLength[\Edge]$ and the
  tangential slope
  $s^\tau_{\Edge}=\tfrac12\big(g^\tau_{\Cell,\Edge}+g^\tau_{\Cell',\Edge}\big)$,
  with
  $g^\tau_{\Cell,\Edge}=g^\tau_{\Cell,\Edge}(\GradSurf\zbottom^h)$,
  along the edge tangent. Local conservation is not guaranteed.  Then,
  the Manning law is applied on the edge,
  \begin{equation}\label{eq:kslope}
    w_{\Cell,\Edge}
    =-\frac{s^\nu_{\Edge}}{\Manning\,|\mathbf s_{\Edge}|^{1/2}},
    \qquad
    |\mathbf s_{\Edge}|=\sqrt{(s^\nu_{\Edge})^2+(s^\tau_{\Edge})^2},
  \end{equation}
  where $|\mathbf s_{\Edge}|$ locally approximates $\sin\AngleSymbol$
  in the intrinsic geometry and $\tan\AngleSymbol$ in the projected
  one.  Where $|\mathbf s_{\Edge}|=0$, the edge velocity takes its
  limit $0$, since $|s^\nu_{\Edge}|\le|\mathbf s_{\Edge}|$.
  \begin{remark}
    The Manning law could also be applied before averaging, i.e., by
    taking the plain average of the kinematic field
    $\bar g_{\Cell,\Edge}(\AdvVel)/\edgeLength[\Edge]$.  Since the
    Manning law is smooth away from zero slope, the two averages
    differ by a term of order
    $|\Delta\mathbf s_{\Edge}|^2/|\mathbf s_{\Edge}|^2$ relative to
    the edge velocity, with
    $\Delta\mathbf
    s_{\Edge}=\Transfer_{\Cell\to\Edge}(\GradSurf\zbottom^h)_{\Cell}-\Transfer_{\Cell'\to\Edge}(\GradSurf\zbottom^h)_{\Cell'}$
    the jump of the slope across $\Edge$.  Thus, the two averages
    differ appreciably only where the jump is comparable with the
    slope itself, on counterflow edges and near stationary points.  We
    therefore consider only \AVG\ as the single uncorrected average.
  \end{remark}
\item[\LNs, \LNw.]  The conservative corrections of the slope and of
  the kinematic field, respectively. The average flux of the field is
  replaced by a corrected flux $F_{\Cell,\Edge}=-F_{\Cell',\Edge}$,
  balanced in every cell so that local conservation is achieved.  For
  \LNs, the Manning law \eqref{eq:kslope} is applied with the
  corrected normal slope
  $s^\nu_{\Edge}=F_{\Cell,\Edge}(\GradSurf\zbottom^h)/\edgeLength[\Edge]$
  in place of the average.  For \LNw, the corrected flux gives the
  edge velocity directly,
  \begin{equation}\label{eq:kln}
    w_{\Cell,\Edge}=F_{\Cell,\Edge}(\AdvVel)/\edgeLength[\Edge].
  \end{equation}
  The construction and analysis of these corrections is the subject of
  Section~\ref{sec:ln}, and one main focus of this work.
\end{description}

\begin{figure}\centering
\includegraphics[width=0.75\textwidth]{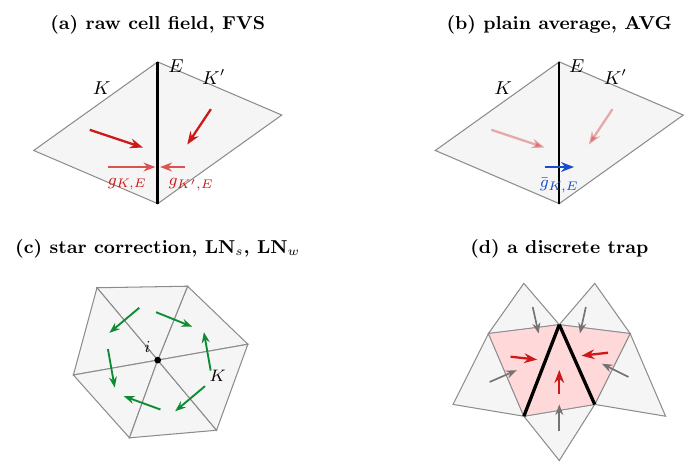}
\caption{Sketch of the edge velocities on cells $\Cell$, $\Cell'$
  sharing an edge $\Edge$.  (a) Raw cell field (\FVS), constant in
  each cell: the one-sided normal components $g_{\Cell,\Edge}$ and
  $g_{\Cell',\Edge}$ both point out of their own cell, so $\Edge$ is a
  counterflow edge and $J_{\Edge}>0$.  (b) Plain average (\AVG): one
  value per edge, the jump is dropped.  (c) Star correction (\LNs,
  \LNw, Section~\ref{sec:lnconstruction}): the jumps of the interior
  edges through the vertex $i$ are shared among the cells of its star,
  and the corrected flux is balanced in every cell (local
  conservation). (d) A discrete trap (Section~\ref{sec:traps}): a
  closed group of cells connected by counterflow edges, which water
  enters and never leaves.}
\label{fig:sketch}
\end{figure}

We call \emph{single-valued edge velocities} those that define one
normal velocity $w_{\Cell,\Edge}=-w_{\Cell',\Edge}$ per edge: \AVG,
\LNs, \LNw. \FVS\ gives a single-valued flux but keeps the two
one-sided velocities.  Fig.~\ref{fig:sketch}, panels (a) and (b)
sketch the raw field and the plain average on one edge.

\subsection{Discrete steady state and discrete traps}\label{sec:traps}

The flux of the kinematic wave model is linear in $\Phi=\Depth^{5/3}$,
with nonnegative weights that depend on the choice of edge velocities.
This means that the flux from a cell $K$ to a neighbor $J$ is
$\lambda_{KJ}\Phi_K$, with $\lambda_{KJ}\ge0$, which depends on the
chosen edge velocities as described in the previous section.  At the
outlet, the outflow from cell $K$ through the outlet edges (boundary)
is given by $\lambda^{\mathrm{out}}_K\Phi_K$.
In this notation, the steady form of \eqref{eq:fvupdate}, multiplied
by cell area $|K|$, states that source plus inflow equals outflow on
each cell $K$ or, explicitly,
$\RainVol{K}+\sum_J\lambda_{JK}\Phi_J=\lambda_K\Phi_K$, with
$\RainVol{K}$ the source in the $K$-th cell.
Under a constant source (e.g., rain constant in time), the steady
state of the scheme satisfies the following sparse linear system:
\begin{equation*}
  \mathsf A\Phi=\RainVol{},
\end{equation*}
where $\RainVol{}$ is the vector collecting the total rainfall volume
on each cell while matrix $\mathsf A$ contains the information on the
flux exchange across edges.  More precisely,
$\mathsf A=\mathrm{diag}(\lambda)-\boldsymbol\Lambda^{\!\top}$, with
$\lambda=\{\lambda_K\}$,
$\lambda_K=\sum_J \lambda_{KJ}+\lambda^{\mathrm{out}}_K$ the total
outflow weight, and $\boldsymbol\Lambda=\{\lambda_{KJ}\}$.
The following proposition states that if every cell of the mesh is
connected to the outlet by at least one directed path, then the steady
state depth arising from the FV numerical solution of the KW model is
unique and necessarily nonnegative.
\begin{proposition}[Steady state]\label{prop:steady}
  The matrix $\mathsf A$ is a nonsingular M-matrix if and only if
  every cell has a directed path, in the transfer directed graph
  $K\to J$ ($\lambda_{KJ}>0$), to a cell with outlet outflow
  ($\lambda_K^{\mathrm{out}}>0$).  In that case the steady state is
  unique and $\Phi=\mathsf A^{-1}\RainVol{}\ge0$.
\end{proposition}
\begin{proof}
  By definition, an \emph{M-matrix} is a matrix in which all the
  off-diagonal elements are less than or equal to zero (the
  fundamental property of a \emph{Z-matrix}) and the eigenvalues have
  all positive real parts. A fundamental property of M-matrices is
  that they have nonnegative inverse, i.e., $\mathsf{A}^{-1}x\ge0$ for
  all $x\ge0$. Here $x\ge0$ means a vector with nonnegative
  elements.

  It is easy to see that the matrix $\mathsf
  A^{\!\top}=\mathrm{diag}(\lambda)-\boldsymbol\Lambda$ is a Z-matrix
  whose row $K$ has excess
  $\lambda_K-\sum_J\lambda_{KJ}=\lambda^{\mathrm{out}}_K\ge0$, with
  strict inequality occurring only for the outlet cells, and whose
  off-diagonal pattern is the transfer directed graph.
  Thus, if every cell has a directed path to a strictly dominant row,
  $\mathsf A^{\!\top}$ is weakly chained diagonally dominant with
  positive diagonal, hence a nonsingular M-matrix
  \citep{ShivakumarChew1974,AzimzadehForsyth2016}.  The same is true
  for $\mathsf A$, and $\mathsf A^{-1}\ge0$.  If some cell $K$ has no
  such path, the set $\mathcal C$ of cells reachable from $K$ is
  closed and has no outlet.  Ordering $\mathcal C$ first,
  $\mathsf A^{\!\top}$ is block lower triangular with a diagonal block
  $\mathsf B$ such that $\mathsf B\boldsymbol 1=0$, hence singular.
\end{proof}

The steady state is thus obtained directly from one sparse linear
solve, up to the solver tolerance, independently of the time step, and
it is a multiple-flow-direction accumulation on the dual graph, with
transfer fractions $\lambda_{KJ}/\lambda_K$ fixed by the edge
velocity.
This places the scheme next to the flow-routing algorithms on grids
and triangulations
\citep{Quinn1991,Freeman1991,Tarboton1997,Tucker2001}.
A \emph{discrete trap} is a closed class of the directed graph without
outlet: water entering this class never leaves.  Traps depend on the
edge velocity and on the mesh only, not on the source (rain), the
Manning coefficient or the time integrator, and they can be counted
without simulation.  We quantify them by computing the \emph{trapped
  rain}, the fraction of the rain falling on the domain that drains
into a trap.  For \FVS, a closed group of cells connected by
counterflow edges is a trap (see Fig.~\ref{fig:sketch}, panel (d)).

\section{Conservative edge fields from elevation data}\label{sec:ln}

This section constructs the corrected fluxes of \LNs\ and \LNw.  In
both cases the average flux \eqref{eq:avgflux} of a piecewise-constant
cell field, the slope $\GradSurf\zbottom^h$ or the kinematic field
$\AdvVel$, is corrected by the post-processing of
\citet{LarsonNiklasson2004} into a single-valued flux that is balanced
in every cell thus ensuring local conservation. The original
correction was designed for the gradient of a $\PONE$ Galerkin
solution of an elliptic problem with a given source, on flat cells.
Instead, we have a field that is derived from data, the cells are
curved, no source is given, and the kinematic field is not even a
gradient.  Here, we just recall the construction, and focus on what
changes: the edge scalars and the curved cells
(Section~\ref{sec:lnconstruction}), the load, which is the main
difference (Section~\ref{sec:load}), and the boundary
(Section~\ref{sec:boundary}).  Section~\ref{sec:divform} states what
the correction imposes on the corrected field.
Our algorithm is summarized in Appendix~\ref{app:algorithm}.

\subsection{The \LNs-\LNw\ construction on curved cells}
\label{sec:lnconstruction}\label{sec:star}\label{sec:compat}

Let $\cellfield$ be a piecewise-constant tangent field on $\TriangH$, with the scalars $g_{K,E}=g_{K,E}(\cellfield)$ and $d_K=d_K(\cellfield)$ of \eqref{eq:gKE}, the average flux $\bar g_{K,E}$ of \eqref{eq:avgflux}, and the jump $J_E$ of \eqref{eq:jump}.
On a boundary edge we set $\bar g_{K,E}=g_{K,E}$, the one-sided value.
For each node $i$, the nodal star $\Star{i}$ is the set of cells
sharing node $i$.  The construction in \citet{LarsonNiklasson2004}
works locally, star by star, to restore the balance between the source
and the discrete cell divergence,
$\sum_{E\subset\partial K}\bar g_{K,E}$.  For every vertex $i$ with
star $\Star{i}$ and $\PONE$ basis function $\nodalBasis{i}$, we call
the load, or forcing, $b_i$ the share of the source attributed to $i$,
i.e., the $\PONE$-projection of the source term on the star and we let
$b_{i,K}$ be the $K$-th element contribution to $b_i$ such that
$b_i=\sum_{K\in\Star{i}}b_{i,K}$. Then both \LNs-\LNw\ corrections can
be organized in the following three steps:
\begin{enumerate}[label=(\arabic*),leftmargin=*]
\item \emph{Split residual.} For $K\in\Star{i}$,
  \begin{equation}\label{eq:Rik}
    R_{i,K}=b_{i,K}+\sum_{\substack{E\subset\partial K\\ i\in E}}
    \frac12\big(\bar g_{K,E}-g_{K,E}\big)
    =b_{i,K}-\frac14\sum_{\substack{E\subset\partial K\ \text{interior}\\ i\in E}}J_E ,
  \end{equation}
  where the second form follows from \eqref{eq:jump} and boundary
  edges give no contribution. The residual $R_{i,K}$ is the share of
  cell $K$ in the balance of vertex $i$ that is left over when the
  one-sided flux $g_{K,E}$ is replaced by the average, with the factor
  $\tfrac12$ corresponding to
  $\int_E\nodalBasis{i}\,ds/\edgeLength[E]$.
\item \emph{Star problem.} Find $V_i$, constant on each cell of
  $\Star{i}$, such that
  \begin{equation}\label{eq:starproblem}
    \sum_{\substack{E\subset\partial K\ \text{interior}\\ i\in E}}
    \frac{\edgeLength[E]}{2}\,(V_{i,K}-V_{i,K'})=-R_{i,K},
    \qquad K\in\Star{i},
  \end{equation}
  where $K'$ is the neighbor of $K$ across the interior edge $E$. The
  left-hand side corresponds to the Laplacian on the graph arising
  from the cells of the star, with weight
  $\edgeLength[E]/2=\int_E\nodalBasis{i}\,ds$ on each interior edge
  through $i$, and the right-hand side is the residual.
\item \emph{Correction.} On every interior edge the average flux is
  corrected as follows:
  \begin{equation}\label{eq:Fcorr}
    F_{K,E}=\bar g_{K,E}+
    \sum_{i\in E}\frac{\edgeLength[E]}{2}\,(V_{i,K}-V_{i,K'}),
  \end{equation}
  and $F_{K,E}=g_{K,E}$ on a boundary edge.  The correction is the
  flux of the star potentials of the two vertices of $E$ across
  $E$. It is antisymmetric in $K$ and $K'$, so the corrected flux is
  single-valued, and it only moves flux between neighboring cells of a
  star (see Fig.~\ref{fig:sketch}, panel (c)).
\end{enumerate}
The graph Laplacian in \eqref{eq:starproblem} has no Dirichlet term
and is thus singular having the constants as kernel. As a consequence,
the star problem is solvable if and only if
$\sum_{K\in\Star{i}}R_{i,K}=0$.  When it is solvable, the additive
constant does not affect the edge flux $F$, which depends only on
differences.

Main adaptations are required in our case since the original work in
\citet{LarsonNiklasson2004} starts from a vector field that is the
gradient of the $\PONE$-Galerkin solution $u_h$ of $-\Lap u=\sigma$
and the element load $b_{i,K}=\int_K\sigma\,\nodalBasis{i}\,dA$ is the
elemental contribution of the forcing $\sigma$.  The term $R_{i,K}$
from~\eqref{eq:Rik} is the residual of the Galerkin equation tested
with $\nodalBasis{i}$, which, when summed over the
  star cells and assuming no quadrature nor geometric errors, vanishes
  by Galerkin orthogonality.
  Two observations arise. One is that the integral
  $\int_K\cellfield_K\cdot\GradSurf\nodalBasis{i}\,dA$ of the Galerkin
  equation cannot be calculated for a general vector field
  $\cellfield_K$ and is thus replaced by its edge form to yield:
\begin{equation}\label{eq:ah}
  a_h(\cellfield;\nodalBasis{i})
  =\sum_{K\in \Star{i}}\int_K\cellfield_K\cdot\GradSurf\nodalBasis{i}\,dA
  =\sum_{K\in \Star{i}}\sum_{\substack{E\subset\partial K\\ i\in E}}
  \frac{g_{K,E}(\cellfield)}{2} .
\end{equation}
This equation is the analogue of integration by parts on the graph associated with the mesh cells.
Note that, for $\cellfield=\GradSurf u_h$ the above equation also
defines the stiffness matrix $\mathsf K_h$,
$a_h(\GradSurf u_h;\nodalBasis{i})=(\mathsf K_h u)_i$, which coincides
with the scalar product of $i$-th row of the matrix-vector product of
the $\PONE$ stiffness matrix and the vector of DOFs.  A second thing
to note is that, as seen in Remark~\ref{rem:discrete-div-thm}, on a
curved cell the edge fluxes of a constant tangent field computed
by~\eqref{eq:avgflux} do not sum to zero. This new term must enter the
cell balance. This is collected in the following Lemma.

\begin{lemma}[Compatibility and element balance]
  \label{lem:compat}
  Let $\cellfield$ be any piecewise-constant tangent field on
  $\TriangH$, and let every edge have a unique length and a
  single-valued average flux \eqref{eq:avgflux}. Then:
  \begin{enumerate}[(i)]
  \item\label{item:i}
    the star  problem of vertex $i$ is solvable if and only if the following sum vanishes:
    \begin{equation*}
      \sum_{K\in \Star{i}}R_{i,K}=b_i-a_h(\cellfield;\nodalBasis{i})
      +\tfrac12\sum_{E\subset\partial\SurfDomain,\,i\in E}g_{K,E}
      =b_i-\tfrac12\sum_{\substack{E\ \text{interior}\\ i\in E}}J_E=0 ;
    \end{equation*}
  \item\label{item:ii} if every star problem containing $K$ is solved,
    the corrected fluxes satisfy the local balance
      equation:
    \begin{equation}\label{eq:elbal}
      \sum_{E\subset\partial K}F_{K,E}=d_K(\cellfield)-\sum_{i\in K}b_{i,K},
    \end{equation}
    i.e., every cell is conservative with respect to the effective
    cell source $\sum_{i\in K}b_{i,K}-d_K(\cellfield)$.
  \end{enumerate}
\end{lemma}
\begin{proof}
  Let $\cellfield$ be any piecewise-constant tangent field.
  \begin{enumerate}[(i)]
  \item Summing the second form of \eqref{eq:Rik} over $\Star{i}$,
    every interior edge through $i$ belongs to two cells of the star
    and it contributes exactly twice with $-\tfrac14J_E$, while
    boundary edges contribute nothing.  Hence,
    $\sum_{K}R_{i,K}=b_i-\tfrac12\sum_{E}J_E$, and the first form
    follows from~\eqref{eq:ah} and~\eqref{eq:jump}.  The star graph is
    connected also at a boundary vertex, where the cells form a path,
    so the star Laplacian is singular with the constants as kernel.
    Hence the star problem is solvable, i.e., the
      right-hand-side is in the image of the star Laplacian matrix,
    if and only if the sum of the residuals vanishes.
  \item In the sum of \eqref{eq:starproblem} over the three vertices
    of $K$, each edge of $K$ collects $\bar g_{K,E}-g_{K,E}$ once, and
    the potential terms assemble into $F_{K,E}-\bar g_{K,E}$ by
    \eqref{eq:Fcorr}.  Hence,
    $\sum_EF_{K,E}-\sum_E\bar g_{K,E}=-\sum_ib_{i,K}-\sum_E\bar
    g_{K,E}+d_K$, which is~\eqref{eq:elbal}.
  \end{enumerate}
\end{proof}

\begin{remark}
  The identity in~\ref{item:i} requires that both cells see the same
  conormal on the shared edge, which is what the unique edge frame
  provides, and that the load be computed with the same scalars
  $g_{K,E}$ (same frames, same transfer) as the average flux.  The
  term $d_K$ in~\ref{item:ii} is the flux of $\cellfield_K$ through
  $\partial K$ as computed by \eqref{eq:gKE}.  The distributional
  divergence of a vector field splits into a part concentrated on the
  edges, carried by the jumps $J_E$, and a part inside the cells.  The
  latter vanishes for a piecewise vector field on flat domains.  By
  the divergence theorem on $\SurfDomain$ and our choice of quadrature
  rules, $d_K$ is the edge-midpoint quadrature of the integral over
  $K$ of the surface divergence of $\cellfield_K$ (and transferred to
  the edge tangent planes), and thus the effective source in (ii)
  collects the divergence carried by the jumps through the loads
  $b_{i,K}$ and the divergence inside the cell through $d_K$.  Its
  size is set by the discretization of the cell geometry rather than
  by the curvature (Section~\ref{sec:transfer}).  Moreover, it
  vanishes on flat cells, where the original construction is
  recovered, and it must be kept on curved cells for the corrected
  flux to be conservative.
\end{remark}

\subsection{The load}\label{sec:load}

For a field derived from elevation data there is no explicit
load. Recall that the load is the source of the equation that the
corrected field would satisfy, i.e., minus the divergence of the
continuous vector field.
The KW equation can be fed by the vector field formed by the
bottom slope or by the Manning-weighted kinematic field, $\AdvVel$,
see~\eqref{eq:kwint}.
The bottom gradient must satisfy~\eqref{eq:LB}, i.e.,
$\LapSurf\zbottom=\cos\AngleSymbol\,H$, and thus the load is
proportional to the mean curvature $H$.
Away from the stationary points of $\zbottom$, the mean curvature
splits through each point along the two orthogonal directions given by
the elevation contour and the steepest descent, $H=k_h+k_v$ (Euler's
theorem).  Here $k_h$ and $k_v$ are the normal curvatures of
$\SurfDomain$, i.e., the curvatures along these two directions, taken
with the sign of $H$ in \eqref{eq:LB} (positive in valleys).  Thus,
since by \eqref{eq:gradG}
$\ABSsurf{\GradSurf\zbottom}=\sin\AngleSymbol$ and assuming a constant
Manning coefficient $\Manning$, the divergence of the bottom gradient
and Manning-weighted kinematic fields can be written as:
\begin{equation}\label{eq:divfields}
  \LapSurf\zbottom=\cos\AngleSymbol\,\big(k_h+k_v\big),
  \qquad
  \DivSurf\AdvVel=
  -\frac{\cos\AngleSymbol}{\Manning\sqrt{\sin\AngleSymbol}}\,
  \Big(k_h+\frac12\,k_v\Big).
\end{equation}
Note that the two fields scale differently with the inclination.
The divergence of the slope stays bounded by the principal curvatures, while the divergence of the kinematic field grows as $(\sin\AngleSymbol)^{-1/2}$ and is unbounded at the stationary points of the terrain.
In both cases the flow field, $-\GradSurf\zbottom$ or $\AdvVel$, has negative divergence, i.e., water accumulates, where the combination of the two curvatures in \eqref{eq:divfields} is positive, e.g., where the flow converges ($k_h>0$) and decelerates ($k_v>0$).

We discuss next three possible choices for the load.
\begin{description}[leftmargin=0pt,font=\normalfont\itshape]
\item[Zero load: $b_i=0$.]
  With this choice the corrected field is forced to be divergence free in every cell, up to $d_K$.
  The case of bottom slope treats $\zbottom^h$ as a harmonic height, i.e., $\LapSurf\zbottom=0$.
  Since $\cos\AngleSymbol>0$, this means that the principal curvature is zero $H=0$, i.e., the bottom surface should be a minimal surface in which every concave bend of the terrain in one principal direction is offset by an equal convex bend in the other.
  For the case of the kinematic field, a uniform water layer would cross every cell without accumulating or depleting.
  In both cases the convergence in the valleys and the divergence on the ridges, which route the water, are removed.
  Moreover, the jumps of the field do not add up to zero around a vertex, condition (i) fails, and the star problems are not even solvable.
\item[Exact load:
  $b_i=-\int_{\Star{i}}\DivSurf\cellfield\,\nodalBasis{i}\,dA$.]  This
  is the theoretically correct load, being the weak divergence of the
  vector field $\cellfield$. However, this choice is practically
  available only for analytical surfaces since it requires the
  evaluation of second derivatives of the elevation function, a
  difficult task for DEM-based elevation data.
  Moreover, even when this information is available, this load satisfies
  condition~\ref{item:i} of Lemma~\ref{lem:compat} only in the limit
  $\meshparam\to0$, while at any finite resolution the star problems are not
  solvable, and thus the exact load cannot be used by the construction.
  In addition, in case of
  correction of the kinematic field, this integral is unbounded at the
  stationary points of the terrain.
\item[Compatible load (Lemma \ref{lem:compat}).]  Our compatible load
  is nothing else than an ad-hoc second order consistent quadrature
  rule for the exact load above, developed by constraining the flux
  balance to be accurate to machine precision. In this case we have:
  \begin{equation}\label{eq:compload}
    b_i=\frac{1}{2}\sum_{\substack{E\ \text{interior}\\ i\in E}}J_E(\cellfield)
    =a_h(\cellfield;\nodalBasis{i})
    -\tfrac12\sum_{E\subset\partial\SurfDomain,\,i\in E}g_{K,E}(\cellfield),
  \end{equation}
  By (i) of Lemma~\ref{lem:compat}, a single total load per vertex makes the star problem solvable and only its split among the cells of the star is free.
\end{description}

For the case of bottom slope, $\cellfield=\GradSurf\zbottom^h$, the
compatible load is $(\mathsf K_h\zbottom^h)_i$ minus the boundary
terms.  The discrete elevation field $\zbottom^h$ can be considered as
the Galerkin solution of the discrete problem, and on flat cells, as
in the projected geometry, Galerkin orthogonality holds exactly and
\LNs\ is the classical construction.

On curved cells the only additional term is $d_K$ and the load can be
given the following geometric interpretation.  On a polyhedral
surface, the stiffness matrix is the cotangent matrix, and
$-(\mathsf K_h\xv)_i/m_i$, $m_i=\sum_{K\ni i}|K|/3$ the lumped mass,
is the discrete mean-curvature vector
\citep{PinkallPolthier1993,Meyer2003} at vertex $i$.  On each cell,
the matrix $\mathsf K_h$ is the stiffness matrix of the
Laplace-Beltrami operator calculated by edge-quadrature. At an
interior vertex, $-b_i/m_i$ is the vertical component of this vector,
the discrete counterpart of $\LapSurf\zbottom=\cos\AngleSymbol\,H$:
terrain curvature is the source that the slope field must balance.  It
converges to $\cos\AngleSymbol\,H$ weakly at $O(\meshparam^2)$ on
every mesh, and pointwise only on regular meshes
\citep{Xu2004,Hildebrandt2006,Wardetzky2007}.
For the kinematic field there is no such reading. The cell field
$\AdvVel=-\GradSurf\zbottom^h/(\Manning\,\ABSsurf{\GradSurf\zbottom^h}^{1/2})$,
evaluated cell by cell, is parallel to the slope in each cell, but its
magnitude changes from cell to cell and there is no stiffness matrix
and no Galerkin problem behind it. Indeed, the compatible load
\eqref{eq:compload} with $\cellfield=\AdvVel$ is the weak divergence
of the cell field assembled from its own jumps, the discrete
counterpart of $-\int\DivSurf\AdvVel\,\nodalBasis{i}\,dA$.  Here
compatibility is not a reinterpretation of Galerkin orthogonality, but
it is the only available definition of the load.

We distribute the compatible load in proportion to the cell areas
(lumped split), $b_{i,K}=\tfrac{|K|}{3}\,b_i/m_i$. By \eqref{eq:Rik},
the star problem of $i$ then reads:
\begin{equation}\label{eq:starrhs}
  \sum_{\substack{E\subset\partial K\ \text{interior}\\ i\in E}}
  \tfrac{\edgeLength[E]}{2}\,(V_{i,K}-V_{i,K'})
  =\frac14\sum_{\substack{E\subset\partial K\ \text{interior}\\ i\in E}}J_E
  -\frac{|K|}{3}\,\frac{b_i}{m_i},\qquad K\in\Star{i}\,.
\end{equation}
The right-hand side sums to zero over the star as the jump of every
interior edge through $i$ is removed from the two cells that share the
edge and redistributed over the star in proportion to the cell areas.

\subsection{Boundary treatment}\label{sec:boundary}

In \citet{LarsonNiklasson2004} the stars of the Dirichlet vertices are
solved too (the one-sided boundary terms make them non-singular), and
element-wise conservation relies on solving the star problems at all
vertices.
Here, we treat every boundary vertex as a Neumann vertex: the boundary
flux is the one-sided value $g_{K,E}$, uncorrected, and the boundary
star is solved with the compatible load \eqref{eq:compload}, which
involves only the interior edges through $i$, since a boundary edge
carries no jump (in the second form of \eqref{eq:compload}, the
one-sided boundary fluxes cancel against the corresponding part of
$a_h$).  With this choice every star problem is solvable, every
element is conservative by Lemma~\ref{lem:compat}, the boundary fluxes
are exactly those of the plain average, and the interior accuracy is
preserved.  For the kinematic wave the boundary fluxes matter mainly
at the outlet, since walls carry no flux or are water divides.

\subsection{What the correction imposes}
\label{sec:divform}\label{sec:what}\label{sec:whichfield}

As mentioned above, the compatible load has a direct reading in terms of discrete divergences.
Let
\begin{equation}\label{eq:Dnodal}
  D_i(\cellfield)=
  -\frac{1}{2m_i}\sum_{\substack{E\ \text{interior}\\ i\in E}}J_E(\cellfield)
\end{equation}
be the weak (lumped) nodal divergence of the piecewise-constant field, i.e.\ the assembly of its normal jumps, so that $b_i=-m_iD_i(\cellfield)$ at every vertex by \eqref{eq:compload}.
The following relation between the corrected flux and the divergence of the field follows directly from Lemma~\ref{lem:compat}.

\begin{corollary}[Imposed cell divergence]\label{cor:divform}
  With the compatible load and the lumped split, the corrected flux satisfies, in every cell,
  \begin{equation}\label{eq:balance}
    \sum_{E\subset\partial K}F_{K,E}=
    \frac{|K|}{3}\sum_{i\in K}D_i(\cellfield)+d_K(\cellfield).
  \end{equation}
\end{corollary}
In hydrological terms, every cell concentrates or spreads the flow
according to the convergence or divergence of the cell gradients. This
is a statement on where the flow is directed by the terrain, not on
mass, which every single-valued cell flux conserves.
Note that the plain average \eqref{eq:avgflux} lacks consistency in
the case of an unstructured mesh. Indeed its cell divergence differs
from the right-hand side of \eqref{eq:balance} by an error term whose
magnitude does not decrease under mesh
  refinement. On the other hand, on a structured mesh the average flux
  achieves first-order convergence.

\paragraph{Which field to correct.}
The corollary holds for both fields, and shows that the choice of the
corrected field is a choice of which constraint is enforced up to
machine error. \LNs\ imposes \eqref{eq:balance} on the slope, but
because the Manning law is nonlinear, the edge velocity
\eqref{eq:kslope} computed from the corrected slope does not inherit
this constraint.  \LNw\ imposes the same constraint but on the field
with which the scheme moves the water, whose divergence weighs the
convergence of the flow lines and the change of steepness along them
as in \eqref{eq:divfields}.  In the numerical tests neither correction
is uniformly more accurate than the other.
We use \LNw\ as the default, but slope correction remains a valid
alternative, reported throughout.

\begin{remark}[Single-cell traps]\label{rem:traps}
  For the corrected edge velocities, a single cell is a trap if and
  only if no edge carries water out of it.  For \LNw\ this means
  $F_{K,E}(\AdvVel)\le0$ on every edge of $K$, which by
  \eqref{eq:balance} requires a convergent imposed divergence,
  $\tfrac{|K|}{3}\sum_{i\in K}D_i(\AdvVel)+d_K(\AdvVel)\le0$.  For
  \LNs\ it requires $F_{K,E}(\GradSurf\zbottom^h)\ge0$ on every edge
  of $K$, hence, by \eqref{eq:balance} with
  $\cellfield=\GradSurf\zbottom^h$,
  $\tfrac{|K|}{3}\sum_{i\in
    K}D_i(\GradSurf\zbottom^h)+d_K(\GradSurf\zbottom^h)\ge0$, which leads
  a nonnegative discrete Laplace-Beltrami operator of the height,
  intra-cell term included, i.e., a valley-like cell.  A single-cell
  trap is therefore possible only where the data themselves converge,
  never on divergent (ridge) cells.  A sufficient condition for the
  absence of multi-cell traps, in terms of the mesh and of the
  piecewise-linear structure of the surface
  \citep{Banchoff1967,Edelsbrunner2003}, remains an open question.
\end{remark}

\begin{remark}[Structured vs unstructured grids]\label{rem:raster}
  On structured meshes, where every square of a grid is split along
  the same diagonal, the plain average is superconvergent. Its
  slope-induced error is of second order, its balance defect decreases
  at first order, and the correction does not improve this accuracy.
  However, in certain situations, approximation of a surface by a
  uniform quadrilateral tessellation is not a convergent process and
  certain geometric quantities, such as, e.g., the Gaussian curvature
  at a node, cannot be asymptotically recovered if the nodal valence
  (number of edges emanating from the node) is smaller than 6 and the
  secant edges connected to the node have different lengths (the
  so-called parallelogram property)~\cite{Morvan2008}.  Even
  triangulated rasters lose the parallelogram property wherever the
  terrain bends at the raster scale, and thus lose the
  superconvergence of the average wherever the curvature changes over
  a few cells, leading to non consistency of the approximation.

  For these reasons, unstructured triangulation of the terrain
  surface should be the preferred choice when working with surface
  partial differential equations.  Moreover, unstructured
  triangulations better conform to terrain features (channels,
  divides, breaklines, catchment boundaries) with locally refined
  cells, and possibly they can be put in direct connection with (be
  subsets or supersets of) surface boundary meshes of 3D meshes
  where coupled surface-subsurface models are
  defined~\citep{Camporese2010,Scudeler2016}.
\end{remark}

\section{Numerical experiments}
\label{sec:experiments}\label{sec:results}
\begin{table}\centering\small
    \setlength{\tabcolsep}{4pt}
  \caption{The four tests. All have $\Manning=0.03$\,s\,m$^{-1/3}$ and $c=1/2$.}\label{tab:tests}
\begin{tabularx}{\textwidth}{@{}>{\raggedright\arraybackslash}p{2.2cm}*{4}{>{\raggedright\arraybackslash}X}@{}}
 
    \toprule
    & \textbf{Test 1} & \textbf{Test 2} & \textbf{Test 3} & \textbf{Test 4} \\ \midrule
    \textbf{Surface} & sloping plane & convergent valley &  fully
                                                           3D surface
                                                           \citep{BachiniPutti2020}
                                                          & Fernow
                                                            catchment
                                                            (data)
                                                            \citep{JensenMcGuirePrince2017,BotterEtAl2021}\\[0.3em]
    \textbf{Domain} & $[0,\pv{t1.L}]\times[0,\pv{t1.W}]$~m & $[0,\pv{t2.L}]\times[0,\pv{t2.W}]$~m & $[0,20]\times[0,8]$~m & \pv{t4.area}\,km$^2$ \\[0.3em]
    \textbf{Parameters} & slopes=\pv{t1.slopes} & slopes=\pv{t2.S}, curv=\pv{t2.a} & slopes up to \pv{t3.slopemax} & slopes up to \pv{t4.slopemax} \\[0.3em]
    \textbf{Outlet} & $x=100$~m & $x=100$~m & whole boundary & whole boundary \\[0.3em]
    \textbf{Test purpose} & intrinsic vs projected & mesh dependence & degenerate stationary point & terrain data \\[0.3em]
    \textbf{Exact reference} & event hydrograph & steady state & steady, valley half & none \\[0.3em]
    \textbf{Meshes} & unstructured  & structured, unstructured  & unstructured & data mesh \\[0.3em]
    \textbf{Rain}  & \pv{t1.rain}~mm/h  & \pv{t2.rain}~mm/h & \pv{t3.rain}~mm/h & \pv{t4.rain}~mm/h \\[0.3em]
    \textbf{$t_r$ / $t_f$} &\pv{t1.train} / \pv{t1.tend}~min &
                                                               \pv{t2.train}
                                                               / \pv{t2.tend}~min
                                        & \pv{t3.train} / \pv{t3.tend}
                                          ~min & \pv{t4.train} / \pv{t4.tend}~min
    \\[0.3em]
    \textbf{$\Dt_{\max}$} & \pv{t1.dtmax}~s & \pv{t2.dtmax}~s & \pv{t3.dtmax}~s & \pv{t4.dtmax}~s \\
    \bottomrule
  \end{tabularx}
\end{table}

We report four numerical tests, summarized in
Table~\ref{tab:tests}. Test~1, a sloping plane with an exact event
hydrograph, compares the intrinsic and the projected kinematic wave,
and in this case the edge velocity plays no role. Test~2, a steep
convergent valley with an exact steady state, shows how the effect of
the locally conservative correction depends on the mesh. Test~3, the
surface of \citet[Sec.~4.4]{BachiniPutti2020}, admissible except at a
degenerate stationary point, shows that a single point of the terrain
can affect the whole response, and that the single-valued edge
velocities keep the scheme robust there. Test~4, a catchment of the
Fernow Experimental Forest, applies the whole method to real terrain data.

\subsection{Common setup}\label{sec:setup}\label{sec:verification}

\begin{table}\centering
\caption{Notation used in the numerical experiments.}\label{tab:notation}
\begin{tabular}{@{}>{\raggedright}p{0.22\textwidth}>{\raggedright\arraybackslash}p{0.72\textwidth}@{}}
\toprule
Paper notation & Meaning \\ \midrule
$\kEta$, $\kVol$ & intrinsic over projected steady depth, and stored volume per unit horizontal area
  \eqref{eq:depthratio}, \eqref{eq:volration}\\
$\bar t$ & mean response time: steady stored volume over steady outflow\\
transport error & error of the first-order upwind scheme, common to all edge velocities\\
slope-induced error & error due to the edge velocity alone, measured against the exact-slope reference\\
transfer directed graph & routing graph of the scheme: $J\to K$ when $J$ sends water to $K$\\
discrete trap & closed class of the directed graph without outflow: a pit of the routing\\
trapped rain & fraction of the rain that ends in a trap\\
\bottomrule
\end{tabular}
\end{table}

This section collects the setup, the reference solutions, and the
metrics that are common to all tests. Table~\ref{tab:notation} recalls
some of the notation.
Every event starts on a dry bed, with constant rain until $t_r$
followed by a recession time.
Unless stated otherwise, the rain intensity is $r=50$\,mm/h, the
Manning coefficient is $\Manning=0.03$\,s\,m$^{-1/3}$ and the time
step is computed by \eqref{eq:dt} with $c=1/2$.  The value of
$\Dt_{\max}$ only bounds the step on a dry bed and sets the time
resolution of the hydrographs.
The FV scheme is used as described in ~\ref{sec:fv} without additional
controls such as limiters or cap and clipping of the depth, whose
positivity follows from the scheme.
The nodal frames are computed from the exact gradient on the
analytical surfaces and from the 2-ring approximation of
Section~\ref{sec:nodalgeom} on data. Steady states are computed
exactly by one sparse solve.
Two kinds of meshes are used: \emph{structured} meshes, which split
every square of an $N_x\times N_y$ grid along the same diagonal, and
\emph{unstructured} meshes, which are generated by
Triangle~\citep{Shewchuk1996} with constrained Delaunay refinement.
Concerning boundary condition, tests~1 and~2 have one outlet and walls
that are parallel to the flow or divides, consistently with the
boundary rule of Section~\ref{sec:admissibility}, while tests~3 and~4
have free outflow on the whole boundary.

\paragraph{References.} For all test cases on an analytical surface we
defined some benchmark references to compare our numerical results.
Test~1 uses the exact event hydrograph of the plane for both the
intrinsic and the projected model. For Test~2, Appendix~\ref{app:l4}
gives the closed-form steady state, with the intrinsic depth
$\Depth=\kEta\DepthP$.
Moreover, Test~2 is also used to compare the edge velocities during
routing.  To measure routing separately from depth, the steady edge
fluxes are interpolated by the lowest-order Raviart-Thomas field,
whose value at the cell centroid is given by
$\Qvect_K=\sum_{E\subset\partial K}Q_{K,E}(\xv_K-\xv_E)/(2|K|)$, with
$\xv_E$ the vertex opposite $E$ and $\xv_K$ the centroid of the flat
cell. This vector field is compared with the exact discharge vector of
magnitude $\ABSsurf{\Qvect}=|\QvectP|=q_x|\Grad\BSM|/S$
(Appendix~\ref{app:l4}) and direction
$-\GradSurf\zbottom/\ABSsurf{\GradSurf\zbottom}$. Its steady state
divergence is exact up to machine precision for every edge velocity,
so this \emph{routing error} does not measure the mass balance of the
cells, but how the edge velocity spatially distributes the discharge,
i.e., its direction and, through the direction, its magnitude.
In Test~3 the steady reference is computed
along flow lines. Under constant rain, the steady unit discharge at a
point equals the rainfall rate times the upslope contributing horizontal
area per unit contour width.
Test~4 does not have any reference solution, it is used to show
applicability in real cases.

In addition, besides the four edge velocities of
Section~\ref{sec:sources} (\FVS, \AVG, \LNs, \LNw), the
\emph{exact-slope reference} computes the kinematic field considering
the exact edge means of the normal and tangential slopes, i.e.
$s^\nu_E=\edgeLength[E]^{-1}\int_E\GradSurf\zbottom\cdot\normalEdge_{K,E}\,ds$
and likewise $s^\tau_E$, computed by Gauss quadrature along the edge
curve with the exact conormal. It needs the analytical surface, and
its solutions are marked by the subscript $\mathrm{ex}$. It carries
the transport error of the scheme and no slope error: since the
transport error is common to all edge velocities, the difference
$\|\Depth_X-\Depth_{\mathrm{ex}}\|/\|\Depth_{\mathrm{ex}}\|$ isolates
the \emph{slope-induced error} of an edge velocity $X$ without a
higher-order scheme.

\paragraph{Metrics.} Steady errors are relative $L^1$ errors over the
valid cells, i.e., where the reference is valid, outside a small disc
around a degenerate stationary point, and outside traps. The time
metrics are the mean response time $\bar t=V_s/Q_s$, which does not
depend on the sampling of the hydrograph and is exact also for the
discrete scheme, the threshold times $t_\alpha$ of the rising limb
(first time with $Q\ge\alpha Q_p$), and the recession time $\tRec$
(time after the end of the rain at which $Q$ falls below $Q_p/2$),
with $Q_p$ the outflow at the end of the rain. Robustness is measured
by the number of traps and by the trapped rain, both computed from the
transfer directed graph without simulation.

\subsection{Test 1: rain event on a sloping plane}\label{sec:t1}

\begin{figure}
\centering
\includegraphics[width=0.49\textwidth]{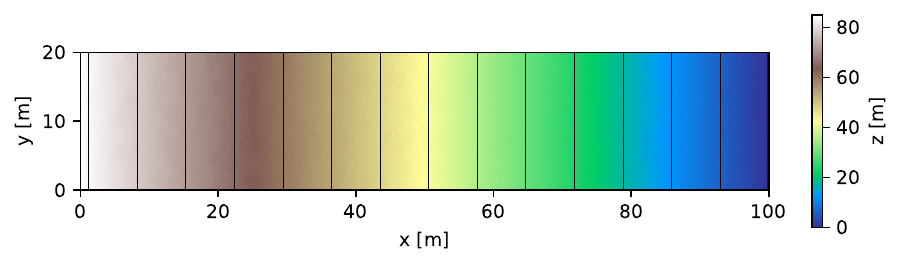}
\includegraphics[width=0.49\textwidth]{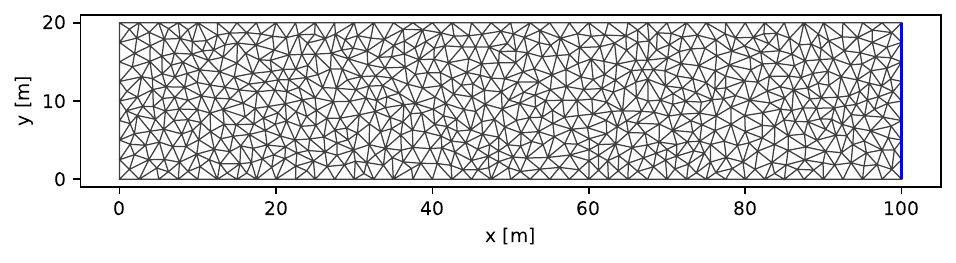}
\caption{Test 1: the plane with $\tan\AngleSymbol=0.85$ (contours, outlet in blue)
and the mesh.}\label{fig:t1setup}
\end{figure}

The domain of Test~1 is a sloping plane with height function defined
by $\BSM=S(L-x)$ on $[0,\pv{t1.L}]\times[0,\pv{t1.W}]$\,m. We simulate
and compare the results considering 
different values for the slope,
$S = \tan\AngleSymbol={}$\pv{t1.slopes}. We always assume an outlet at
$x=L=100$~m.
We consider an \pv{t1.kind} mesh with $\meshparam=\pv{t1.h}$\,m
(\pv{t1.nc} cells, Fig.~\ref{fig:t1setup}). The cell velocity field is spatially constant
on a sloping plane, and thus the four edge velocities give the same
flux on any mesh. Their steady states coincide up to \pv{t1.coincide},
and \LNw\ is used.

\paragraph{Exact solution.}
Intrinsic and projected kinematic wave models both have a closed-form
event hydrograph on the plane, which is provided in
Appendix~\ref{app:plane}.
By Remark~\ref{rem:timescale}, the intrinsic hydrograph is the
projected one stretched in time by $\kVol$, the rise from $t=0$ and,
from equilibrium, the recession from $t_r$. The exact steady storage
gives the exact $\bar t$ of each model.

\begin{figure}\centering
\includegraphics[width=0.6\textwidth]{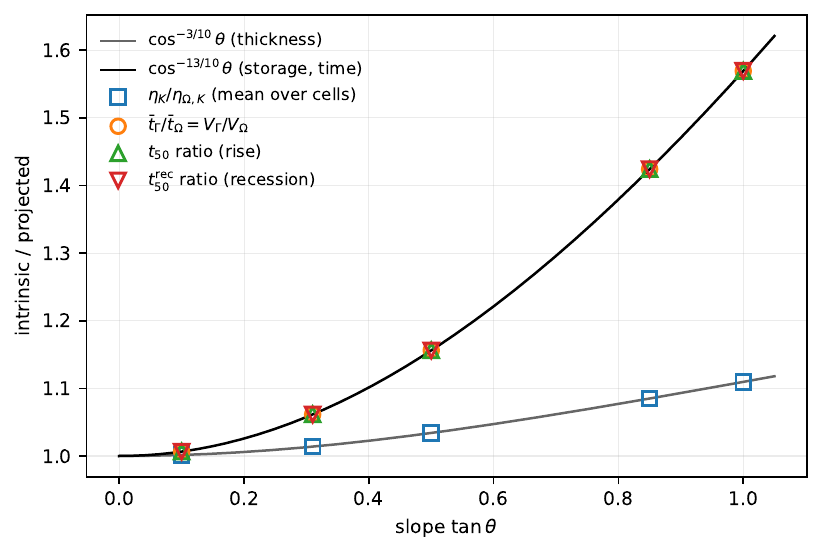}
\caption{Test 1: measured intrinsic/projected ratios against the slope, with
the closed forms $\kEta=\cos^{-3/10}\AngleSymbol$ and
$\kVol=\cos^{-13/10}\AngleSymbol$.}\label{fig:t1ratios}
\end{figure}

\begin{figure}\centering
\includegraphics[width=\textwidth]{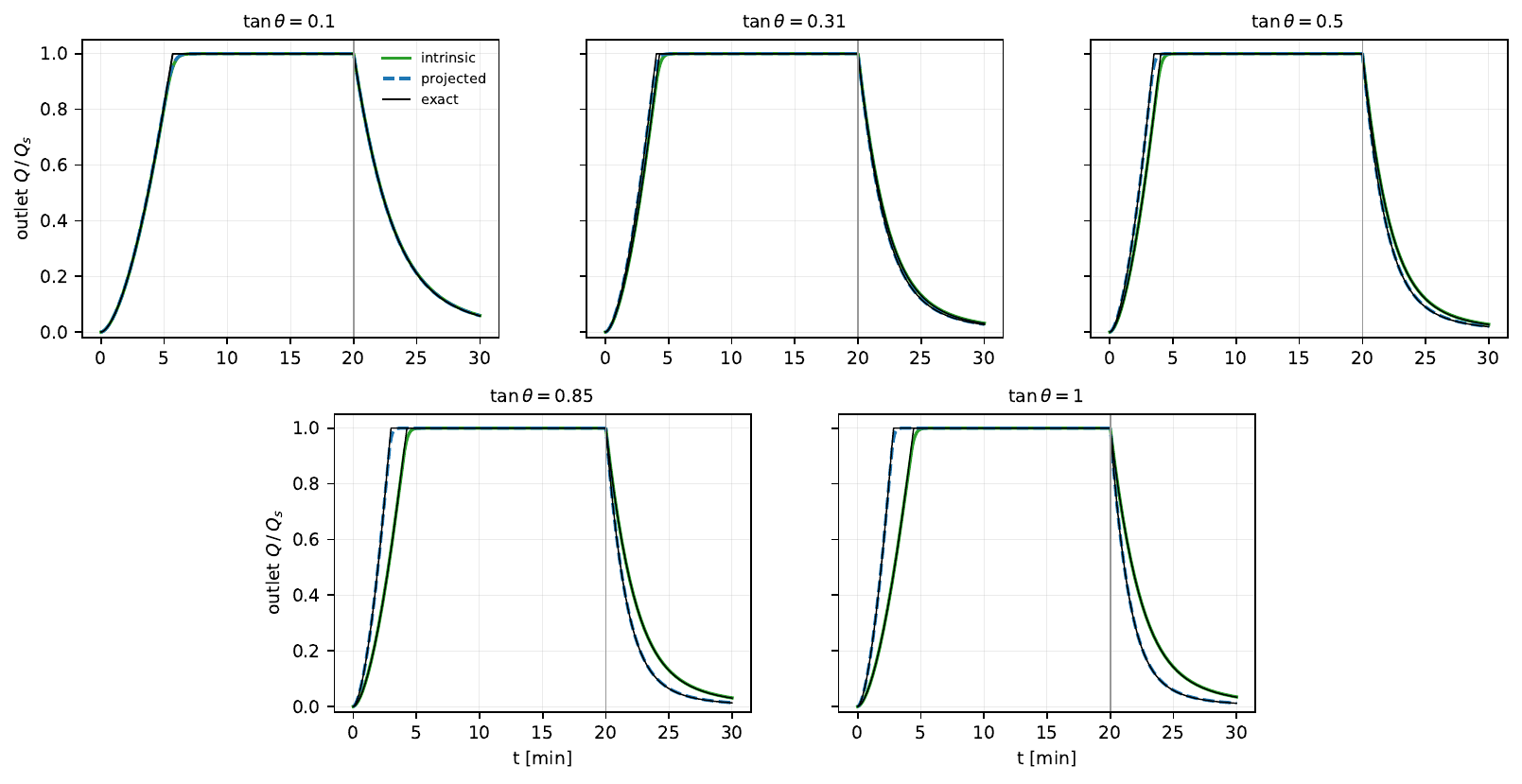}
\caption{Test 1: outlet hydrographs $Q/Q_s$, intrinsic (solid) and
  projected (dashed), with the exact hydrographs (thin black). The
  gray line marks the end of the rain.}\label{fig:t1hyd}
\end{figure}

\paragraph{Results.} The depth ratio equals $\kEta$ in every cell, and
the storage ratio, i.e.\ the ratio of the mean response times, equals
$\kVol$, both up to machine precision (see Fig.~\ref{fig:t1ratios}). On a
plane the scheme reproduces the stretch: with time step
$\kVol\,\Dt_{\max}$ in the intrinsic run, the two runs take the same
number of steps and the intrinsic hydrograph is the projected one
stretched by $\kVol$ to \pv{t1.matched.max}. 
The intrinsic and the projected hydrographs are within
\pv{t1.exacterr.range}\,$Q_s$ of the exact ones, and
$\bar t$ is off by \pv{t1.tbarerr.range} in both: this is the common
transport error of the first-order scheme, concentrated at the corner
of the rising limb. This is visible in Fig.~\ref{fig:t1hyd}.

The exact solution holds in both frameworks and against it the two
schemes have the same error, the transport error of the first-order
scheme, which vanishes under refinement. The difference between the
two \emph{models} does not: the projected model misses the intrinsic
mean response time by \pv{t1.modelgap.range} from
$\tan\AngleSymbol=0.1$ to $1$. On slopes above about $0.3$ the model
error is larger than the discretization error of the scheme on this
mesh, and it is the only error that refinement does not remove: where
the intrinsic model is the consistent one (thin layer along the normal
on a steep bed), the projected model is not.

\subsection{Test 2: convergent valley, structured and unstructured meshes}\label{sec:t2}

\begin{figure}\centering
\includegraphics[width=0.6\textwidth]{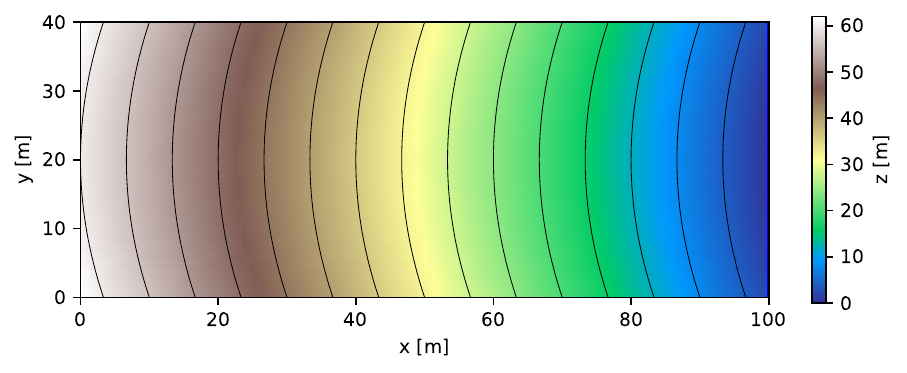}

\includegraphics[width=0.49\textwidth]{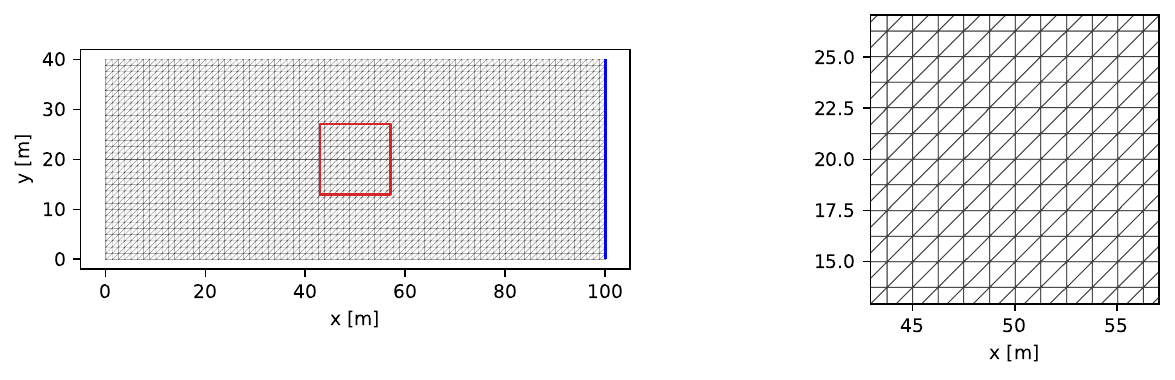}\hfill
\includegraphics[width=0.49\textwidth]{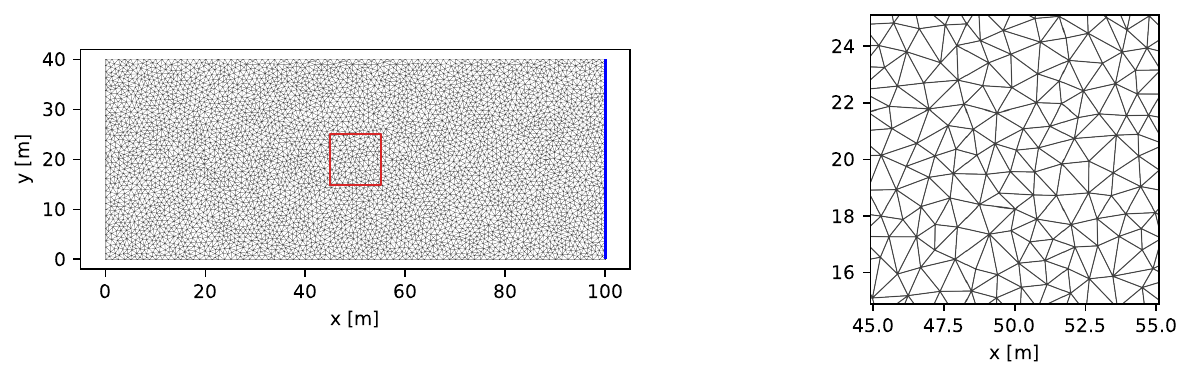}
\caption{Test 2: the valley (contours, outlet in blue) and the meshes with
$\meshparam=\pv{t2.hrun}$\,m, structured (left panel) and unstructured
(right panel), with a
zoom.}\label{fig:t2mesh}
\end{figure}

We consider here the surface given by the height function:
\begin{equation*}
  \BSM=S(L-x)+a(y-W/2)^2,
\end{equation*}
with $S=\pv{t2.S}$ and $a={}$\pv{t2.a}, on the domain
$[0,\pv{t2.L}]\times[0,\pv{t2.W}]$\,m. The outlet is located at $x=L=100$~m,
while the walls $y=0$ and $y=W=\pv{t2.W}$\,m are divides.
Its elevation contour lines along the y direction are parabolas and
the flow converges on the valley axis $y=W/2=20$\,m. With
$\mathbf m=\Grad\BSM/|\Grad\BSM|$, the two curvatures of
\eqref{eq:divfields} are $k_h=2a\,m_x^2\cos\AngleSymbol$ and
$k_v=2a\,m_y^2\cos^3\AngleSymbol$: $k_v$ vanishes on the axis and is
at most $6\cdot10^{-4}$\,m$^{-1}$ at the walls, more than ten times
smaller than $k_h$ ($7.6$--$8.6\cdot10^{-3}$\,m$^{-1}$), so the
convergence of the flow lines dominates.
Two mesh families with $\meshparam=\pv{t2.hlist}$\,m are used
(Fig.~\ref{fig:t2mesh}). One is structured, with the valley axis
aligned with mesh edges on every level, and the other
  fully unstructured. On the structured mesh the two cells on the
sides of an axis edge have velocities pointing into each other, and
thus the axis is a line of counterflow edges (\pv{t2.stru.cf} on the
finest mesh). The rain event is run on $\meshparam=\pv{t2.hrun}$\,m.

\paragraph{Exact solution.} The steady state has the closed form of
Appendix~\ref{app:l4}, and the intrinsic depth is
$\Depth=\kEta\DepthP$. The exact mean response times are
$\bar t=\pv{t2.ex.tbar.int}$\,s and \pv{t2.ex.tbar.proj}\,s,
respectively for the intrinsic and projected model, with ratio
\pv{t2.ex.tbar.ratio}.

\begin{figure}\centering
\includegraphics[width=0.85\textwidth]{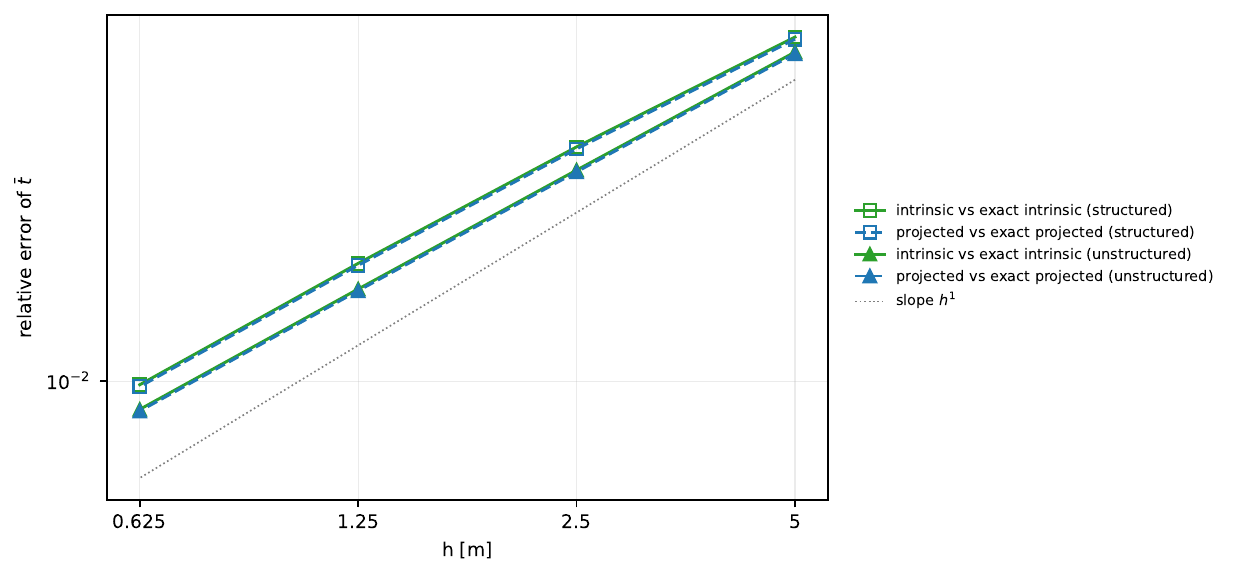}
\caption{Test 2: relative error of $\bar t$ against $\meshparam$. Green and blue: each
geometry against its own exact value (discretization error). Squares:
structured, filled triangles: unstructured.}\label{fig:t2tbar}
\end{figure}

\begin{table}\centering
  \caption{Intrinsic versus projected model, with \LNw: ratio of the
    mean response times $\bar t_{\SurfDomain}/\bar t_{\Domain}$
    against the storage-weighted and the area means of the storage
    factor $\kVol=\cos^{-13/10}\AngleSymbol_K$. Test~2: unstructured
    mesh, $\meshparam=\pv{t2.hrun}$\,m; Test~3:
    $\meshparam=\pv{t3.proj.h}$\,m; Test~4: data
    mesh.}\label{tab:intproj}
\begin{tabular}{@{}lccc@{}}\toprule
  & Test 2 & Test 3 & Test 4\\\midrule
  slope $\tan\AngleSymbol$, mean / max & 0.61 / 0.63 
           & \pv{t3.slopemean} / \pv{t3.slopemax} & \pv{t4.slopemean} / \pv{t4.slopemax}\\
  $\bar t_{\SurfDomain}/\bar t_{\Domain}$, exact & \pv{t2.ex.tbar.ratio} & -- & --\\
  $\bar t_{\SurfDomain}/\bar t_{\Domain}$, computed & \pv{t2.proj.tria.ratio} & \pv{t3.proj.ratio} & \pv{t4.proj.ratio}\\
  storage-weighted mean of $\kVol$ & \pv{t2.proj.tria.storage} & \pv{t3.proj.storage} & \pv{t4.proj.storage}\\
  area mean of $\kVol$ & \pv{t2.proj.tria.area} & \pv{t3.proj.area} & \pv{t4.proj.area}\\
  \bottomrule
\end{tabular}
\end{table}

\paragraph{Intrinsic versus projected.}
As in Test~1, we first compare the two models and confirm also on this
curved surface that the model error, which does not depend on the
mesh, dominates the discretization error.  On a curved surface there
is no single time stretch: the ratio of the mean response times is the
storage-weighted mean of $\kVol$ (Section~\ref{sec:equiv}), and the
computed ratio follows it in every test
(Table~\ref{tab:intproj}). Test~2 is the curved case with an exact
solution, and it separates the two errors. The computed mean response
times converge to their own exact values at first order, with the same
error in the two geometries (\pv{t2.tbar.tria.errint.list} on the
unstructured meshes, similar on the structured ones, see
Fig.~\ref{fig:t2tbar}), and this is the transport error of the
scheme. The model error is different in kind: the exact projected
$\bar t$ differs by \pv{t2.tbar.gap} from the exact intrinsic one,
about \pv{t2.tbar.tria.gapovererr} times the discretization error on
the finest unstructured mesh, and it does not depend on the mesh.

\begin{figure}\centering
\includegraphics[width=0.49\textwidth]{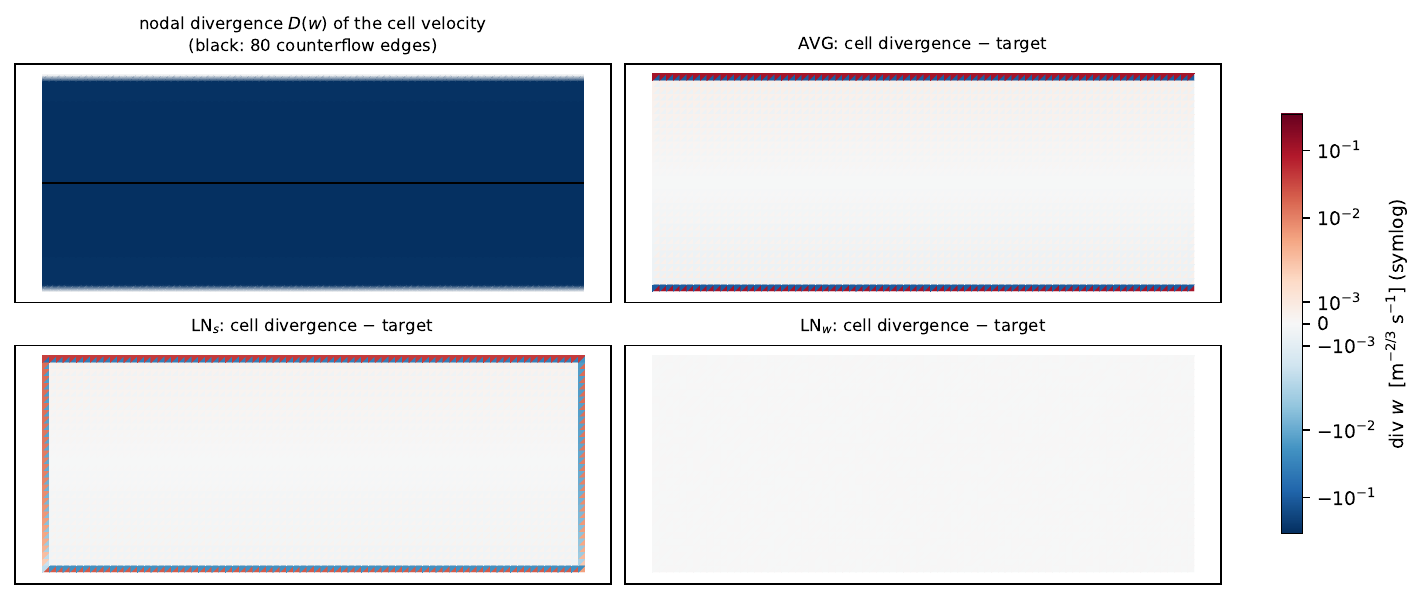}\hfill
\includegraphics[width=0.49\textwidth]{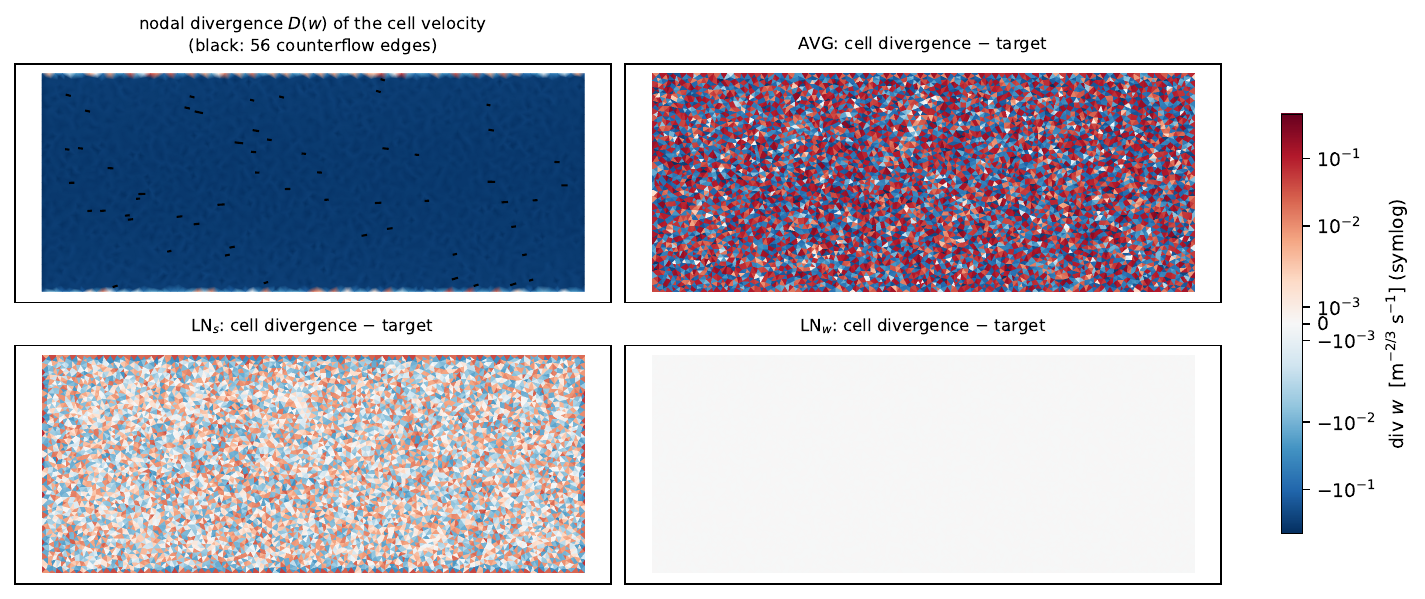}
\caption{Test 2: nodal divergence $D(\AdvVel)$ of the raw field (red:
  divergent, blue: convergent) with the counterflow edges (black), and
  the cell-divergence defects of \AVG, \LNs, \LNw, on one symmetric
  logarithmic scale. The defect is the cell divergence of the edge
  velocity, i.e.\ its net outflow per unit cell area, minus the target
  of Corollary~\ref{cor:divform}, i.e. the cell average of
  $D(\AdvVel)$ plus $d_K/|K|$ (\LNw\ matches it by
  construction). Structured (left panel) and unstructured (right
  panel) mesh.}\label{fig:t2div}
\end{figure}

\begin{figure}\centering
\includegraphics[width=0.49\textwidth]{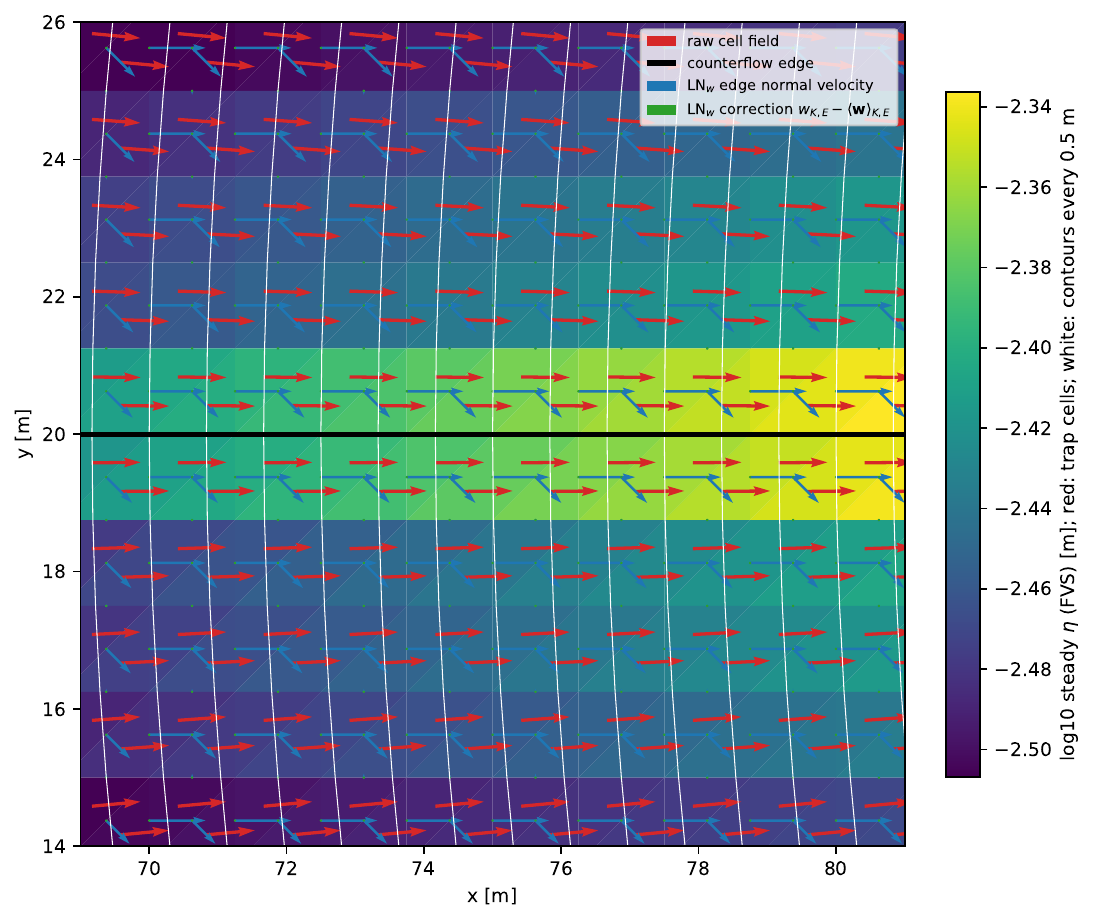}\hfill
\includegraphics[width=0.49\textwidth]{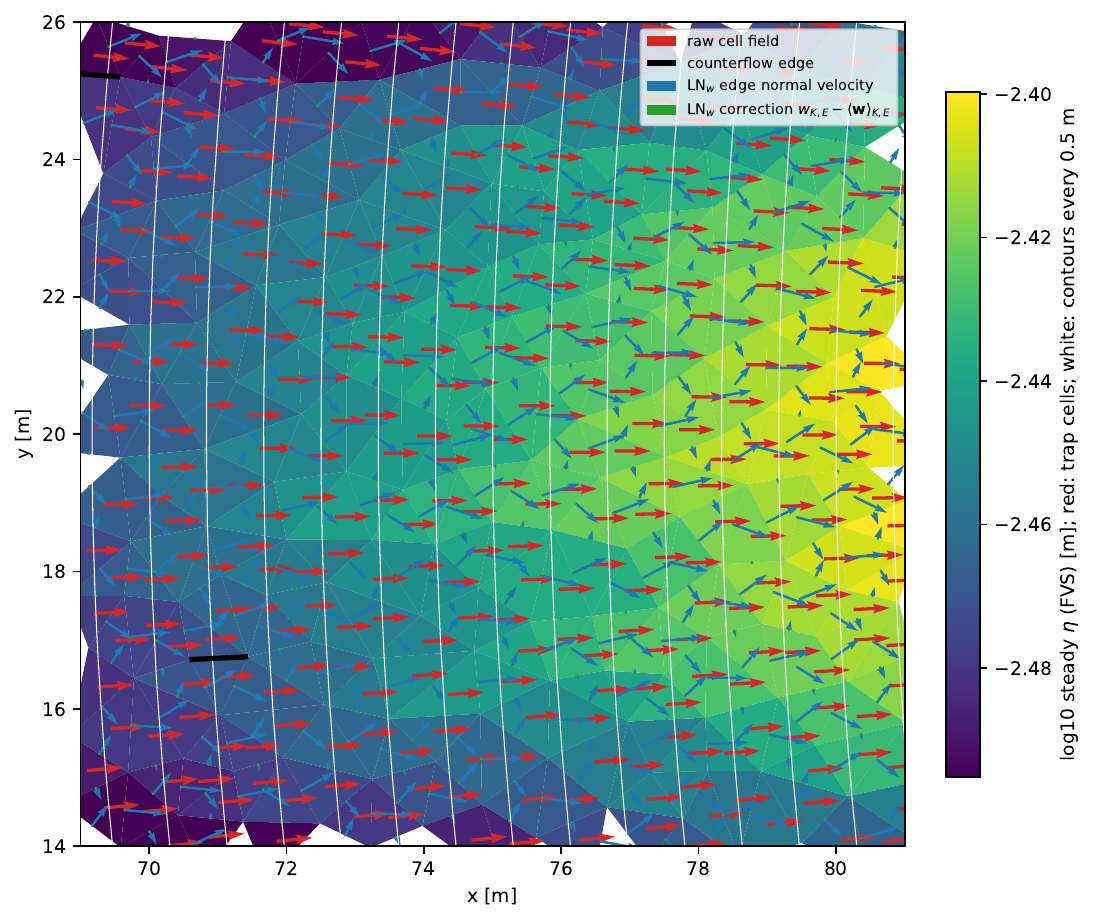}
\caption{Test 2: zoom on the valley axis, structured (left panel) and
  unstructured
(right panel) mesh. Color: steady depth of \FVS. Red arrows: raw cell
field. Black: counterflow edges. Blue: \LNw\ edge normal velocities.
Green: the \LNw\ correction
$w_{K,E}-\bar g_{K,E}(\AdvVel)/\edgeLength[E]$.}\label{fig:t2zoom}
\end{figure}

\begin{figure}\centering
\includegraphics[width=0.49\textwidth]{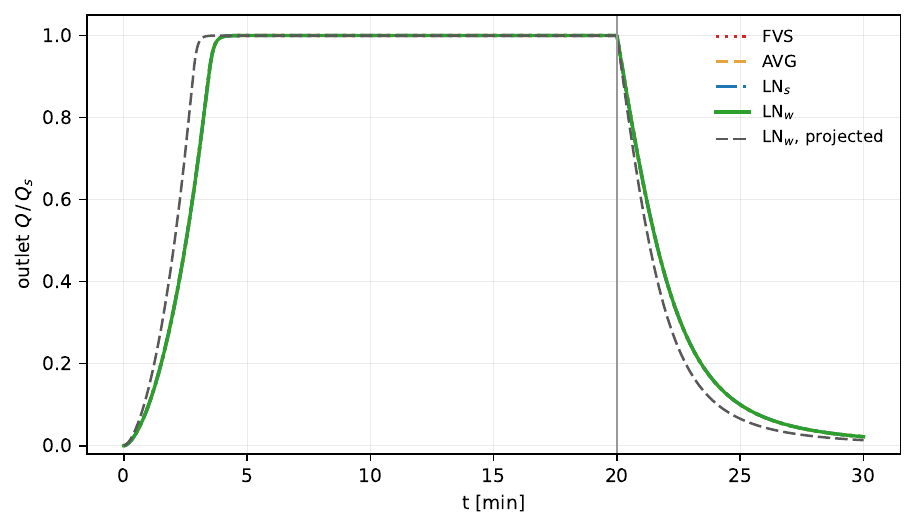}
\caption{Test 2: outlet hydrographs $Q/Q_s$ on the unstructured
mesh. The gray dashed line is with \LNw\ in the projected model.}\label{fig:t2hyd}
\end{figure}

\begin{figure}\centering
\includegraphics[width=\textwidth]{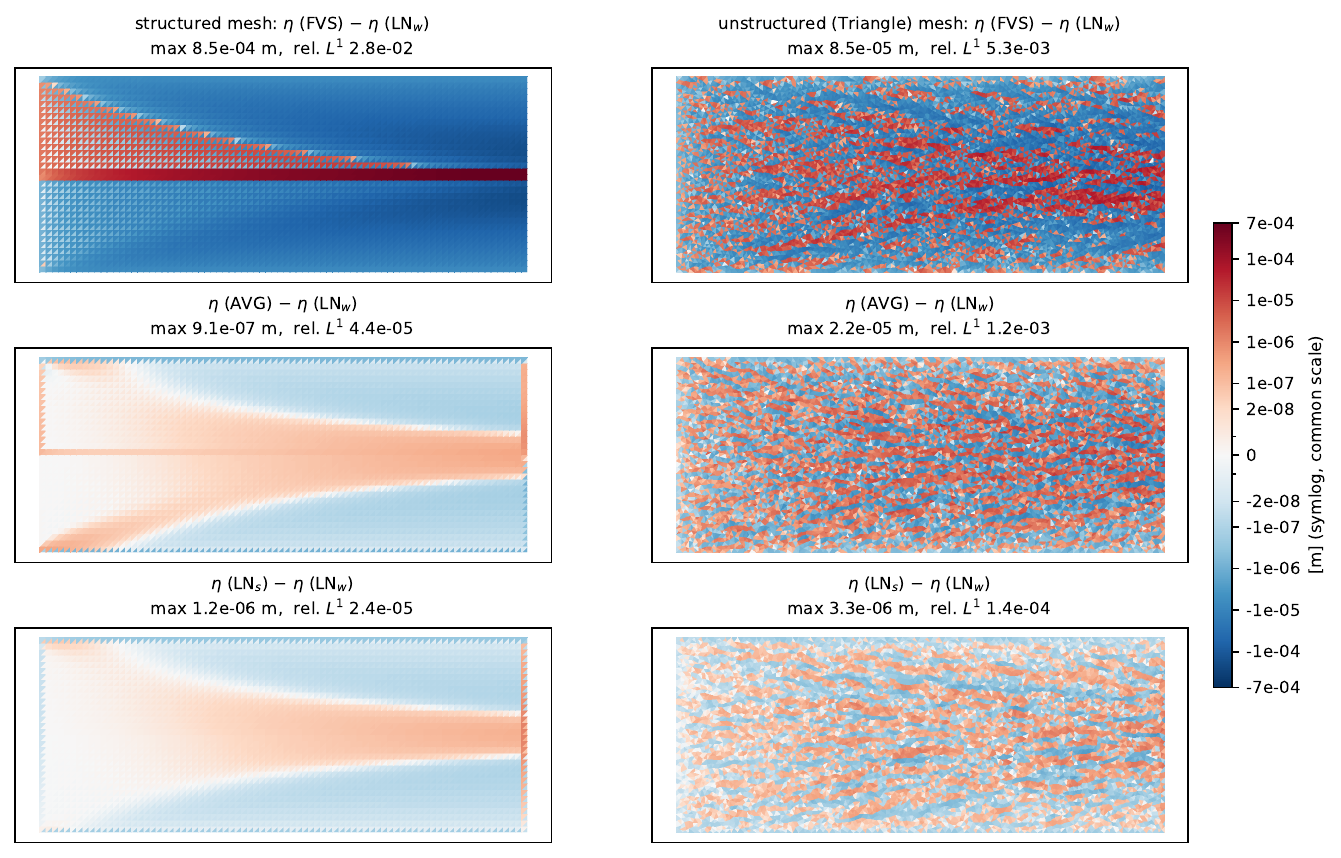}
\caption{Test 2: depth differences $\Depth_X-\Depth_{\LNw}$, $X\in\{\FVS,\AVG,\LNs\}$, at the
end of the rain, structured (left panel) and unstructured (right panel) mesh, on one
common symmetric logarithmic color scale.}\label{fig:t2diff}
\end{figure}

\paragraph{Results.}
In this case, on the structured mesh the counterflow edges of the raw
field are exactly the valley axis, while on the unstructured mesh they
are scattered over the domain (Figs.~\ref{fig:t2div}
and~\ref{fig:t2zoom}).  No edge velocity has a trap on any mesh of
either family.
On the run mesh the outlet response does not separate the edge
velocities: the mean response times lie in
\pv{t2.tria.run.tbar.range}\,s on the unstructured mesh and
\pv{t2.stru.run.tbar.range}\,s on the structured one, the threshold
times are $t_{10}=\pv{t2.tria.run.t10}$, $t_{50}=\pv{t2.tria.run.t50}$
and $t_{90}=\pv{t2.tria.run.t90}$\,s, with a spread of the rise times
of at most \pv{t2.stru.run.trise.spread}\,s, and the recession time is
$\tRec=\pv{t2.tria.run.trec}$\,s (Fig.~\ref{fig:t2hyd}). The depth
differences sit on the axis on the structured mesh and are scattered
on the unstructured one (Fig.~\ref{fig:t2diff}). The mesh study below
isolates the slope-induced error, which the outlet response hides.

\begin{figure}\centering
\includegraphics[width=\textwidth]{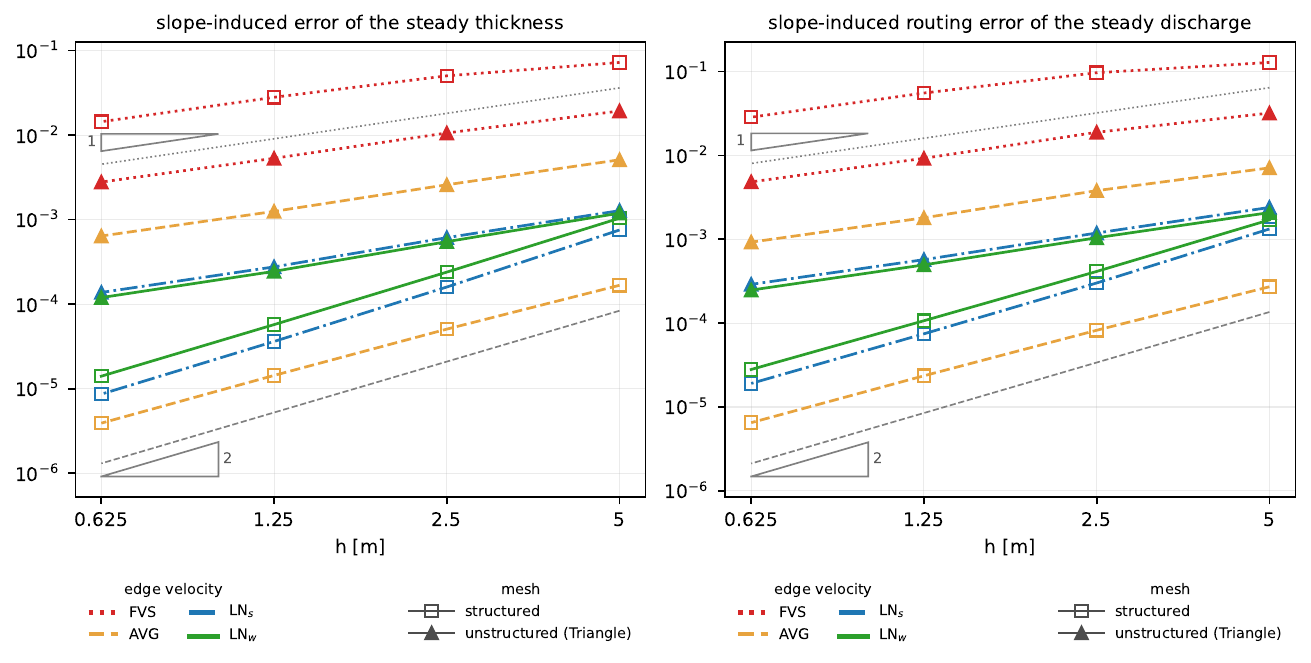}
\caption{Test 2: slope-induced error against $\meshparam$ of the
  steady depth (left panel) and of the routing of the steady discharge
  (right panel), on the structured (squares) and unstructured (filled
  triangles) meshes. Dotted line: slope $\meshparam^1$. Dashed line:
  slope $\meshparam^2$. Note that, as discussed in
  Remark~\ref{rem:raster}, a superconvergent behavior was expected on
  the structured mesh.}\label{fig:t2meshdep}
\end{figure}

We also test for convergence. On the unstructured meshes the
corrections reduce the slope-induced error of the plain average by a
factor \pv{t2.tria.avgoverlnw} on the finest level (\AVG\
\pv{t2.tria.avg.slopeerr}, \LNs\ \pv{t2.tria.lns.slopeerr}, \LNw\
\pv{t2.tria.lnw.slopeerr}), with the same first-order rate, so that
the ratio is stable under refinement (Fig.~\ref{fig:t2meshdep}, left
panel).
On the structured meshes the plain average is superconvergent (order
\pv{t2.stru.avg.slopeerr.rate}), all single-valued edge velocities
converge at second order with small constants, and the correction does
not improve on the average. As expected, \FVS\ has the largest
slope-induced error on both families (\pv{t2.stru.fvs.slopeerr}
structured, \pv{t2.tria.fvs.slopeerr} unstructured): on the structured
mesh it exchanges water across the axis edges, and its error against
the exact solution on the axis band is \pv{t2.stru.fvs.err.axis},
against \pv{t2.stru.lnw.err.axis} for \LNw\ and
\pv{t2.stru.avg.err.axis} for \AVG\ whereas off the band all are close
to \pv{t2.stru.lnw.err.off} (Fig.~\ref{fig:t2zoom}).
The transport error (\pv{t2.stru.transport} structured,
\pv{t2.tria.transport} unstructured, observed order
\pv{t2.stru.transport.rate} and \pv{t2.tria.transport.rate}) is at
least one order of magnitude larger than every slope-induced error of
the single-valued edge velocities, and the total error of \LNw\
coincides with it to three digits. The routing error of the steady
discharge gives the same ranking as the depth error, with a smaller
factor: against the exact-slope reference, the slope-induced routing
error of \AVG\ is \pv{t2.tria.route.avgoverlnw} times that of \LNw\ on
the unstructured meshes (\pv{t2.tria.avgoverlnw} for the depth), and
\pv{t2.stru.route.avgoverlnw} times on the structured ones
(\pv{t2.stru.avgoverlnw} for the depth). See Fig.~\ref{fig:t2meshdep},
right panel, for the routing error convergence.

\subsection{Test 3: a degenerate stationary point}\label{sec:t3}\label{sec:ex1}

\begin{figure}\centering
\begin{minipage}[t]{0.49\textwidth}\centering
\includegraphics[width=\linewidth]{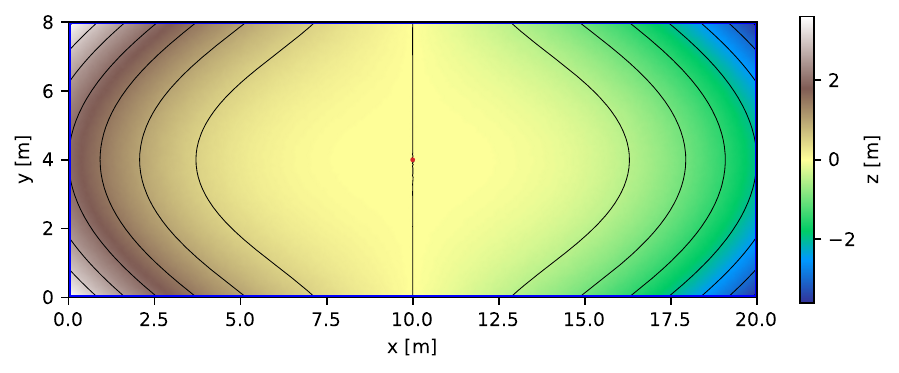}
\end{minipage}\hfill
\begin{minipage}[t]{0.49\textwidth}\centering
\includegraphics[width=\linewidth]{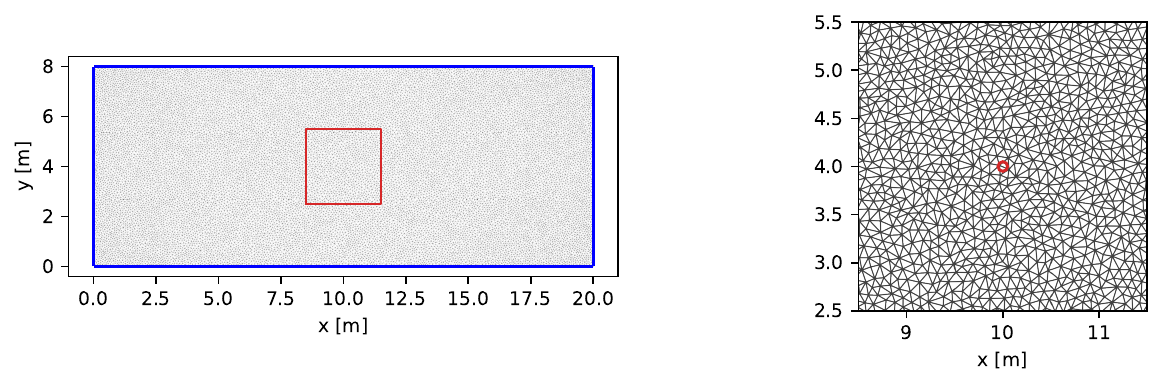}
\end{minipage}
\caption{Test 3: the surface (contours, free outflow everywhere) and
  the mesh $\meshparam=\pv{t3.hrun}$\,m with a zoom on the stationary
  point. In red, the disc of radius \pv{t3.refdisc}\,m excluded from
  the reference.}\label{fig:t3setup}\label{fig:ex1terrain}
\end{figure}

Test~3 considers the surface  used in
\citet[Sec.~4.4]{BachiniPutti2020}, given by the height function:
\begin{equation*}
  \BSM=-\frac{1}{500}(x-10)^3-\frac{1}{100}(x-10)(y-4)^2
\end{equation*}
on $[0,20]\times[0,8]$\,m, with free outflow on the whole boundary.
For this surface, the slopes are up to \pv{t3.slopemax} (mean slope
\pv{t3.slopemean}, with $\AngleSymbol\le\pv{t3.thetamax}^\circ$). We
consider unstructured meshes with $\meshparam=\pv{t3.hlist}$\,m
(\pv{t3.nclist} cells, respectively. See Fig.~\ref{fig:t3setup}), and
the rain event is run on $\meshparam=\pv{t3.hrun}$\,m.  The reference
is the exact steady flow-line solution on the valley half
(\pv{t3.refarea} of the area), outside a disc of radius
\pv{t3.refdisc}\,m around the stationary point.

The surface is a cross-slope valley for $x<10$ and a ridge for $x>10$,
separated by a straight contour. Its gradient,
\begin{equation*}
\Der_x\BSM=-\tfrac{3}{500}(x-10)^2-\tfrac{1}{100}(y-4)^2, \qquad
\Der_y\BSM=-\tfrac{2}{100}(x-10)(y-4),
\end{equation*}
is homogeneous of degree two in $(x-10,y-4)$ and vanishes only at
$(10,4)$, where all second derivatives vanish too. Flow lines with
$|y-4|\ll|x-10|$ behave like $|y-4|\propto|x-10|^{10/3}$: they reach
the stationary point tangentially from the valley side and leave it
tangentially on the ridge side. Tracing the flow lines shows that
$83\%$ of the valley half drains into the stationary point and $83\%$
of the ridge half is fed from it, so that the point carries $41.6\%$
of the total rain. Downstream of it the steady problem is ill
posed. Away from the point the surface is admissible: accuracy is
therefore measured on the valley half, where every point is regular,
and the ridge half serves as a robustness test.

\paragraph{Intrinsic versus projected.} The comparison of the two
models confirms the observations of Tests~1 and~2. Here the
storage-weighted and the area means of $\kVol$ differ most
(Table~\ref{tab:intproj}): storage concentrates where the slope is
gentle, and the computed ratio follows the storage-weighted mean to
\pv{t3.proj.agree} for the single-valued edge velocities.

\begin{figure}\centering
\includegraphics[width=\textwidth]{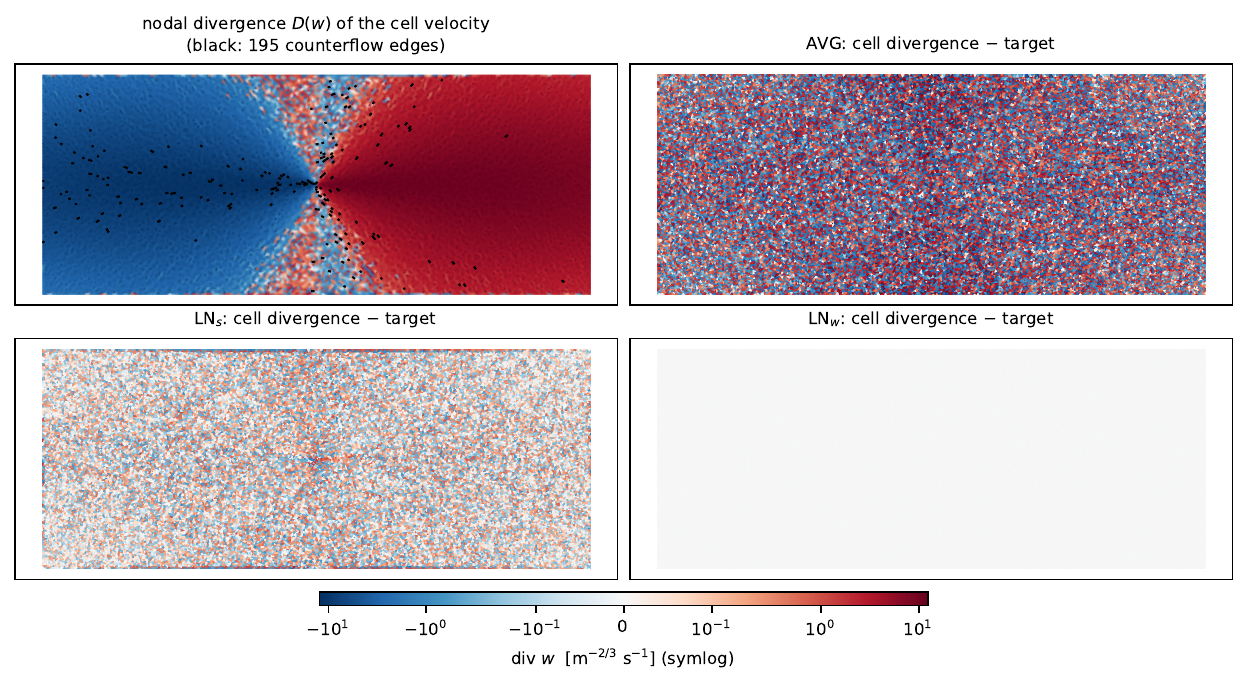}
\caption{Test 3, $\meshparam=\pv{t3.hrun}$\,m: nodal divergence
  $D(\AdvVel)$ of the raw field (red: divergent, blue: convergent)
  with the counterflow edges (black), and the cell-divergence defects
  of \AVG, \LNs, \LNw, on one symmetric logarithmic
  scale.}\label{fig:t3div}\label{fig:ex1div}
\end{figure}

\begin{figure}\centering
  \includegraphics[width=0.49\textwidth]{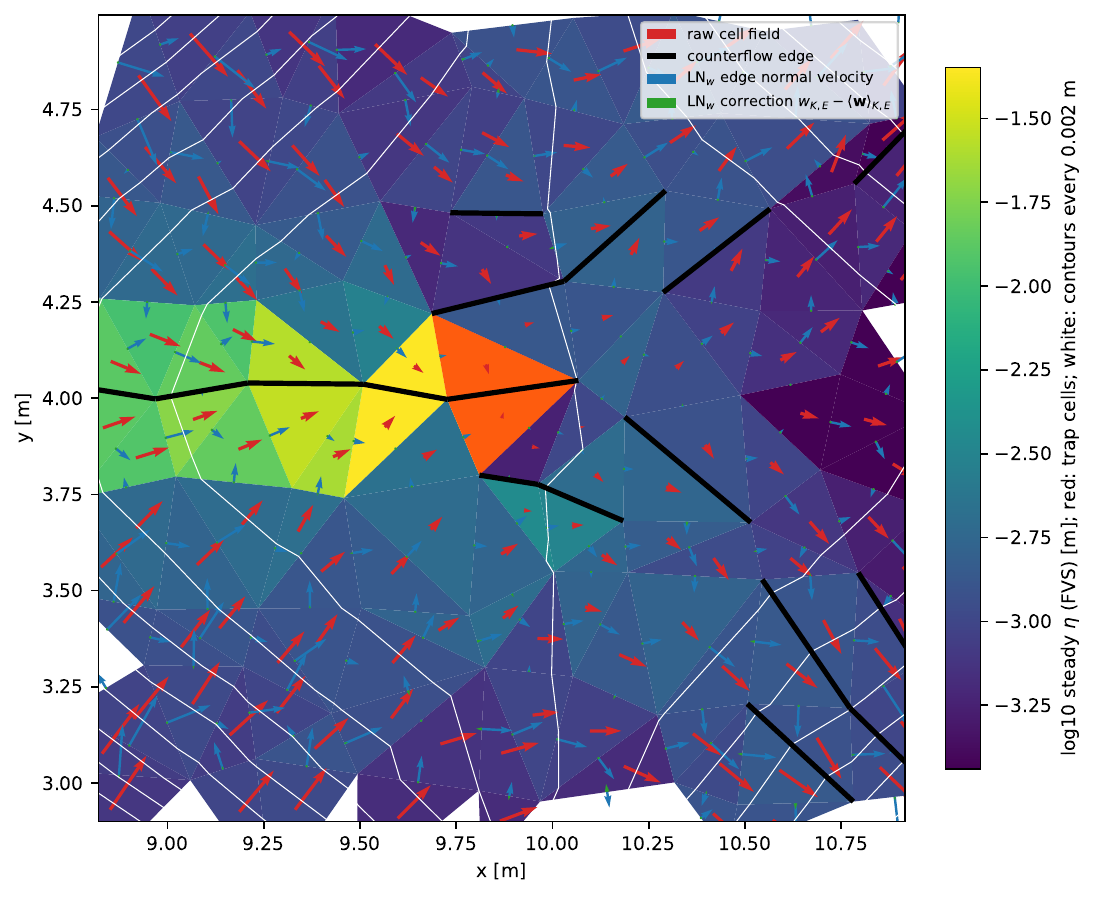}
  \hfill
  \includegraphics[width=0.49\textwidth]{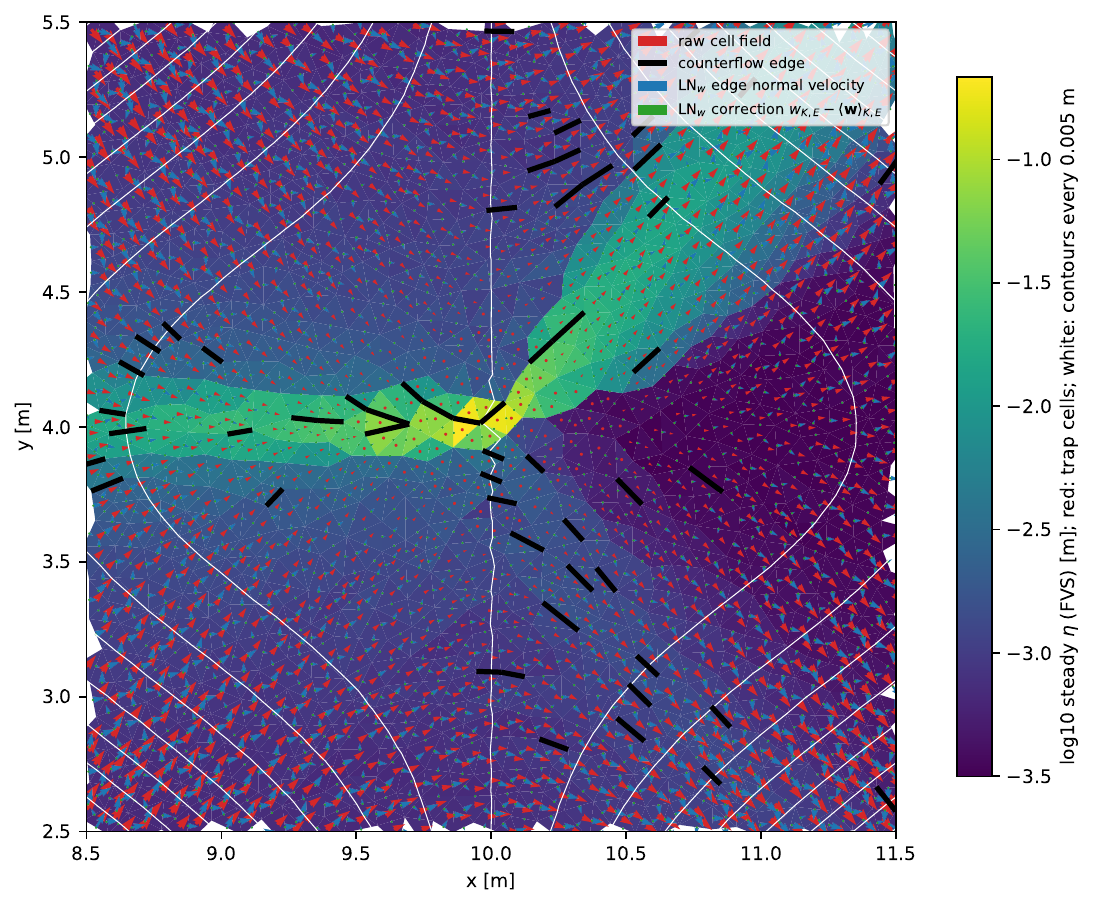}
  \caption{Test 3: zoom at the stationary point of the \FVS\ trap on
    the unstructured meshes with
    $\meshparam=\pv{t3.fam.pert.fvs.trap.hlist}$\,m (left panel), and
    $\meshparam=\pv{t3.hrun}$\,m (right panel). {Color: steady depth
      of \FVS; red arrows: raw cell field; black: counterflow edges;
      blue: \LNw\ edge normal velocities; green: the \LNw\
      correction.}}\label{fig:t3trap}\label{fig:ex1trap}
\end{figure}

\paragraph{Results.}
In this case, the normal components of the raw kinematic field jump
across the edges with a relative jump of \pv{t3.div.jump} on the
finest mesh that decreases like $\meshparam$. The jumps carry the
divergence of the field, and where they are convergent the two cells
push water into each other: \pv{t3.div.cf} edges are counterflow
edges, on the valley axis and around the point
(Fig.~\ref{fig:t3div}). On the unstructured meshes no trap forms with
\AVG, \LNs\ and \LNw. The exact-slope reference, instead, traps on the
unstructured mesh with $\meshparam=\pv{t3.exact.trap.hlist}$\,m, where
it keeps \pv{t3.exact.trap.rain} of the rain: exact slopes reproduce
the vanishing velocity of the continuum at the stationary point, where
the flow of the valley half converges, and the cells around it lose
their outflow.
A single-valued edge velocity is therefore not trap free by itself. On
the meshes where it does not trap, the exact-slope reference has a
steady error comparable with that of the other edge velocities
(\pv{t3.exact.err.all.list} on all valid cells,
\pv{t3.exact.err.off.list} off the axis, for
$\meshparam=\pv{t3.exact.notrap.hlist}$\,m respectively).
Fig.~\ref{fig:t3trap} shows instead the flux-vector splitting traps on
the mesh with $\meshparam=\pv{t3.fam.pert.fvs.trap.hlist}$\,m, where a
group of cells connected by counterflow edges around the point closes
into a trap that keeps \pv{t3.fam.pert.fvs.trap.rain} of the rain, the
exact-slope reference on
$\meshparam=\pv{t3.fam.pert.exact.trap.hlist}$\,m
(\pv{t3.fam.pert.exact.trap.rain}), and the single-valued edge
velocities on none. The mesh is admissible everywhere else, yet a
single degenerate point is enough to capture the water. The trap is
not visible in the vertex heights, but the analysis of the transfer
directed graph finds it before any simulation, as discussed in
Section~\ref{sec:traps}.

Against the reference, the steady error is dominated by the valley
axis, where the flow concentrates: off the axis ($|y-4|>0.5$\,m) it
converges at first order for every edge velocity, to about
\pv{t3.lnw.err.off} on the finest mesh for \LNw, while on all valid
cells it is not yet asymptotic at these resolutions, from
\pv{t3.err.all.coarse.range} on the coarsest mesh to
\pv{t3.err.all.fine.range} on the finest, for every edge velocity. The
exact-slope reference has comparable errors where it does not trap, so
these are transport errors on the axis, not slope errors. In this
regard, Test~3 is used for robustness for the imposed divergence, not
to rank the edge velocities by accuracy.

The correction does what Corollary~\ref{cor:divform} states: \LNw\
reproduces in every cell the target divergence of $\AdvVel$ to
\pv{t3.div.lnw.defect}, and \LNs\ that of the slope to roundoff.  The
plain average misses the target by \pv{t3.div.avg.defect}, and this
defect does not decrease under refinement.
The nodal divergence of the \LNw\ edge field, i.e., its cell
divergences averaged back to the vertices with lumped weights, differs
from the raw nodal divergence $D(\AdvVel)$ of \eqref{eq:Dnodal} by
\pv{t3.nod.w.b.lnw}. Against the exact divergence this averaging helps
pointwise, from \pv{t3.nod.w.c.D} for $D(\AdvVel)$ to
\pv{t3.nod.w.c.lnw} for \LNw\ (\pv{t3.nod.w.c.avg} for \AVG; interior
vertices outside the disc), and all fields agree in block integrals
($2$\,m blocks: \pv{t3.nod.w.d.D} for $D(\AdvVel)$,
\pv{t3.nod.w.d.lnw} for \LNw, \pv{t3.nod.w.d.avg} for \AVG).  The
defect of \AVG\ is a defect with respect to the imposed target, not an
error of order one against the exact divergence: pointwise \AVG\ is
less accurate than \LNw, in block averages the two are equally
accurate (Fig.~\ref{fig:t3div}).

\begin{figure}\centering
\includegraphics[width=0.85\textwidth]{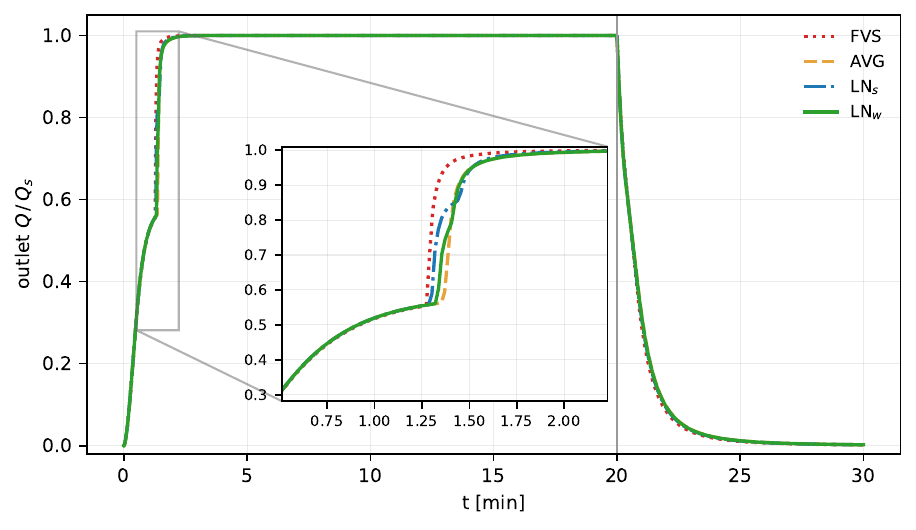}
\caption{Test 3: outlet hydrographs $Q/Q_s$ on
$\meshparam=\pv{t3.hrun}$\,m (inset: zoom on the rise to the plateau,
$t=\pv{t3.inset.t0}$-\pv{t3.inset.t1}\,s).}\label{fig:t3hyd}\label{fig:ex1hyd}
\end{figure}

\begin{figure}\centering
\includegraphics[width=\textwidth]{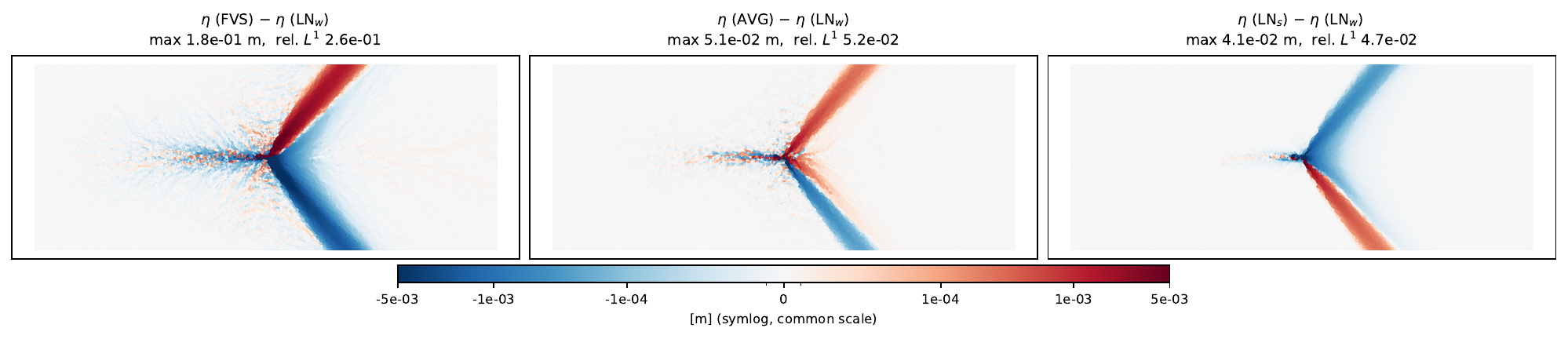}
\caption{Test 3: depth differences $\Depth_X-\Depth_{\LNw}$, $X\in\{\FVS,\AVG,\LNs\}$, at the
end of the rain, on one common symmetric logarithmic color
scale.}\label{fig:t3diff}
\end{figure}

For the rain event on $\meshparam=\pv{t3.hrun}$\,m the outlet response
does not discriminate between the single-valued edge velocities: their
mean response times differ by \pv{t3.tmean.spread}
($\bar t=\pv{t3.lnw.tmean}$\,s for \LNw) and the hydrographs nearly
overlap (Fig.~\ref{fig:t3hyd}), but the differences sit on the valley
axis and downstream of the point (Fig.~\ref{fig:t3diff}). \FVS\ drains
faster ($\bar t=\pv{t3.fvs.tmean}$\,s) and its steady depth differs
from that of \LNw\ by \pv{t3.fvs.vslnw}, against \pv{t3.avg.vslnw} for
\AVG\ and \pv{t3.lns.vslnw} for \LNs. Mass is conserved to
\pv{t3.mass.max} in all cases.

\subsection{Test 4: a catchment}\label{sec:t4}\label{sec:ex2}

\begin{figure}\centering
\includegraphics[height=0.45\textheight]{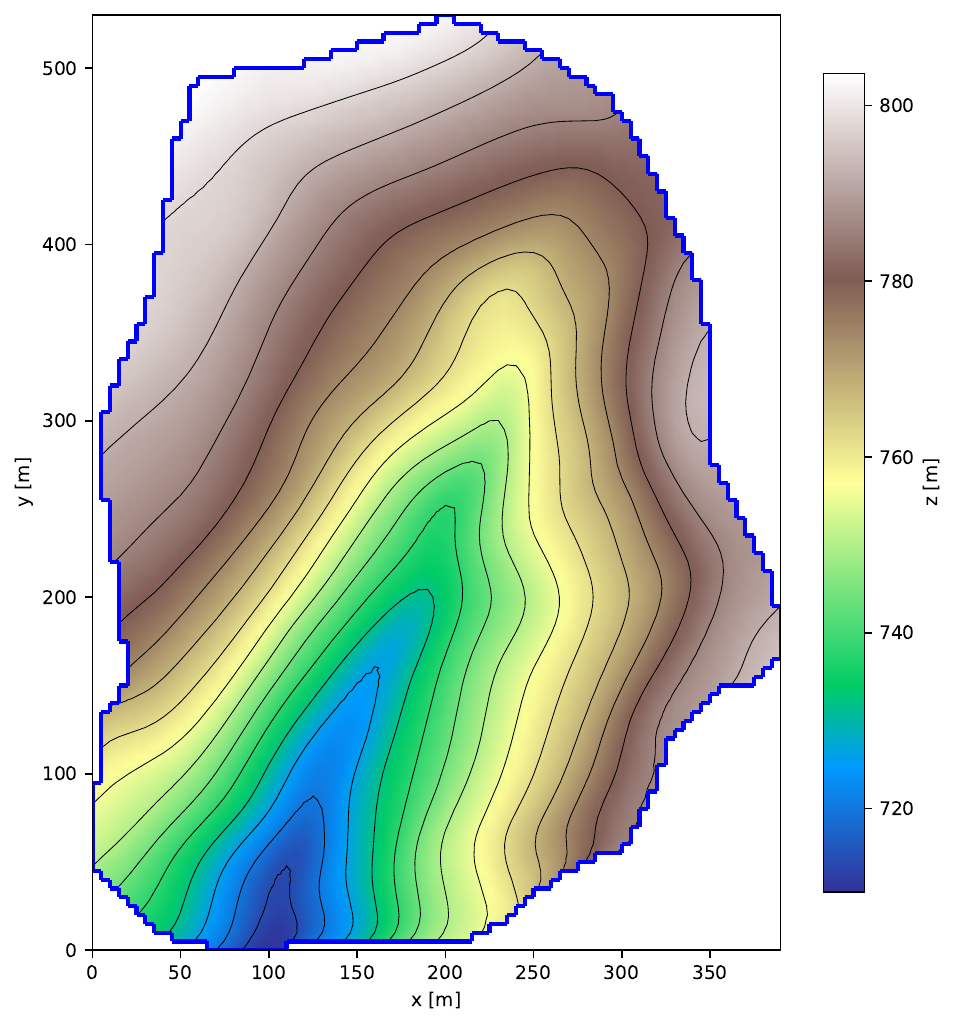}
\includegraphics[height=0.45\textheight]{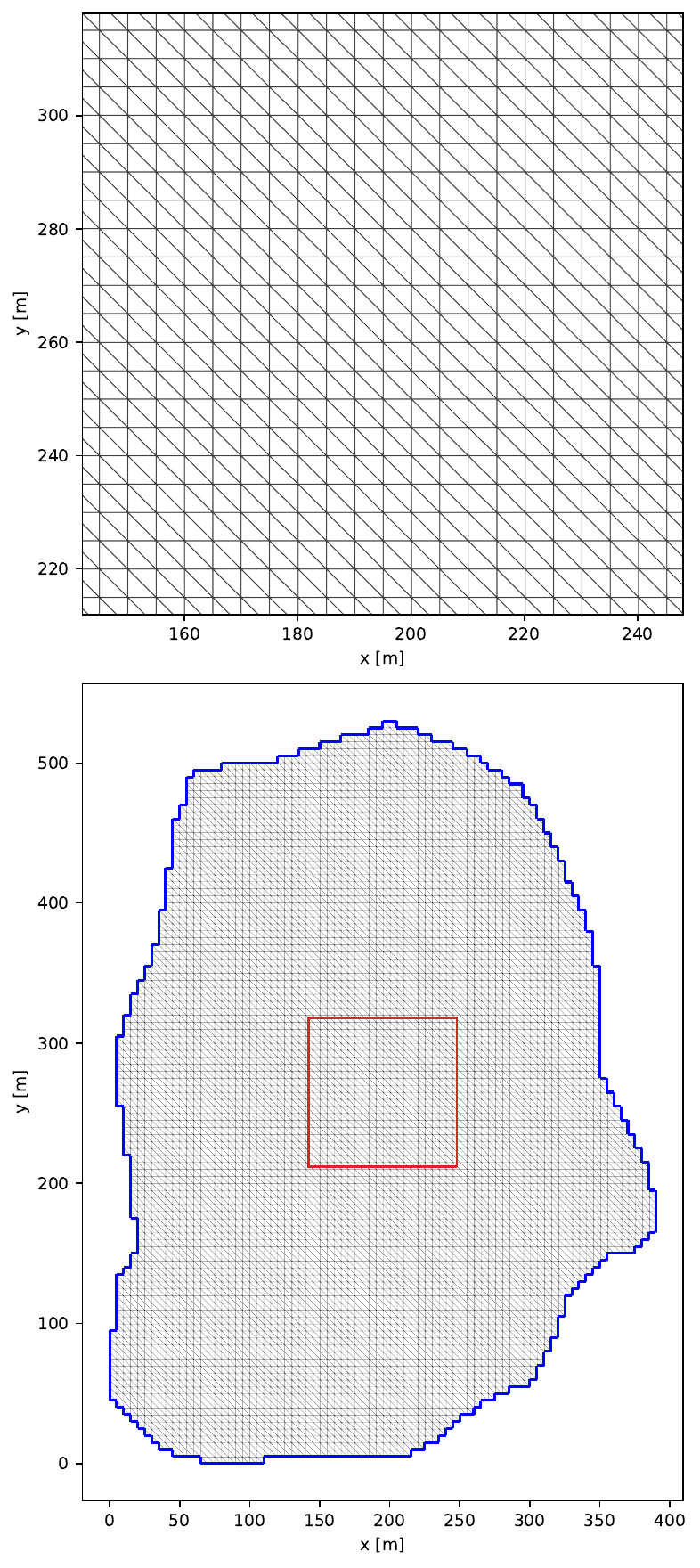}
\caption{Test 4: elevation of the Fernow catchment and the mesh, with a
zoom.}\label{fig:t4setup}\label{fig:ex2terrain}
\end{figure}

The real-terrain test considers a catchment in the Fernow Experimental
Forest, West Virginia, USA, which has been extensively investigated in
previous hydrological studies
\citep{JensenMcGuirePrince2017,BotterEtAl2021}. The mesh is the
triangulation of a depitted and smoothed 5\,m elevation raster, with
\pv{t4.nv} vertices and \pv{t4.nc} cells. The preprocessing of the
original DEM is not investigated here. The resulting terrain mesh is
used as input to assess the numerical properties of the proposed flux
reconstruction on real topographic data.
We assume free outflow on the whole boundary. See
Fig.~\ref{fig:t4setup} for the topography elevation and mesh. The
slope has mean \pv{t4.slopemean} and its maximum value is up to
\pv{t4.slopemax}, with $\cos\AngleSymbol\ge\pv{t4.cosmin}$.
The slope is resolved on most pairs of cells, while the curvature is
not resolved anywhere by the 5\,m raster.
The geometry is built from the heights only, and approximated in the
interior with the \pv{t4.rings}-ring fit. Boundary vertices use the
same 2-ring fit. The intra-cell term is small but not zero, with
$\max_K|d_K|/\max|F_{K,E}|={}$\pv{t4.dk}.
The rain event lasts \pv{t4.train}\,min, \pv{t4.tend}\,min are
simulated, and $\Dt_{\max}=\pv{t4.dtmax}$\,s.

\begin{table}\centering
  \caption{Test 4: traps / trapped rain at steady state on meshes of the same
    Fernow data (intrinsic geometry, free outflow on the whole
    boundary). The last row
    (smoother, depitted) is the mesh of Test~4.}\label{tab:t4census}
  \begin{tabular}{lcccc}
    \toprule
    data & {\FVS} & {\AVG} & {\LNs} & {\LNw} \\
    \midrule
    original & 12 / 98.5\% & 16 / 98.6\% & 9 / 98.3\% & 4 / 94.1\% \\
    smoothed, not depitted & 4 / 98.6\% & 5 / 94.5\% & 3 / 94.5\% & 3 / 94.5\% \\
    depitted & 13 / 90.4\% & 9 / 89.9\% & 2 / 89.9\% & 4 / 89.8\% \\
    smoothed, depitted & 1 / 91.6\% & 0 / 0.0\% & 0 / 0.0\% & 0 / 0.0\% \\
    \bottomrule
  \end{tabular}
\end{table}
\begin{figure}\centering
  \includegraphics[width=\textwidth]{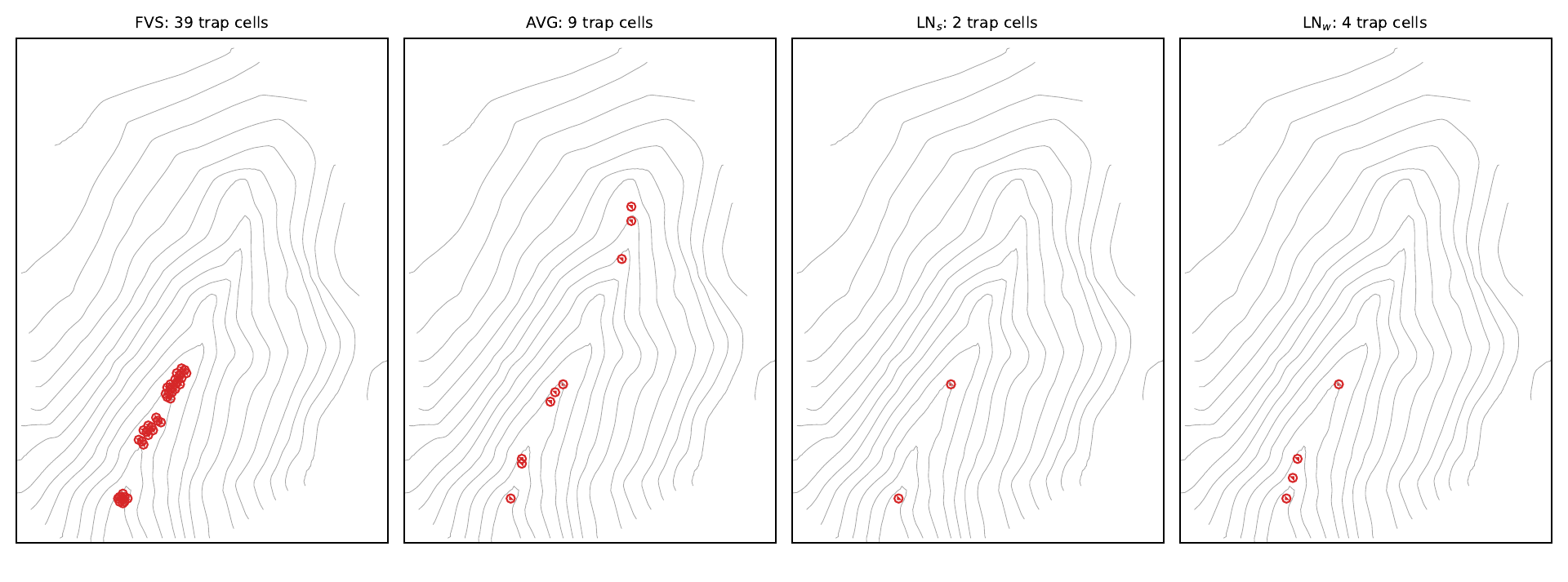}
  \caption{Test 4: trap cells (red, circled) of the four edge velocities on the
    \pv{t4c.fig.data} mesh. In gray, the contour lines of
    the heights.}\label{fig:t4traps}
\end{figure}

\paragraph{Mesh variants of the same data.} The mesh of Test~4 is trap
free for the three single-valued edge velocities, but this is a
property of that mesh, not of the data. Table~\ref{tab:t4census}
repeats the steady trap census on three other meshes of the same
catchment: the original data, the smoothed but not depitted data, the
depitted mesh without smoothing. Not all of them satisfy the
assumption of a resolved input (Section~\ref{sec:nodalgeom}), which is
the usual situation of an application.  \AVG\ has more traps than both
\LNs\ and \LNw\ on these variants.
Every trap of the three single-valued edge velocities is a single
cell, as in Remark~\ref{rem:traps}, whereas the traps of \FVS\ are
groups of cells that exchange water across counterflow edges. The
correction helps, in the number of traps, but is not sufficient: in
this case, the traps localize in the central portion of the channel,
where water converges, and apart from the mesh of Test~4,
approximately 90\% or more of the rainfall is trapped for all four
edge-velocity reconstructions considered. Figure~\ref{fig:t4traps}
shows the trap cells on the depitted input mesh. All but one of the
trap cells associated with the corrected fields are also traps for
\AVG. Removing them requires conditioning the data at the mesh scale,
not another edge velocity.

\paragraph{Intrinsic versus projected.} On terrain data too the
computed ratio follows the storage-weighted mean of $\kVol$, to
\pv{t4.proj.agree}, below the area mean, since storage concentrates on
the gentle valley floors (Table~\ref{tab:intproj}; \FVS, which traps,
has no steady state).

\begin{figure}\centering
\includegraphics[width=0.6\textwidth]{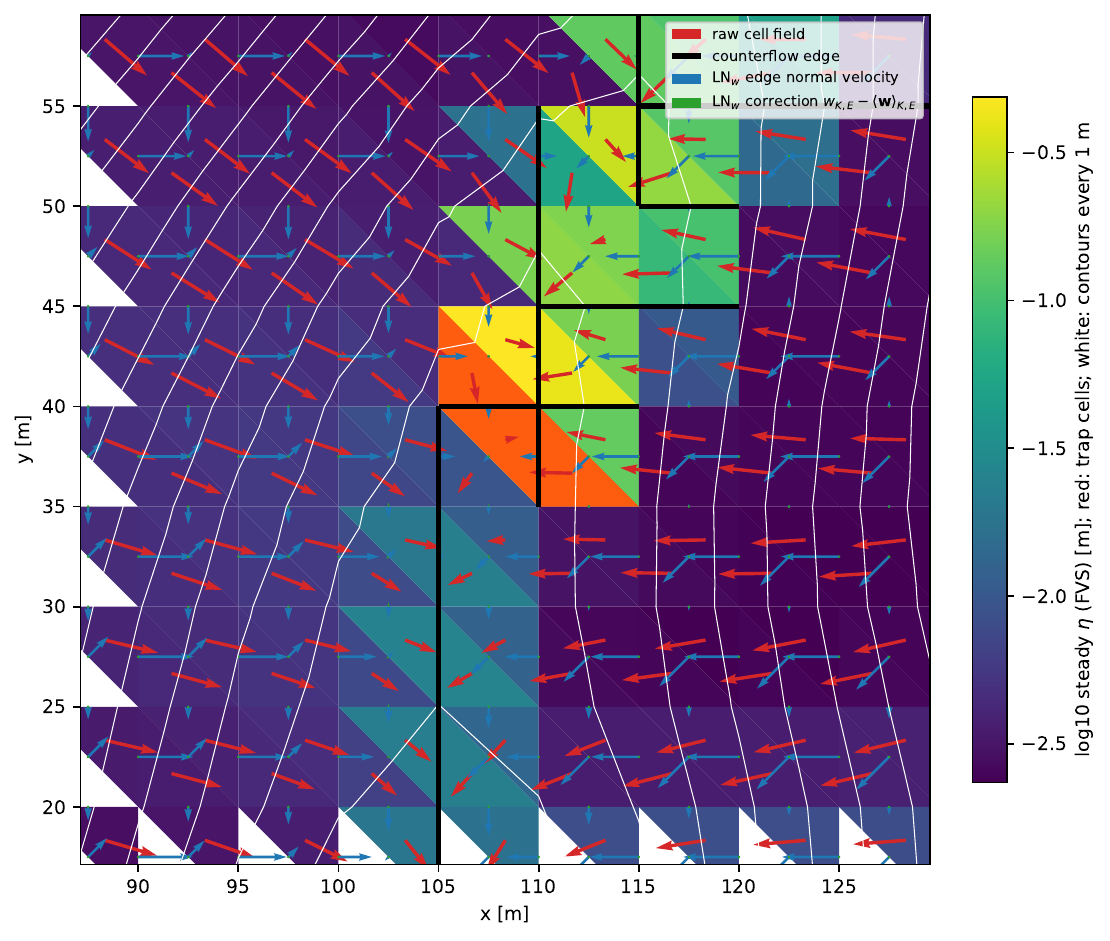}
\caption{Test 4: the \FVS\ trap on the main channel (red cells). {Red arrows:
raw cell field; black: counterflow edges; blue: \LNw\ edge normal
velocities; green: the \LNw\ correction.}}\label{fig:t4zoom}\label{fig:ex2trap}
\end{figure}

\begin{figure}\centering
\includegraphics[width=0.75\textwidth]{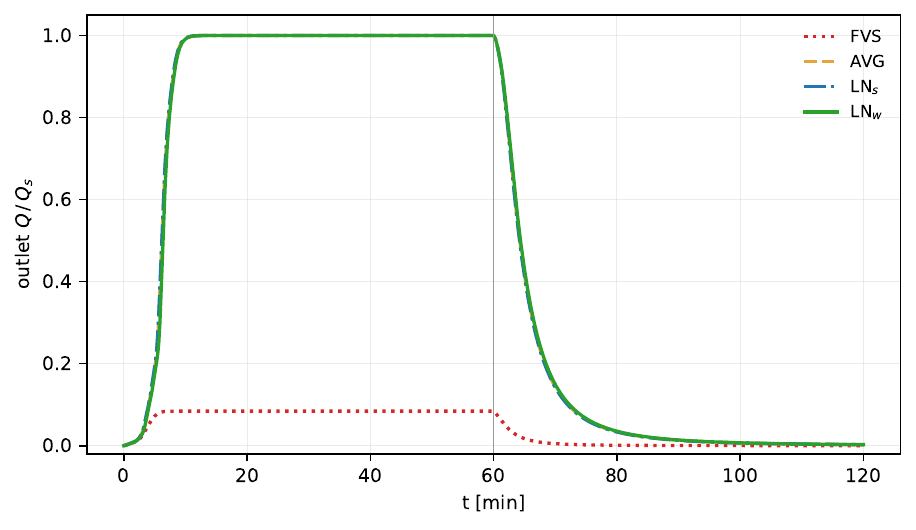}
\caption{Test 4: outlet hydrographs $Q/Q_s$.}\label{fig:t4hyd}\label{fig:ex2hyd}
\end{figure}

\begin{figure}\centering
\includegraphics[width=0.8\textwidth]{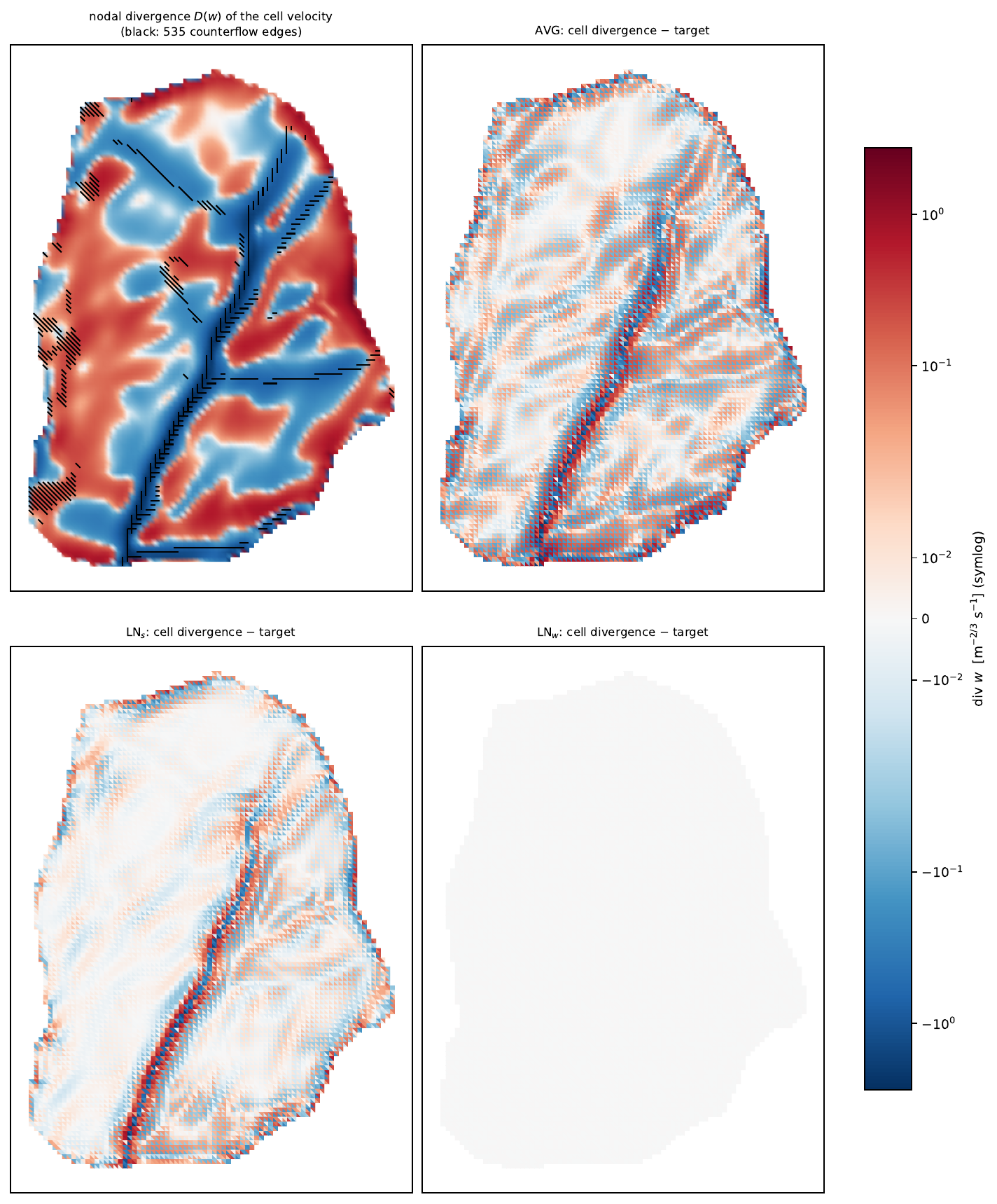}
\caption{Test 4: nodal divergence $D(\AdvVel)$ of the raw field (red:
  divergent, blue: convergent) with the counterflow edges (black), and
  the cell-divergence defects of \AVG, \LNs, \LNw, on one symmetric
  logarithmic scale.}\label{fig:t4div}
\end{figure}

\paragraph{Results.}
On real data the raw kinematic field jumps more across the edges than
on the analytical surfaces (relative normal jump \pv{t4.div.jump}),
and the \pv{t4.div.cf} counterflow edges (\pv{t4.div.cffrac} of the
edge length) lie along the channels (see Fig.~\ref{fig:t4div}). With
flux-vector splitting, one closed group of \pv{t4.fvs.trapcells} cells
on the main channel, \pv{t4.fvs.trap.dz}\,m above the lowest point of
the catchment, keeps \pv{t4.fvs.trapped} of the rain
(Fig.~\ref{fig:t4zoom}): the pond grows without bound
(\pv{t4.fvs.maxh}\,m at the end of the event), the time step collapses
with it, and the outlet receives only \pv{t4.fvs.peak} of the rain
(Fig.~\ref{fig:t4hyd}).
The same mesh is trap free for \AVG, \LNs\ and \LNw.
The divergence defects repeat those of Test~3 (see
Fig.~\ref{fig:t4div}): by construction \LNw\ matches the target of
Corollary~\ref{cor:divform} to machine precision, whereas the plain
average misses it by \pv{t4.div.avg.defect}, as on any unstructured mesh.
Although structured in projection, the lifted raster does not retain
the superconvergence of the average (Remark~\ref{rem:raster}): its
defect is \pv{t4.div.avg.overstru} times that of the structured mesh
of Test~2 (quadratic analytic surface) with the same spacing
(\pv{t2.stru.avg.defect.coarse}) and close to that of the unstructured
one (\pv{t2.tria.avg.defect.coarse}). We attribute it to the
curvature of the data, which varies at the scale of the cells, and
therefore, for the imposed divergence, the whole mesh behaves as an
irregular one.

\begin{figure}\centering
\includegraphics[width=\textwidth]{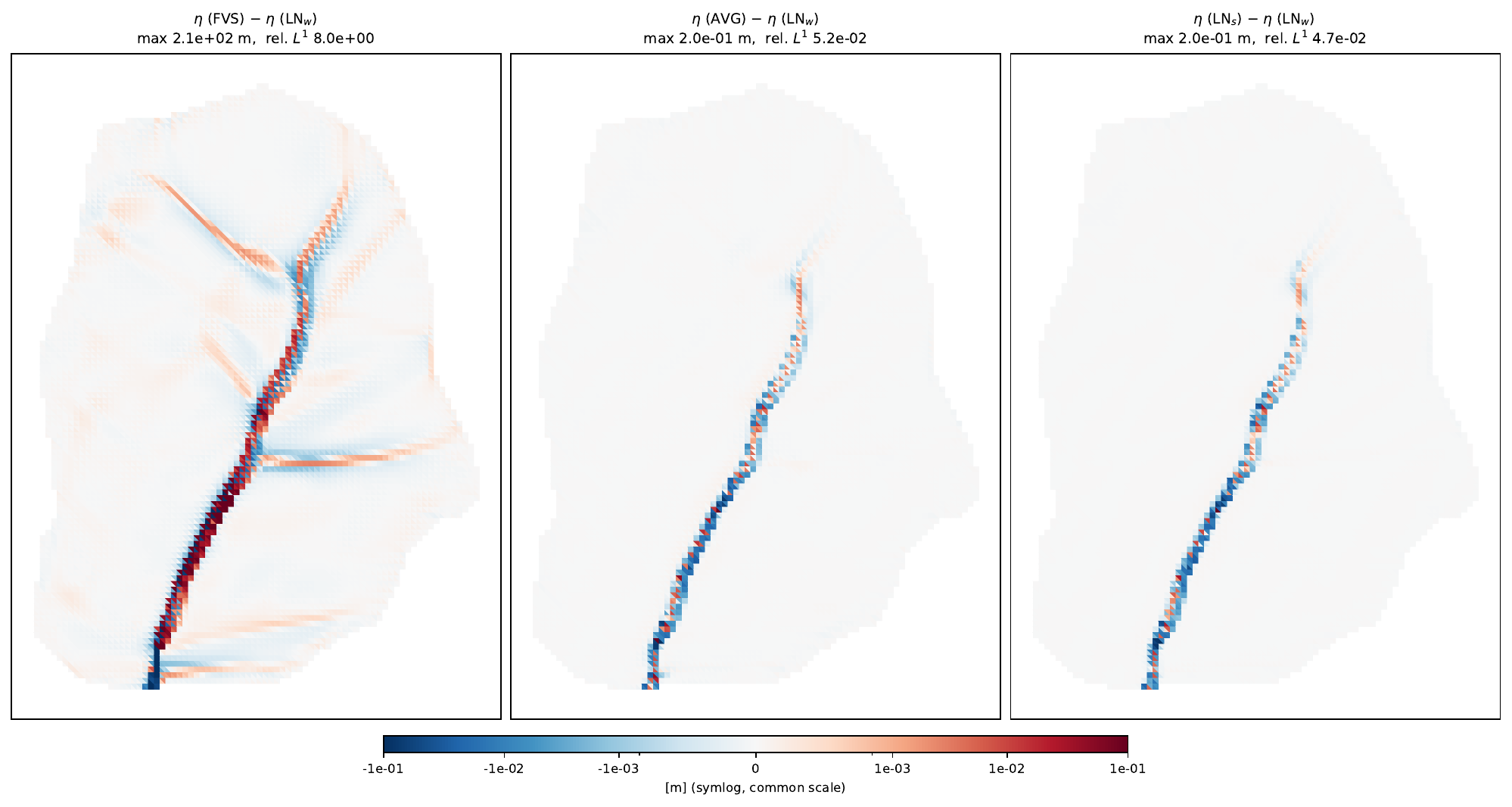}
\caption{Test 4: depth differences $\Depth_X-\Depth_{\LNw}$, $X\in\{\FVS,\AVG,\LNs\}$, at the
end of the rain, on one common symmetric logarithmic color
scale.}\label{fig:t4diff}\label{fig:ex2diff}
\end{figure}

On the outlet response the three single-valued edge velocities are
close (Fig.~\ref{fig:t4hyd}): the mean response times are
\pv{t4.avg.tmean}, \pv{t4.lns.tmean} and \pv{t4.lnw.tmean}\,s for
\AVG, \LNs\ and \LNw\ (a \pv{t4.tmean.spread} range), the steady depth
of \AVG\ and \LNs\ differs from that of \LNw\ by \pv{t4.avg.vslnw} and
\pv{t4.lns.vslnw}, the hydrographs overlap, and the differences sit in
the channel cells (Fig.~\ref{fig:t4diff}).

\subsection{Discussion of the tests: \AVG\ and \LNw\ on terrain
  meshes}
\label{sec:testdiscussion}

On the meshes of the four tests \AVG\ is trap free, like \LNw\, and at
first order its depths and hydrographs are close to those of \LNw: the
steady depths differ by \pv{t3.avg.vslnw} (Test~3) and
\pv{t4.avg.vslnw} (Test~4) and the mean response times by at most
\pv{t4.tmean.spread}, the largest local difference being the maximum
depth of Test~4 in the channel cells (\pv{t4.avg.maxh} against
\pv{t4.lnw.maxh}\,m). The transport error of the scheme is at least
one order of magnitude larger than every slope-induced error of the
single-valued edge velocities (see Test~2), and against the exact
steady depth of Tests~2 and~3 \LNw\ is not more accurate than
\AVG. The two differ in what the edge field satisfies. \LNw\ is
defined by star problems that are solvable on every mesh, imposes the
divergence of the data to roundoff, and its slope-induced error, like
that of \LNs, is \pv{t2.tria.avgoverlnw} times smaller on the
unstructured meshes of Test~2. 
\AVG\ satisfies the imposed divergence only asymptotically, on meshes
that are structured on the surface, where its defect decreases at
first order. However, on unstructured meshes its defect does not
decrease.

Surface triangulations of elevation data are not structured on the
surface, even when they are structured in projection.  This effect is
not caused by the irregularity of the lifted raster but by the
curvature, which for example the 5\,m raster does not resolve. No
refinement is available at the fixed spacing of the data, and finer
data are not always available, nor wanted.

On the terrain meshes considered here, the plain
average produces an edge field whose divergence
misses the imposed target, i.e.\ the cell average
of the nodal divergence of the data, by an amount
of order one. The correction, instead, matches
this target in every cell, independently of the
resolution of the data. This is a property of
the cell balance, not of the accuracy of the
divergence.

Since the kinematic field depends only on the terrain, the star
problems are solved once, before any simulation (on Fernow, the four
edge velocities and all exact steady states take \pv{t4.tsteady}\,s on
a personal laptop): the correction keeps the reason for choosing the
kinematic wave, its low cost. At first order this consistency does not
change depths or hydrographs, because the transport error hides it. A
higher-order transport discretization would be needed to determine
whether the edge-velocity error becomes dominant as the transport
error is reduced: \LNw\ removes the order-one divergence defect, but
its slope-induced error is still of first order (observed order
\pv{t2.tria.lnw.slopeerr.rate}, against \pv{t2.tria.avg.slopeerr.rate}
for \AVG, Test~2), with a constant \pv{t2.tria.avgoverlnw} times
smaller on the unstructured meshes (\LNs: similar). 
If so, a higher-order transport scheme on terrain meshes needs, beyond
the correction, higher-order edge velocities, for instance the
correction of a recovered field.

Thus, at first order the choice between \AVG\ or \LNw\ is not on the
accuracy in depths and hydrographs, but the corrected field gives two
guarantees that the plain average does not: every star problem is
solvable and every cell carries the imposed divergence, on every mesh,
and a single-cell trap can only form where the data converge. In
addition, \LNw\ reduces the slope-induced error on irregular meshes by
a factor \pv{t2.tria.avgoverlnw} (see Test~2), and this is the error
that would remain once the transport error is reduced.

\section{Conclusions and outlook}\label{sec:conclusions}\label{sec:discussion}

In a kinematic wave model, the edge velocity computed from the
elevation data determines the flow paths, and no term of the model can
redirect a flow that the edge velocity misroutes. We formulated the
kinematic wave on the bed surface, in the intrinsic framework of
\citet{BachiniPutti2020}, discretized it by first-order finite volumes
on a surface triangulation of elevation data, and built a
single-valued edge velocity by a conservative correction of
Larson-Niklasson type.  Three findings emerge.

First, the intrinsic and the projected kinematic waves are not
interchangeable. Under vertical rain they share the steady discharge
field, on any graph surface, but they differ in storage and timing,
through the depth factor $\kEta=\cos^{-3/10}\AngleSymbol$ and the
storage factor $\kVol=\cos^{-13/10}\AngleSymbol$. On a plane the
difference is a uniform stretch of time (Remark~\ref{rem:timescale}),
while on curved terrain there is no single stretch, because the
response time weights $\kVol$ by where the water is stored, not by
area. The resulting gap in timing is a model difference, not a
discretization error: it does not vanish under refinement, it exceeds
the discretization error on all but the mildest slopes, and it grows
with the slope in the cases considered; the applicability of the
kinematic approximation requires a separate assessment of its
neglected terms.

Second, the consistency of the edge velocity with the data has a
precise discrete meaning. The star problems of the correction are
solvable only with a compatible load, the weak discrete divergence of
the field, assembled from its normal jumps across the edges
(Lemma~\ref{lem:compat}). The corrected field then carries the
divergence of the data in every cell, at any resolution
(Corollary~\ref{cor:divform}). The plain average does so only on
meshes that are structured on the surface and resolve the curvature of
the terrain, and surface triangulations of elevation data,
triangulated rasters included, are not such meshes
(Section~\ref{sec:testdiscussion}).

Third, robustness is a property of the transfer directed graph, not of
accuracy. The steady state of the upwind scheme solves an M-matrix
system if and only if every cell has a directed path to an outlet
(Section~\ref{sec:traps}), so spurious sinks are closed classes of the
graph and can be counted before any simulation.  This is a property of
the mesh and is not due to the accuracy of the edge values.  The
correction can eliminate some numerical traps introduced by the
edge-velocity reconstruction, but does not guarantee trap-free routing
when the terrain data are insufficiently resolved or conditioned.

Summarizing, our analysis shows that, for a kinematic wave on
terrain data, it is recommended to:
\begin{itemize}
\item provide a resolved input surface with its pits filled, since
  without a pressure term nothing redirects the flow along a crease or
  out of a pit;
\item count the traps of the chosen edge velocity in the chosen
  geometry (mesh)
  before simulating;
\item use a single-valued, compatibly corrected edge velocity (\LNs\ or
  \LNw): it is computed once, since the slope and the kinematic field depend
  only on the terrain. %
\end{itemize}

Three questions remain open. The first is how to bring rough terrain
data, such as rasters, meshes with stripes or triangulated V-creases,
to a resolved input without distorting the timing, possibly building
on terrain partitioning approaches that preserve topographic features
\citep{MorettiOrlandini2023}. The second is how to find, if possible,
a sufficient condition for the absence of traps of the compatible
corrections, in terms of the mesh and of the piecewise-linear
structure of the surface \citep{Banchoff1967,Edelsbrunner2003}. The
third is which of the two kinematic models is closer to the physical
flow, which requires a higher-fidelity reference: the diffusive wave
and the shallow water equations in the intrinsic framework of
\cite{BachiniPutti2020}, with a two-dimensional, geometric extension
of the classical applicability criteria
\citep{WoolhiserLiggett1967,Ponce1978}.
The terrain part of the driver is common to the three models, and its
consistency with the data, at the resolution of the data, is what this
first step establishes.

\section*{Data availability}
The code and input data used in this study are available from the
authors upon reasonable request.

\section*{CRediT authorship contribution statement}
Elena Bachini: Conceptualization, Methodology, Formal analysis,
Software, Validation, Writing -- original draft.

Matteo Camporese: Conceptualization, Methodology, Writing -- review \&
editing.

\section*{Declaration of generative AI and
AI-assisted technologies in the manuscript
preparation process}

The authors used AI-assisted tools in a limited capacity to support
code development (Claude-Anthropic) and the writing process
(Claude-Anthropic and ChatGPT-OpenAI). All mathematical content,
numerical results, and final text were reviewed and validated by the
authors, who take full responsibility for the manuscript.

\section*{Acknowledgments and funding}
E. Bachini received funding from the PRODIGI (``Processi e modelli
digitali innovativi per il controllo, la sostenibilità e la
valorizzazione delle produzioni dell’industria agroalimentare in
Veneto'') project, CUP C23C25001260002, funded by Regione del Veneto,
and from the project ``Metodi numerici integrati per la simulazione e
la prevenzione del dissesto idrogeologico'', funded by INdAM-GNCS
2026. E. Bachini is a member of the Gruppo Nazionale Calcolo Scientifico
- Istituto Nazionale di Alta Matematica (GNCS-INdAM).

\appendix
\section{The correction step by step}\label{app:algorithm}
This appendix collects the construction of
Sections~\ref{sec:lnconstruction}--\ref{sec:boundary}. The input is
the surface triangulation with its cell and edge frames, and a
piecewise-constant tangent field $\cellfield$: the kinematic field
$\AdvVel$ for \LNw, the slope $\GradSurf\zbottom^h$ for \LNs. Steps~3
and~4 are steps (1)--(3) of \citet{LarsonNiklasson2004}, while Step~2
replaces the given source of the classical setting by the compatible
load.
\begin{enumerate}[label=\textit{Step \arabic*.},leftmargin=*,align=left]
\item \emph{Edge scalars.} For every cell $K$ and every edge
  $E\subset\partial K$, transfer $\cellfield_K$ to the frame of $E$
  \eqref{eq:transfer} and compute $g_{K,E}$ and $d_K$ \eqref{eq:gKE}. On every
  interior edge, compute the average flux $\bar g_{K,E}$ \eqref{eq:avgflux} and
  the jump $J_E$ \eqref{eq:jump}, and on a boundary edge set $\bar g_{K,E}=g_{K,E}$.
\item \emph{Compatible load.} For every vertex $i$, interior or boundary,
  compute the lumped mass $m_i=\sum_{K\ni i}|K|/3$ and the load
  $b_i=\tfrac12\sum_{E\ \text{interior},\,i\in E}J_E$ \eqref{eq:compload}, i.e.\
  $D_i=-b_i/m_i$ \eqref{eq:Dnodal}.
\item \emph{Star problems.} For every vertex $i$, assemble on the cells of
  $\Star{i}$ the problem $\mathsf L_iV_i=\rho_i$, where
  the graph Laplacian is computed as
  \begin{equation*}
    (\mathsf L_iV)_K=\sum_{E\subset\partial K\ \text{interior},\,i\in E}
    \tfrac{\edgeLength[E]}{2}(V_K-V_{K'})
  \end{equation*}
  and the right-hand side  \eqref{eq:starrhs} is
  \begin{equation*}
  \rho_{i,K}=-R_{i,K}=\tfrac14\sum_{E\subset\partial K\ \text{interior},\,i\in E}J_E
  -\tfrac{|K|}{3}\,b_i/m_i \ .
\end{equation*}
 The right-hand side sums to zero over the star, so the
  system is solvable (Lemma~\ref{lem:compat}). Its solution is unique up to a
  constant, which does not affect the fluxes, and fixing one unknown is
  exact.
\item \emph{Corrected flux.} On every interior edge $E$ shared by $K$
  and $K'$, the corrected flux is given by \eqref{eq:Fcorr}:
  \begin{equation*}
  F_{K,E}=\bar g_{K,E}+\sum_{i\in
    E}\tfrac{\edgeLength[E]}{2}(V_{i,K}-V_{i,K'}),
\end{equation*}
and on a boundary edge we set $F_{K,E}=g_{K,E}$. The edge velocity is
  $w_{K,E}=F_{K,E}(\AdvVel)/\edgeLength[E]$ for \LNw\ \eqref{eq:kln}, and the
  Manning law \eqref{eq:kslope} with $s^\nu_E=F_{K,E}(\GradSurf\zbottom^h)/\edgeLength[E]$
  for \LNs.
\end{enumerate}
The cell balance \eqref{eq:balance} provides a check to roundoff. Each star
problem has as many unknowns as cells around the vertex, and the star
problems are independent of each other. Since the kinematic field depends
only on the terrain, the four steps are carried out once, before any
simulation.

\section{Exact event hydrograph on the sloping plane}\label{app:plane}
In this section, we recall the exact outlet hydrograph of Test~1. On
the plane $\BSM=S(L-x)$, with $S=\tan\AngleSymbol$, the intrinsic and
the projected kinematic wave equations can be written in the same way
per unit contour width:
\begin{equation*}
  \Der_t\Depth+\alpha\,\Der_\xi(\Depth^{5/3})=r_*
\end{equation*}
with a general coordinate $\xi\in[0,\Lambda]$. For the projected model
we have $\Depth=\DepthP$ and
$(r_*,\alpha,\Lambda)=(r,\tan^{1/2}\AngleSymbol/\Manning,L)$, while
$(r\cos\AngleSymbol,\sin^{1/2}\AngleSymbol/\Manning,L/\cos\AngleSymbol)$
in the intrinsic one.
Starting from a dry bed, the outlet discharge is
$q=\alpha(r_*t)^{5/3}$ until
\begin{equation*}
  t_e=\Big(\frac{\Lambda}{\alpha\,r_*^{2/3}}\Big)^{3/5},
\end{equation*}
 and $q=rL$ afterwards
\citep{OvertonMeadows1976}, so that the
steady discharge is the same in the two models. After the end of the rain at $t_r\ge t_e$,
the characteristic leaving the equilibrium profile at $\xi_0$ carries
$q=r_*\xi_0$ and reaches the outlet when
\begin{equation}\label{eq:recession}
t-t_r=\frac{\Lambda-q/r_*}{\tfrac53\,\alpha^{3/5}q^{2/5}},
\end{equation}
which defines the recession $q(t)$ implicitly (the same relation holds
for $t_r<t_e$, after the plateau $q=\alpha(r_*t_r)^{5/3}$).

\section{Exact steady solution on the convergent valley}\label{app:l4}
In this section, we compute the steady solution on the surface of
Test~2.  We consider a domain $[0,L]\times[0,W]$ with height function
$\BSM=S(L-x)+a(y-y_c)^2$, where $y_c=W/2$ is the central line of the
valley. For the projected model, the flow direction is
$-\Grad\BSM=(S,\,-2a(y-y_c))$, so flow lines satisfy
$dy/dx=-2a(y-y_c)/S$, i.e.
\begin{equation*}
y-y_c=(y_1-y_c)\,e^{2a(x_1-x)/S}
\end{equation*}
through $(x_1,y_1)$. All flow lines start either at $x=0$, where no
water enters, or on the side walls, where the flow points inward, so
that the walls are divides. The strip between the flow lines through
$y_1$ and $y_1+dy_1$ has width $e^{2a(x_1-x)/S}dy_1$ at abscissa $x$,
so the steady unit discharge in $x$ is
\begin{equation*}
q_x(x_1,y_1)=r\int_{x_s}^{x_1}e^{2a(x_1-x)/S}\,dx=\frac{rS}{2a}\Big(e^{2a(x_1-x_s)/S}-1\Big),
\quad x_s=\max\Big(0,\;x_1-\frac{S}{2a}\ln\frac{W}{2|y_1-y_c|}\Big),
\end{equation*}
with $r$ the rain intensity per unit horizontal area of \eqref{eq:kw}
and with $x_s$ the abscissa where the flow line leaves the wall. With
$q_x=\DepthP^{5/3}S/(\Manning|\Grad\BSM|^{1/2})$, the exact steady depth is
$\DepthP=\big(\Manning\,q_x\,|\Grad\BSM|^{1/2}/S\big)^{3/5}$, and by
\eqref{eq:depthratio} the intrinsic depth is
$\Depth=\kEta\DepthP$. The discharge vector,
$\QvectP=q_x\,\big(1,\,-2a(y-y_c)/S\big)$, is parallel to $-\Grad\BSM$, with
magnitude $|\QvectP|=q_x|\Grad\BSM|/S$, and $\ABSsurf{\Qvect}=|\QvectP|$ at the
corresponding point of $\SurfDomain$ (Section~\ref{sec:equiv}).


\bibliographystyle{abbrvnat}
\bibliography{strings_paper,biblio_paper}

\end{document}